\documentclass[11pt]{elsarticle}
\usepackage{comment,listings,framed,graphicx,url,xcolor}
\usepackage{amsmath,amssymb,amsthm}
\usepackage{fullpage} % changes the margin
\definecolor{mild}{rgb}{1,0.98,0.8}
\usepackage{definitions}
\newcommand{\ba}{\begin{eqnarray}}
\newcommand{\ea}{\end{eqnarray}}
\newcommand{\bsub}{\begin{subequations}}
\newcommand{\esub}{\end{subequations}}
\newcommand{\bas}{\begin{eqnarray*}}
\newcommand{\eas}{\end{eqnarray*}}
\newcommand{\ourcite}[2]{\cite{#1}~(#2)}
\newcommand{\Root}{Root}

\newcommand{\id}{\mathrm{id}}
\newcommand{\uz}{u^{init}}
\newcommand{\vz}{v^{init}}
\newcommand{\wz}{w^{init}}
\newcommand{\adm}{representable}
\newcommand{\htop}{H_{t}}
\newcommand{\hbot}{H_{b}}
\newcommand{\hleft}{H_{l}}
\newcommand{\hright}{H_{r}}
\newcommand{\Nave}{N_{ave}}
\newcommand{\hH}{\mathcal{H}}
\newcommand{\sS}{\mathcal{S}}
\newcommand{\dD}{\mathcal{D}}
\newcommand{\hup}{\gamma_{r}}
\newcommand{\hdown}{\gamma_{l}}
\newcommand{\fl}{\hdown}
\newcommand{\fr}{\hup}
\newcommand{\flk}{\gamma_{l,k}}
\newcommand{\frk}{\gamma_{r,k}}
\newcommand{\bh}{b^h}
\newcommand{\bkh}{b^{h_k}}
\newcommand{\bhstar}{\bh_*}
\newcommand{\bheps}{\bh_{\eps}}
\newcommand{\bhz}{\bh_0}
\newcommand{\bhinf}{b^{\infty}}
\newcommand{\bhh}{\bh_h}

\newcommand{\tT}{\mathcal{T}}

\newcommand{\cg}{\mathcal{C}}
\newcommand{\calbe}{\cg_{\alpha,\beta}}

\newcommand{\ver}{vert}
\newcommand{\bd}{bd}
\newcommand{\uU}{{\mathcal{U}}}
\newcommand{\pP}{\Pi}
\newcommand{\ppreisach}{\Pi^{K,0}}
\newcommand{\ppreisache}{\Pi^{K,\eps}}
\newcommand{\ppreisachs}{\Pi^{K,*}}
\newcommand{\pknon}{\Pi^{K,1}}
\newcommand{\pklin}{\Pi^{K}_{\infty}}
\newcommand{\pgen}{\Pi^{\gamma}}
\newcommand{\Q}{\mathbb{Q}}
\newcommand{\vbar}{\overline{V}}
\newcommand{\ubar}{\overline{u}}

\newcommand{\PL}{\mathcal{L}}

\newcommand{\Hi}{V}
\newcommand{\vs}{v^*}
\newcommand{\cab}[1]{\cg(a,b;#1)}
\newcommand{\ccab}{\cg(a,b)}
\newcommand{\rab}[1]{R(a,b;#1)}
\newcommand{\resab}[3]{R(#1,#2;\;#3)}
\newcommand{\resk}{\mathcal{R}_k}
\newcommand{\resg}{\mathcal{R}_{\gamma}}

\newcommand{\ggres}{{\mathcal{G}}}
\newcommand{\vit}{\tilde{v}}
\newcommand{\Dx}{\tfrac{\partial}{\partial x}}
\newcommand{\Dt}{\tfrac{\partial}{\partial t}}
\newcommand{\dtt}{\tfrac{d}{dt}}

\newcommand{\Dom}{\operatorname{Dom}}
\newcommand{\Rg}{\operatorname{Rg}}
\newcommand{\sgn}{\operatorname{sgn}}
\newcommand\myskip[1]{}
\newcommand{\mpcomment}[1]{{#1}}

\newtheorem{remark}{Remark}
\newtheorem{lemma}{Lemma}
\newtheorem{proposition}{Proposition}

\newtheorem{definition}{Definition}
\newcommand{\linear}{linear play}
\newcommand{\klinear}{$K$-linear play}
\newcommand{\knonlinear}{$K$-nonlinear play}
\newcommand{\kpreisach}{$K$-Preisach}
\newcommand{\kgeneral}{$K$-generalized play}
\newcommand{\general}{generalized play}
\newcommand{\nonlinear}{nonlinear play}
\newcommand{\ghysteron}{{}trapezoidal hysteron{}}
\newcommand{\uhysteron}{{}unit hysteron{}}
\begin{document}
%%%%%%%%%%%%%%%%%%%%%%%%%%%%%%%%%%%%%%%%%%%%%%%%%
\begin{frontmatter}

\author{Malgorzata Peszynska and Ralph E. Showalter}
\address{Department of Mathematics, Oregon State University, Corvallis, OR, 97331, USA}
\ead{mpesz@math.oregonstate.edu, show@math.oregonstate.edu}


\title{Approximation of Hysteresis Functional\tnoteref{t1}}
\tnotetext[t1]{This work was partially supported by the National Science Foundation DMS-1912938 and DMS-1522734, and by the NSF IRD plan 2019-21 for M.~Peszynska while serving at the National Science Foundation. Any opinion, findings, and conclusions or recommendations expressed in this material are those of the authors and do not necessarily reflect the views of the National Science Foundation.
}
%%%
\begin{abstract}
We develop a practical discrete  model of hysteresis based on nonlinear play and generalized play, for use in first-order conservation laws with applications to adsorption-desorption hysteresis models. The model is easy to calibrate from sparse data, and offers rich secondary curves. We compare it with discrete regularized Preisach models. We also prove well-posedness and numerical stability of the class of hysteresis operators involving all those types, describe implementation and present numerical examples using experimental data. 
\end{abstract}
%%%%%%%%%%%%%%
\end{frontmatter}
%%%%%%%%%%%%%%%%%%%%%%%%%% HEADER
\noindent
%%
%\textbf{PS-AHF paper} \hfill \hfill Date: \today
%\\
%%%%%%%%%%%%%%%%%%%%%%%%%%
%\medskip
%\hrule
%\medskip
%\tableofcontents
%%%%%%%%%%%%%%%%%%%%%%%%%
%\newpage
\section{Introduction}
%%%%%%%%%%%%%%%%%%%%%%%
In this paper we describe and analyze a new robust and fairly simple algorithm for approximation and calibration of hysteresis functionals $u \to w=\hH(u)$ which can be used in numerical schemes for PDEs arising in the applications. This paper extends the results in \cite{PeszShow20} to a broader class of hysteresis models. We explain how the model is calibrated, provide details on the solver, and compare the advantages and disadvantages of the different hysteresis constructions.  
Our work is motivated  by the applications to flow and transport in porous media, and specifically by the adsorption--desorption hysteresis \cite{AdsHyst02,SarkMon00L,KierlikMon01,MonsonLibby,MP16} which is significant and important in modeling of carbon sequestration \cite{JessenK,Czerw,Prusty2008,PIMA11,Assef19,Firozaabadi} and wood science and engineering \cite{Salin,Fredriksson}. We consider the PDE model
\ba
\label{eq:pde}
\Dt (a(u)+ \hH(u)) + A(u)=f,
\ea
in which $u$ is the unknown, $f$ is an external source, $A$ is a transport operator (advective and/or diffusive, generally nonlinear), and $a(\cdot)$ is a strongly monotone function. The problem is posed in the sense of distributions, in a functional space to be made precise below, and with some boundary and initial data.  

Hysteresis is a well known nonlinear phenomenon in which the output of a process depends not only on the independent variable,  but also on the history of the process, in a rate independent way. Hysteresis is well known to occur in electromagnetism \cite{Kadar88,HoffMeyer89,Visintin94}, plasticity \cite{Mielke10,SimoHughes,HoffGotz98}, phase transitions \cite{LittShow94,HoffKenm07,ShowLittHorn96}, multiphase flow in porous media \cite{Topp,Philip,Mualem,KilloughCarlson,BeliaevHass,Schweizer,CaoPop}, and many other applications \cite{HoffBotkin,HornShow93} including food processing and ecology \cite{food2002,Masud18,Beisner}. The hysteresis models have been well studied, and the models range from simple to complex, with the latter requiring detailed data; see, e.g., the monographs and reviews in \cite{KrasPokr89,Mayergoyz,Visintin94,BrokateSprekels96,Macki93}.  In particular, the ingenious well-known and well analyzed Preisach model considers 
a collection of (input, output) pairs from data  $\dD=(u(t),w(t))$ for $u \in \dD_u=\{u: u (u,w) \in \dD\}$, with $\dD_u$ dense in $ C([0,T])$, and records the hysteretic output $u\to w=\hH(u)$ in the so-called Preisach plane. This record is then used to build $w=\hH(u)$ as an integral over a continuum of parameters, for an arbitrary input $u$.
See \cite{VerdiVis85,VerdiVis89,HoffMeyer89,Krejci13,HSV88,Macki93}. 

However, experimental data $\dD$ for hysteresis is frequently sparse rather than dense \cite{JessenK,Salin,Prusty2008}; this limits the use of the Preisach model; in addition, its discrete form produces a very rough output $w(t)$. Our aim here is to approximate $\hH(u)$ with a practical tunable hysteresis model producing a piecewise smooth $w(t)$ when $\dD$ is only modest. An alternative is to ignore the hysteretic nature of $\hH(u)$, but this may lead to substantial modeling errors in predictive simulations of \eqref{eq:pde} \cite{Assef19,KilloughCarlson}.

 The data in $\dD$ includes, at the minimum, the boundary $H$ of the graph $(u,\hH(u))$; see Fig.~\ref{fig:teaser} for illustration. In particular, $H$ contains the ``left'' and ``right'' bounding curves $\fl(u)$ and $\fr(u)$ called {\em primary scanning curves}; here $\hup(u)\leq \hdown(u)$ are piecewise smooth monotone increasing functions. When the input $u(t)$ is increasing, the output $w(t) \in \hH(u(t))$ eventually reaches  the curve $\hup$, which it then follows upward.  Similarly, when $u(t)$ decreases, the output $w(t) \in \hH(u(t))$ eventually reaches and descends along the curve $\hdown$. When the input changes direction, $w(t)$ switches between $\fl$ and $\fr$ along the secondary scanning curves prescribed by the particular model. Since the models we consider are approximate, we usually obtain $H^* \approx H$ and $\hH^*(u) \approx \hH(u)$.   

%%%%%%%%%%%%%%%%%%%%%%%%%
\subsection*{\mpcomment{Example}}
%%%%%%%%%%%%%%%%%%%%%%%%%%

\begin{figure}[ht]
\centering
\includegraphics[height=28mm]{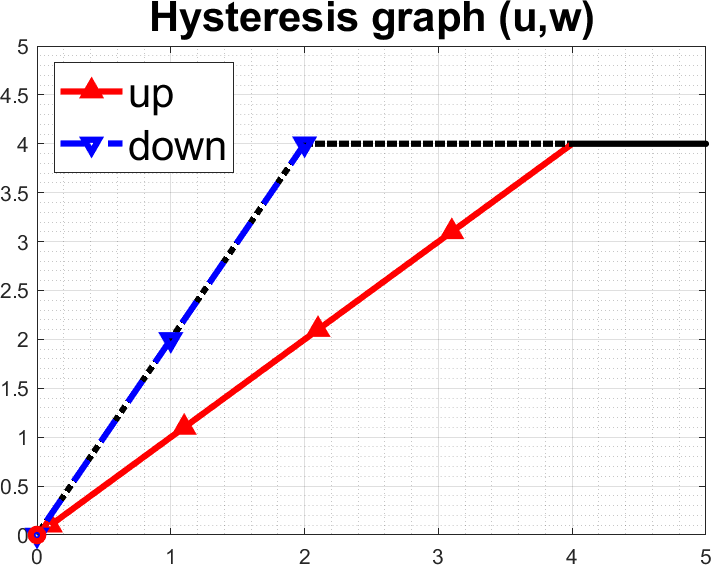}\;\;
\includegraphics[height=28mm]{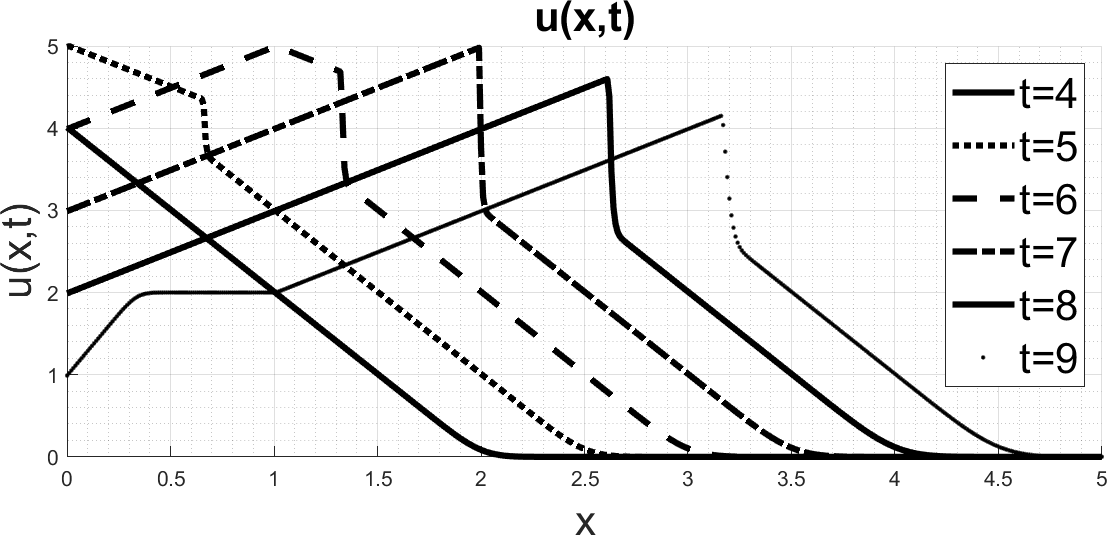}\;\;
\includegraphics[height=28mm]{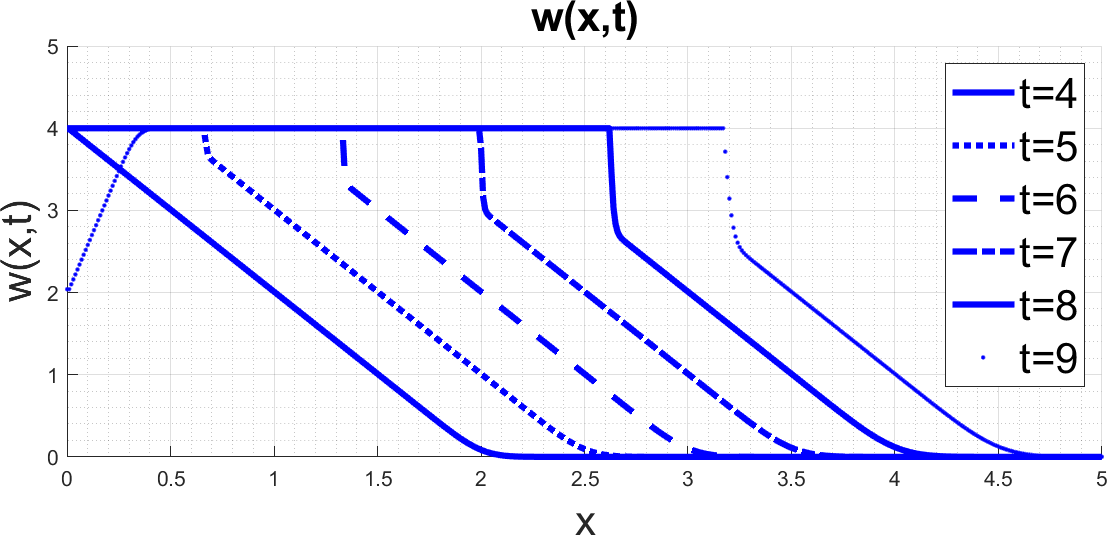}
\caption{Adsorption hysteresis example. Left: graph $(u,w)$.  
Middle: plot of $u(x,t)$. Right: plot of $w(x,t)$ at selected times as shown.\label{fig:example:PDE}}
\end{figure}

\mpcomment{
Consider the adsorption of a chemical of concentration $u(t)$ in the fluid at a point within a porous medium, and let $w(t)$ be the concentration of that chemical that is adsorbed onto or desorbed from the particles of the porous medium. Classical models assume these are related by a function $w = b(u)$ of {\em Langmuir} type. However, $u$ and $w$ are related more generally by a {\em hysteresis} relationship: they follow one path when they increase and another when they decrease. For an explicit example we assume the amount of solute $w$ adsorbed by the porous medium increases according to $w = u$ up to a maximum adsorbed concentration of $4$, but it desorbs from there only after $u$ has decreased to $2$ and thereafter is given by $w = 2u$. (See Figure \ref{fig:example:PDE}, Left.)
Such a relationship can be described with the truncation function $b(s) = s^+ -(s-4)^+$: it increases along the right scanning curve $w = \fr(u) = b(u)$ and decreases along the left scanning curve $w = \fl(u) = b(2u)$. We assume further that $w = b(v)$ is constant between these curves, {\em i.e.}, when $u < v < 2u$. For example, if the fluid concentration $u(t)$ at a point increases from $0$ to $5$, the amount adsorbed onto the medium at that point is $w(t) = b(u(t))$. As the concentration decreases from $5$ down to $0$, the adsorbed amount decreases according to $w(t) = b(2 u(t))$. 
The adsorbed concentration $w(u) = \hH(u)$ and the total concentration $m(u) = u + w(u)$ are given by  
} %% END MPCOMMENT

\mpcomment{
\begin{equation} \label{examplefcns}
w(u) = \begin{cases} 
u,\ & 0 \le u \le 4, \text{ increasing}, \\
4,\ & 4 < u \le 5, \text{ increasing}, \\
4,\ & 5 > u \ge 2, \text{ decreasing}, \\
2u,\ & 2 > u \ge 0, \text{ decreasing}, \end{cases}
\qquad
m(u) = \begin{cases} 
2u,\ & 0 \le u \le 4, \text{ increasing}, \\
u + 4,\ & 4 < u \le 5, \text{ increasing}, \\
u + 4,\ & 5 > u \ge 2, \text{ decreasing}, \\
3u,\ & 2 > u \ge 0, \text{ decreasing}. \end{cases}
\end{equation}
}
\mpcomment{
These relations are {\em rate independent}.
Since the adsorbed amount depends not just on the current value of the fluid concentration but on its history, it is a {\em hysteresis functional} of the fluid concentration denoted by $w = \hH(u)$.
Note that the path \eqref{examplefcns} would be followed for instance by the solution of the initial-value problem
}
\mpcomment{
\begin{equation}  \label{exampleode}
\begin{array}{cc}
\tfrac{d}{dt}(u + \hH(u)) = f(t), \\
(u + \hH(u))(0) = 0,
\end{array}
\quad f(t) = 
\begin{cases}
1,\ 0 \le t \le 9, \\ -1,\ 9 \le t \le 18.
\end{cases}   
\end{equation}
The same path would be followed for any source function $f(t)$ which causes $u$ to increase monotonically from $0$ to $5$ and then to decrease monotonically from $5$ to $0$.
}

\mpcomment{
Suppose an adsorbing porous medium occupies the narrow tube $\{x: x \ge 0\}$, it is fully-saturated with fluid, and initially neither contains any solute.
If fluid enters the medium at $x = 0$ with a solute concentration of $\varphi(t)$ and it flows rightward with unit velocity, the solute concentration in the fluid $u(x,t)$ and the adsorbed solute concentration $w(x,t)$ satisfy the initial-boundary-value problem
\begin{subequations} \label{ibvp}
\begin{eqnarray}
\tfrac{\partial}{\partial t}(u + w) + \tfrac{\partial}{\partial x}u = 0, \ w = \hH(u),
\\
u(x,0) = w(x,0) = 0, \ 
u(0,t) = \varphi(t),\ x > 0,\ t > 0.
\end{eqnarray}
\end{subequations}
Let the boundary values of the incoming fluid concentration be given by
$\varphi(t) = t,\  0 \le t \le 5, \textbf{ and } \varphi(t) = 10 - t,\ 5 < t \le 10.$
The fluid concentration $u(x,t)$ and the adsorbed concentration $w(x,t)$ 
within the medium satisfy this nonlinear transport equation 
and are given in the Figure~\ref{fig:example:PDE} below at discrete times $t \in \{4,5,6,7,8,9\}$. The boundary-values are translated rightward with speed $1/2$ when $0 \le t \le 4$. The characteristic speed jumps to $1$ at $x=0$ when $t=4$, and a shock develops then. (Note $u(x,5)$ and later.) A rarefraction wave is initiated at $x=0$ when $t=8$ and the characteristic speed drops to $1/3$.  (Note $u(x,9)$.)
} %% end of mpcomment

%%%%%%%%%%%

%%%%%%%%%%
\subsection*{\kgeneral\ family of models}
The approach called \kgeneral\ is an umbrella for a family of flexible models $\hH(\pP^K;u)$ calibrated from $H$ alone. Overall, the difficulty of approximating $H^* \approx H$ is not very different from that in the approximation of $\fl(\cdot),\fr(\cdot)$ by continuous piecewise linear or by step functions. The model $\hH(\pP^K;u)$ has $K$ components with parameters encoded in an array $\pP=\pP^K \in \R^{K \times 4}$ calibrated from $H$. The model can be enriched if $\dD$ includes data on internal loops. The class of \kgeneral\ models includes (i) {\em generalized play} as well as 
(ii) a discrete version \kpreisach\ of the Preisach model, as well as the most useful subclass called (iii) {\em \knonlinear}; 
these are known from the literature \cite{KrasPokr89,Mayergoyz,Visintin94}, but our calibration efforts and theory for discrete models for \eqref{eq:pde} is new. More generally, one can construct $\hH(\pP;\cdot)$ with
some infinite dimensional $\pP$ calibrated from a dense $\dD$, e.g., $\pP$ may represent the information in Preisach plane found from $\dD$.
We present a brief overview of (i-iii) now.  

 (i) The well-known generalized play model follows $\fl(\cdot)$ and $\fr(\cdot)$ exactly. It is given by an auxiliary evolution equation with time-dependent constraints
\ba
\label{eq:gplay}
\dtt {v}(t)+\cg(\fr(u(t)),\fl(u(t));v(t)) \ni 0,\;\; v(0)= \vz \in \R,
\ea
where $\cg(a,b;r)$ is a constraint graph which enforces $a\leq r \leq b$; see details in Sec.~\ref{sec:prelim}. The output can be further transformed by $w=\mu b(v(t))$ with some monotone $b(\cdot)$, and $\mu>0$. In fact, one can consider a family of $K$ \general\ models, each expressed by $w_k=\mu_kb_k(v_k)$ with $v_k$ given by \eqref{eq:gplay} with primary curves $\flk,\frk$. These are added together so $w=\sum_k w_k$; the parameters $\mu_k,\flk,\frk,b_k$ are recorded in $k$'th row of $\pP$. 
The \general\ model is conceptually simple and fairly easy to implement and is amenable to analyses. In addition, the output $\hH(\pP;u)$ exactly matches $H$ if $u$ is designed to sweep it. Its disadvantage is that the secondary scanning curves are only horizontal lines.   
 
In (ii-iii), the \kpreisach\ and \knonlinear\ models have the same functional form $w=\hH(\pP^K;u)=\sum_k \mu_k b_k(v_k)$ where each $v_k$ solves an auxiliary linear play problem problem of the form \eqref{eq:gplay} but with the constraint graphs $\fr(u)=u-\beta_k$, $\fl(u)=u-\alpha_k$, where $\alpha_k\leq \beta_k$. However, \kpreisach\ and \knonlinear\  have different properties. (ii) The discrete version \kpreisach\ of the Preisach model is built with $K$ step functions approximating the curves $\fl$ and $\fr$ thus it feature discontinuities; therefore, regularization and extra effort by nonlinear solvers is required, while a rather rugged approximation $H^*$ of $H$ emerges even when $K=O(100)$. Its advantage is that it requires very little effort in calibration.    

As a middle ground, we propose to calibrate the (iii)  \knonlinear\ model which aims to adhere to the piecewise linear interpolants of $\fl$ and $\fr$. The model has some restrictions, and may require $K=O(100)$; we give details in Sec.~\ref{sec:model}. However, the quality of $H^* \approx H$ is high, while the model can be enhanced when $\dD$ is more rich; we provide an outlook in Sec.~\ref{sec:secondary}. 
%
  
%%%%%%%%%%%
\subsection*{Numerical analysis of \kgeneral}  

The analysis of numerical schemes for parabolic PDEs given by \eqref{eq:pde} with (primarily Preisach) hysteresis was considered in many works; e.g., \cite{VerdiVis85,VerdiVis89}. 
In turn, in \cite{PeszShow20} we developed rigorous numerical analysis for the \knonlinear\ model when $A$ in \eqref{eq:pde} represents nonlinear advection. In this paper we extend these results to the \kgeneral\ model while proving some subtle auxiliary results. When combined with the analysis in \cite{PeszShow20}, these give $TV_T$ stability of an explicit upwind scheme combined with a nonlinear solver for nonlinear advection only; see Sec.~\ref{sec:analysis} and \ref{sec:numerical}. We also confirm experimentally convergence of the scheme in the $u$ variable, and stability in $w$.  Throughout, we compare the advantages and disadvantages  of the models (i-iii). \mpcomment{We discuss the computational complexity of accounting for hysteresis with our models in Sec.~\ref{sec:complexity}.}

%%%%%%%%%%%%
\begin{figure}[ht]
\centering
\begin{tabular}{ccc}
\includegraphics
[height=30mm]{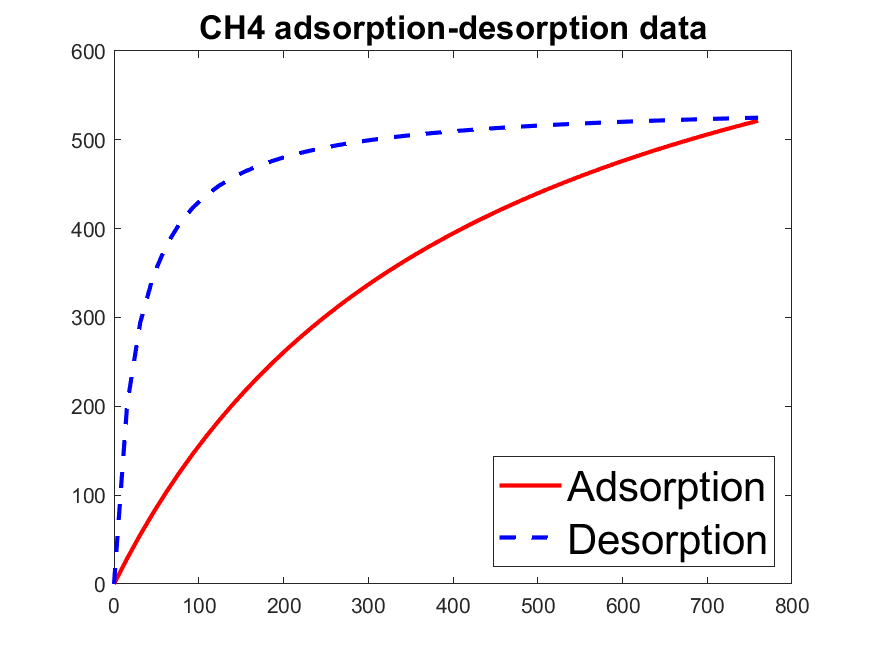}
&
\includegraphics
[height=30mm]{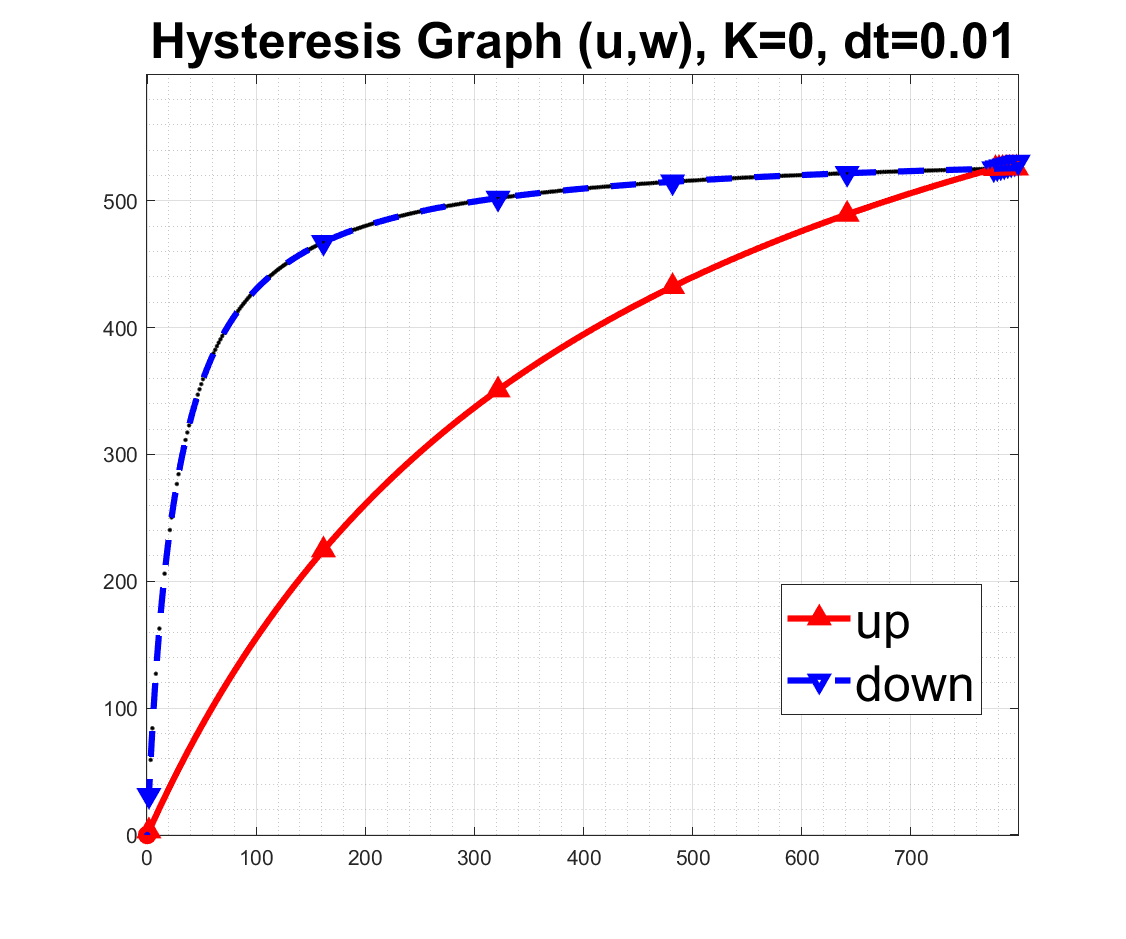}
&
\includegraphics
[height=30mm]{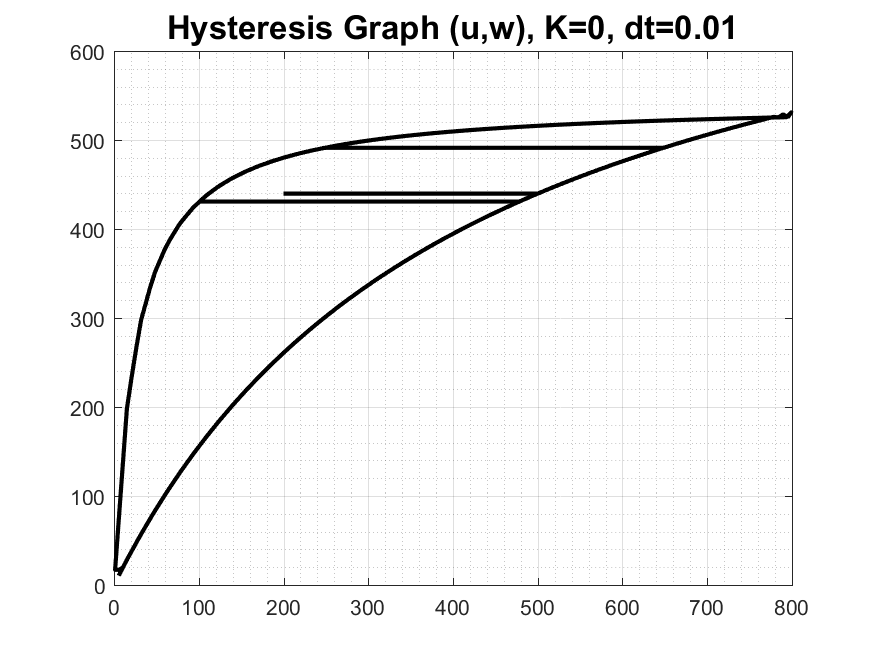}
\\
(a)&(b)&(c)
\\
\end{tabular}
\caption{
\label{fig:teaser}
Modeling hysteresis given sparse data $\dD$: (a) Experimental data $\dD$ for  $CH4$ adsorption--desorption hysteresis
with only primary scanning curves $H$; (b) The output $\hH(u)$ obtained with  \general\ hysteresis model for a particular input $u(t)$ which sweeps $H$. (c) The output $\hH(u)$with the same \general\ model  for $u(t)$  which sweeps $H$ and includes a few secondary curves.  
Details for this example along with \knonlinear\ and \kpreisach\ graphs are in Sec.~\ref{sec:adsorption}. $K=0$ denotes \general\ model.
}
\end{figure}
 
%%%%%
\subsection*{Notation}
\label{sec:notation}
Let $t$ denote time, and let $u(t)$ be an input function. We also allow $u=u(x,t)$ with $x\in \Omega$ where $\Omega \subset \R^d$ is some spatial domain; we drop $x$ when it is not relevant to the discussion. Consider the output $w(t)$ (or $w(x,t)$)  obtained by some hysteresis model $w=\hH(\pP;u;\vz)$ parametrized by a collection of parameters in $\pP$. 
The output $w(t)$ also depends on the history through the auxiliary variable $\vz$. We drop $\pP$ and $\vz$ when there is no need to single these out, and when this does not lead to confusion. 

We distinguish between the operator $\hH(\pP;\cdot)$, and its graph $\hH(\pP)=(u,\hH(\pP;u))\subset \R^2$, when the inputs $u$ are from some family $\dD_u$. At times, of interest is a fixed $u$ and the resulting trace $\hH(\pP;u)\subset \R^2$. We also denote the boundary of $\hH(\pP)$ by $H(\pP)$. In particular most useful is the family $\PL(\uU)$  of continuous piecewise linear functions on $[0,T]$ identified by their peak values (local minima and maxima) given in the sequence 
$
\uU=[\uU^0,\uU^1,\ldots \uU^S]
$
corresponding to some (increasing) collection of time steps $\tT=[T^0,T^1,\ldots T^S]$. Clearly $u\in \PL(\uU)$ is not differentiable at $t=T^m$. 
Note that the particular set $\tT$ is unimportant since hysteresis is a rate-independent process. 

Since only the derivative $\Dt{w}$ occurs in the PDE \eqref{eq:pde}, a constant can be added to $w$ without change, so we can assume without loss of generality that the hysteresis output $\hH(u)$ is non-negative.

We will consider the evolution on the time interval $(0,T]$ partitioned into discrete time steps $t_n = n \tau$ with uniform time step $\tau=T/N$. We will set $u(t_n)=u^n$ and denote the approximations $U^n \approx u_n$, with similar notation for $w(\cdot)$ and other functions. 

We will denote the identify function $x \to x$ with $\id$, and $x_+=\mathrm{max}(x,0)$. 

\subsection*{Asssumptions}
We proceed under the following conditions:
%%%%%%%%%%%%%%%%%%%%%%%%%%%%%%%%%%%%%%%%%%%%%%%%%%%%%%%%%%%%%%%%%
\bsub   \label{assumptions}
\ba     \label{assumption-a}
    a(\cdot), b(\cdot), \fl(\cdot), \fr(\cdot) 
\text{ are {continuous monotone functions} on } \R,
\\      \label{assumption-b}
    a(\cdot) \text{ is strongly monotone, }
    b(\cdot) \text{ is Lipschitz continuous, and }
\\      \label{assumption-c}
\label{eq:assumlr}
\fr(u) \le \fl(u) \text{ for } u \in \R.
\ea
\esub
Without loss of generality, we assume that $a = \id + a_0$, where $a_0$ is continuous and non-decreasing on $\R$.
\subsection*{Plan of the paper}
We discuss preliminaries in Sec.~\ref{sec:prelim}; we follow up in Sec.~\ref{sec:model} with a discussion of \kgeneral\ hysteresis models including literature notes. In Sec.~\ref{sec:inverse} we show how to calibrate $\pP$ so that a given $H$ is the boundary of $\hH(\pP)$. In Sec.~\ref{sec:analysis} we analyze the \kgeneral\ model. Section~\ref{sec:numerical} contains a discussion of a solver, stability, convergence, \mpcomment{and computational cost} of an explicit--implicit numerical scheme for \eqref{eq:pde}. In Sec.~\ref{sec:secondary} we provide an outlook towards calibration of $\hH(\cdot)$ with respect to secondary curves, and we summarize in Sec.~\ref{sec:summary}. We also provide an Appendix with additional details.

%\end{document}
%%%%%%%%%%%%%%%%%%%%%%%%%%
\section{ODE with constraint graphs and numerical approximation}
%%%%%%%%%%%%%%%%%%%%%%%%%%%%%%%%%%%%%%
\label{sec:prelim}
In this section we provide the necessary definitions and references to relevant theory for \eqref{eq:gplay} and its finite difference approximations. These are useful later in Sec.~\ref{sec:analysis} for the study of  \eqref{eq:pde}, and for $\ell^1$ or $\ell^2$ estimates for ODE systems in $\R^n$. 

%%%%%%%%%%
\subsection{ODE with maximal monotone graphs on a Hilbert space $\Hi$}

Consider a Hilbert space $\Hi$, with inner product $(\cdot,\cdot)$, and norm $\norm{\cdot}{}$. Let $\cg(\cdot)$ be a multivalued operator, i.e., a relation on $\Hi$: $\cg \subseteq \Hi \times \Hi$. Its domain $\Dom(\cg)=\{v:[v,w]\in \cg \} \subseteq \Hi$ and range $\Rg(\cg)=\{w:[v,w]\in \cg\}$ and inverse $\cg^{-1}=\{[w,v]: [v,w]\in \cg\}$ are defined as usual. We recall that $\cg$ is {\em monotone} if $(v_1-v_2,w_1-w_2)_\Hi \geq 0$ for any $v_j,w_j\in \Hi$ such that $w_j \in \cg(v_j)$; $\cg$ is {\em maximal monotone} if also $\Rg(\cg+I)=\Hi$, and then it follows that $\Rg(\lambda \cg+I)=\Hi$ for all $\lambda > 0$. Here $I$ is the identity operator.  If $\cg$ is maximal monotone and $\lambda>0$, the resolvent $R_{\cg}^{\lambda}(\cdot)=(I+\lambda \cg)^{-1}(\cdot)$ is Lipschitz continuous on all of $\Hi$, and the Yosida approximation of $\cg$ is the function $\cg^{\lambda} \equiv \frac{1}{\lambda}(I-R^{\lambda}_{\cg})$. With $\tau>0$, the range of $(I+\tau \cg)$ is $V$, and the stationary problem,
\ba
\label{eq:stationary}
v + \tau \cg(v) \ni f \in \Hi
\ea
has a unique solution given by $v=R_{\cg}^{\tau}(f)$.  The  symbol $\ni$ is used in \eqref{eq:stationary} because $\cg(v)$ is, in general, a set.  Once $v$ is found, the particular selection $c^* \in \cg(v)$ is unique and equals $\tfrac{1}{\tau}(f-v)$.  
%

%%%%%%%%%%%%
\subsubsection{Abstract Cauchy problem with a maximal monotone $\cg$ on $\Hi$} 
Let data $f \in L^1(0,T;\Hi)$ and $\vz \in \overline{\Dom(\cg)}$,
\ba
\label{eq:ode}
\tfrac{d}{dt}v(t)+\cg(v(t)) \ni f(t),\ t \in (0,T],\  v(0) = \vz. 
\ea
We choose approximations $F^n \approx f(t_n)$ and approximate $v(t_n)$ by successive finite difference solutions $V^n$ to 
%%%%%%%%%%%%%%
%
\ba \label{eq:fd}
\frac{V^n - V^{n-1}}{\tau} + \cg(V^n) \ni F^n, \quad 1 \le n \le N,\ V^0=\vz.
\ea
These solutions $V^n$ are uniquely determined since they are given by the resolvent \eqref{eq:stationary},
\ba
\label{eq:resolvent}
V^n=R_{\cg}^{\tau}(V^{n-1}+\tau F^n).
\ea
The selection $c^n\in \cg(V^n)$ is unique at each $t_n$, with $\tau c^n=V^{n-1}+\tau F^n-V^n$.

The Cauchy problem \eqref{eq:ode} is well-posed and has a unique solution $v(\cdot) \in C([0,T];\Hi)$ obtained as a limit of step-functions $v_{\tau}(t)$ built from $V^n$; see Sec.~\ref{sec:pde}.
Additionally, if $f \in W^{1,1}(0,T;\Hi)$ and $\vz \in {\Dom(\cg)}$, then $v \in W^{1,\infty}(0,T;\Hi)$  with $v(t) \in \Dom(\cg),\ t \in [0,T]$.
If $\cg$ is a subgradient, $f \in L^2(0,T;\Hi)$ and $\vz \in \overline{\Dom(\cg)}$, then $v \in W^{1,\infty}(\varepsilon,T;\Hi)$ for $\varepsilon > 0$ with $v(t) \in \Dom(\cg),\ t \in (0,T]$.
See \cite{Brezis73, Showalter97}. 

%%%%%%%%%%%%%%%%%%%%%%%%
The meaning of inclusion $\ni$ in the ODE \eqref{eq:ode} is as follows. Rewritten as $\tfrac{dv}{dt} \in f-\cg(v)$, the inclusion is replaced by an equality when $f-\cg(v)$ is replaced   by $(f(t+)-\cg (v(t)))^0$, the element of $f-\cg(v)$ of minimal norm; see \ourcite{Brezis73}{pp. 66, 28}.

%%%%%%%%%%%
\subsubsection{Convergence of \eqref{eq:fd}} 
%%%%%%%%%%%%%%%%%%%%%%%%%%%%%%%%%%%%%%%%
The proof of well-posedness discussed above relies on convergence of the step functions $v_{\tau}\to v$ as well as that of the piecewise linear interpolator $V_{\tau}$ of $V^0,V^1,\ldots, V^N$.  We note that the solutions $v(t)$ need not be smooth, even if the input $u(t)$ is smooth. Generally $O(\tau)$ is the best rate in Hilbert space, otherwise the rate is $O(\sqrt{\tau})$. In the more general context of Banach space the rate depends on data $f$, e.g., whether $f \in H^1(0,T;\Hi)$, and whether $\vz \in \Dom(\cg)$, and whether $\cg$ is a subgradient. See \ourcite{Rulla96}{Example 3}; see also \cite{NochettoSavareVerdi2000,Mielke10}, \ourcite{SimoHughes}{1.4, p41} for a-priori and a-posteriori analyses, also in application contexts such as in plasticity. 

%%%%%%%%%%
\subsection{ODE with a fixed constraint graph on $\R$}

Now we set $\Hi=\R$, let $a \leq b$, and consider the non-empty closed interval $[a,b] \subset \R$. If $a < b$ we define
\ba
\label{eq:defr}
\cab{s} = \begin{cases}
(- \infty,0] \text{ if } s = a, \\
\{0\} \text{ if } a < s < b,  \\
[0,\infty) \text{ if } s = b,  \end{cases}
\;\;
\rab{s} = \begin{cases}
a \text{ if } s \leq a, \\
s \text{ if } a <s < b, \\ 
b \text{ if } s \geq b,
\end{cases}
\; s \in \R.
\ea
This definition indicates that $\cab{\cdot}$ is set-valued; its graph will be denoted by $\ccab=\{a\} \times (-\infty,0] \cup (a,b) \times \{0\}  \cup \{b\} \times [0,\infty]$, a maximal monotone relation  on $\R \times \R$ with domain $\Dom{(\ccab)}=[a,b]$. We recall  $\cab{\cdot} = \partial I_{[a,b]}(\cdot)$ in $\R \times \R$ is the subgradient of the 
{\em indicator function} $I_{[a,b]}$ for the interval $[a,b]$: 
$I_{[a,b]}(x)= 0$ if $x \in [a,b]$, and $= + \infty$ otherwise.
It is clear that for any $\tau>0$ and any $s$ we have the equality of sets $\cab{s}=\tau \cab{s}$, so the resolvent $\rab{\cdot}=(I+\tau \cab\cdot)^{-1} =(I+\cab \cdot)^{-1}$ is independent of $\tau > 0$; it can be written as 
$\rab{s}= \mathrm{min}\{\mathrm{max}\{a,s\},b\}$. The function $\rab{\cdot}$ is a monotone piecewise linear continuous function defined on $\R$ with range $[a,b]$, differentiable except at $\{a,b\}$; it is also Lipschitz continuous with a unit Lipschitz constant. 
When $a=b \in \R$,  $\cab{s} = (-\infty,\infty)$ if $s =a=b$, $\Dom{\cab{\cdot}}=\left\{a\right\}$ is a single point, and the graph $\cab{s}$ enforces $\rab{s}=a=b$ for any $s$.  

With $\cg(\cdot)=\cab{\cdot}$, we obtain the finite difference solution to \eqref{eq:ode} using \eqref{eq:resolvent} and the resolvent $R_{\cg}^{\tau}(\cdot)=\rab{\cdot}$. Alternatively, one can replace $\cab{\cdot}$ in \eqref{eq:ode} by a smoother ``penalty functional'' which enforces $v(t)\in \Dom(\cab{})=[a,b]$. Another possibility is to use a Lagrange multiplier, but elimination of Lagrange multiplier typically gives exactly \eqref{eq:resolvent}. 

\subsection{ODE with a time-dependent constraint}
\label{sec:tgraph}
In {generalized play} models of hysteresis \eqref{eq:gplay}  the constraints in $\cab\cdot$ are time dependent and in fact depend on the input function $u(t)$, namely, $a(t) = \fr(u(t))$, $b(t) = \fl(u(t))$. Here $\fr,\fl$ are {continuous} monotone functions, and $u(\cdot) \in C([0,T])$.  
We consider the IVP for \eqref{eq:gplay}
\ba
\label{eq:odeg}
\tfrac{d}{dt}v(t)+\cg(\fr(u(t)),\fl(u(t));v(t)) \ni 0,\ t \in (0,T],\  v(0) = \vz. 
\ea
The approximation of \eqref{eq:odeg} requires that we know $U^1,U^2,\ldots$ and then solve successively for $n \ge 1$
\begin{equation}  \label{eq:fdgenplay}
\frac{1}{\tau}(V^n-V^{n-1})  + \cg(\fr(U^n),\fl(U^n);V^n) \ni 0\,,
\quad V^0 = \vz\,.
\end{equation}
Given $U^n$, and a fixed pair $\fl,\fr$, we can write out the solution $V^{n}$  of this stationary problem, adapting \eqref{eq:defr} to define $\resg(\vbar;\cdot)$ with $\vbar=V^{n-1}$
\ba
\label{eq:vsol}
V^n = \resg(\vbar;U^n) \equiv \resab{\fr(U^n)}{\fl(U^n)}{\vbar} = \begin{cases}
\fr(U^n)\text{ if } \vbar \leq \fr(U^n), \\
\vbar \text{ if } \fr(U^n) <\vbar < \fl(U^n), \\ 
\fl(U^n) \text{ if } \vbar \geq \fl(U^n).  \end{cases}
\ea
Various properties of \eqref{eq:odeg}--\eqref{eq:vsol} are needed in Sec.~\ref{sec:analysis} and \ref{sec:numerical} when \eqref{eq:odeg} is coupled with an evolution problem for $u(t)$.
In particular, each $\resg(\vbar;u)$ is differentiable at the points of differentiability of $\fl(\cdot)$ and $\fr(\cdot)$ except at
$u=\fl^{-1}(\vbar)$ and at $u=\fr^{-1}(\vbar)$. In addition,
we have the following monotonicity result
\begin{lemma}
\label{lem:gres}
Assume $V^{n-1} \in \Dom (\cg(\fr(U^{n-1}),\fl(U^{n-1});\cdot)$. Then  
%%%
\ba
\label{eq:order}
U^n\geq U^{n-1} \implies V^n=\ggres(V^{n-1};U^n) \geq V^{n-1},\\
U^n\leq U^{n-1} \implies V^n=\ggres(V^{n-1};U^n) \leq V^{n-1}.
\nonumber
\ea
\end{lemma}

\begin{proof}
This property might or not be obvious, and is easiest to prove when $\fl$ and $\fr$ are injective. From the assumption $V^{n-1} \in \Dom (\cg(\fr(U^{n-1}),\fl(U^{n-1});\cdot))$ which means 
$\fl^{-1}(V^{n-1}) \leq U^{n-1} \leq \fr^{-1}(V^{n-1})$. Now  the point $(U^n,V^n)$ with $V^n=\ggres(V^{n-1};U^n)$ is on the the graph of the monotone nondecreasing function \eqref{eq:vsol} determined by $\fl(\cdot)$ on the left, $\fr(\cdot)$ on the right, with a ``flat'' connector at $v=V^{n-1}$. Thus $V^n\geq V^{n-1}$ whenever $U^n\geq U^{n-1}$ (the second part follows analogously).  
In the non-injective case we replace $\fl^{-1}(V^{n-1})$, 
$\fr^{-1}(V^{n-1})$ by
$\min(\fl^{-1}(V^{n-1}))$ and  
$\max(\fr^{-1}(V^{n-1}))$, respectively.
\end{proof}

The special case of \eqref{eq:odeg} with $\fr(u)=u-\beta$, $\fl(u)=u-\alpha$ for some $\alpha\leq \beta$ gives 
\ba
\label{eq:odelp}
\tfrac{d}{dt}v(t)+\cg(u(t)-\beta,u(t)-\alpha;v(t)) \ni 0,\ t \in (0,T],\  v(0) = \vz,
\ea
with the approximation
%%
%\begin{equation}  \label{eq:fdlplay}
$
\frac{1}{\tau}(V^n-V^{n-1})  + \cg(U^n-\beta,U^n-\alpha;V^n) \ni 0\,,
\quad V^0 = \vz$.
%\end{equation}
%%
The counterpart of \eqref{eq:vsol} with $\vbar=V^{n-1}$ can be written in two equivalent ways,
\ba
\label{eq:resab}
V^n=U^n+\resab{-\beta}{-\alpha}{\vbar-U^n}=
\resab{U^n-\beta}{U^n-\alpha}{\vbar}.
\ea

%%%%%%%%%%%%%%%
\subsection{Auxiliary implicit ODE: from v to $w=b(v)$}
\label{sec:b}
We recall now the following subtle relationship. 
\begin{lemma}
\label{lem:b}
Assume $u \in C([0,T])$ and that $v \in W^{1,1}(0,T)$ is a strong solution of \eqref{eq:odeg}. Then $b(v(t))$ is the unique solution  determined by $b(v(0))$ of
\ba
\dtt{b(v(t))}+\cg(\fr(u(t)),\fl(u(t)); v(t)) &\ni& 0.    \label{eq:odeb}
\ea

\end{lemma}
%%%%%%%%%%%%%%%%%%%%%%%%%%%%%%%%%%%%%%%%%%%
\begin{proof}
Let $v \in W^{1,1}(0,T)$ be a strong solution of \eqref{eq:odeg}.
Since $b(\cdot)$ is Lipschitz, $w(t) \equiv b(v(t))$ is differentiable a.e., and the chain rule gives
(by \ourcite{KindStam}{Cor. A.6})
\begin{eqnarray*}
w'(t) = b'(v(t))v'(t) \in - b'(v(t))\, \cg(\fr(u),\fl(u); v(t)) \subset - \cg(\fr(u),\fl(u); v(t)),
\end{eqnarray*}
where the last relation follows from $b'(\cdot) \ge 0$.

If $v_j$ are strong solutions for $j = 1,2$, set $w_j(t) = b(v_j(t))$, so  $w_j'(t) + \cg(\fr(u(t)),\fl(u(t)); v_j(t)) \ni 0.$ From Theorem A.1 of \cite{KindStam}, the absolutely continuous function $|w_1(t) - w_2(t)|$ is a.e. differentiable and satisfies
\begin{eqnarray*}
\tfrac{d}{dt} |w_1(t) - w_2(t)|  
=  \sgn_0(w_1(t) - w_2(t)) (w_1'(t) - w_2'(t)) 
\\
\in - \sgn_0(w_1(t) - w_2(t)) \big( \cg(\fr(u),\fl(u); v_1(t)) - \cg(\fr(u),\fl(u); v_2(t)) \big) .
\end{eqnarray*}
If $w_1(t) - w_2(t) \ne 0$, then the last term is non-positive for any choices from the constraint relations, and so $|w_1(t) - w_2(t)| \le |w_1(0)-w(0)|$ for $t \ge 0$.
\end{proof}

We also consider a discrete analogue of Lemma~\ref{lem:b} for \eqref{eq:odeb}, and apply Lemma~\ref{lem:gres}.
\begin{lemma}
\label{lem:bn}
If \eqref{eq:fdgenplay} holds, then $W^n=b(V^n)$ is the unique solution of
\ba
\label{eq:wsol}
\frac{1}{\tau}(b(V^n)-b(V^{n-1}))  + \cg(\fr(U^n),\fl(U^n);V^n) \ni 0\,,
\quad b(V^0) = b(\vz)\,.
\ea
In addition, $U^n\geq U^{n-1}$ implies $W^n\geq W^{n-1}$. 
\end{lemma}
%%

%%%%%%%%%%%%%%%
\section{Hysteresis models: generalized play, \knonlinear, \kpreisach, and related}  \label{sec:model} 
%%%%%%%%%%%%%%%
Mathematical models of hysteresis have a long history, and much work has been devoted to their development and analysis. We refer to the monographs \cite{Visintin93,Visintin94,Mayergoyz,KrasPokr89} and the review paper \cite{Macki93} for overview and the detailed history of a large variety of hysteresis models. 

%%%%%%%%%%%%%
\begin{table}
\begin{tabular}{|l|c |c |c |c|}
\hline
&Linear&Nonlinear&Preisach&Generalized
\\
\hline
\hline
$\cg$&
\multicolumn{3}{c|}{$\cg(u-\beta,u-\alpha;v)$}
&$\cg(\fr(u),\fl(u);v)$
\\
\hline
$u\to v$&\multicolumn{3}{c|}{$u-\beta\leq v\leq u-\alpha$}
%&$u-\beta\leq v\leq u-\alpha$
&
$\fr(u)\leq v\leq \fl(u)$
\\
&\multicolumn{3}{c|}{$\alpha \leq u-v\leq \beta$}
&
$\fl^{-1}(v)\leq u \leq \fr^{-1}(v)$
\\
\hline
%%%%%%%%%%%%%%%%%%%%%%
%$v\to w=\mu \bh(v)$
%\\
truncation $b(v)$
&id
&$\bhh$ 
&$\bhz$, $\bheps$ or $\bhstar$ 
&id
\\
\hline
%%%%%%%%%%%%%%%%%%
primary
&$v=u-\beta$ or $v=u-\alpha$
&
\multicolumn{2}{c|}{translate of $b(\cdot)$}
& $\fl(u)$ or $\fr(u)$
\\
\hline
secondary&horizontal&
\multicolumn{2}{c|}{horizontal}
&
horizontal
\\
\hline
parameters in $\pP$&
$\mu,\alpha,\beta,\infty$
&
$\mu,\alpha,\beta,h$
&
$\mu,\alpha,\beta,\bh_s$
&
$1,\fl(\cdot)$, $\fr(\cdot),\id$
\\
\hline
%%%
reference
&\cite{PeszShow98}
&
\multicolumn{2}{c|}{III.2, p64}
&III.2, p65
\\
\hline\hline
\end{tabular}
\caption{Hysteresis (hysteron) operator 
choices. The input is $u(t)$. The output $v(t)$ satisfies 
$\fr(u)\leq v\leq \fl(u)$, and is scaled and truncated as in $\hH(\pP;u;\vz)\ni w(t) = \mu b(v(t))$. For \nonlinear, we have $b(\cdot)=\bh(\cdot)$ and we record its parameter $h$. For Preisach hysteresis, $(s)=$0, or $\eps$, or $*$, which indicates $b(\cdot)=\bhz,\bheps,\bhstar$, respectively.  For \general\ in principle we can have any $\mu$ and $b(\cdot)$, but   these can be subsumed in the definition of $\fl,\fr$, thus for simplicity we set $\mu=1$ and $b(\cdot)=\id$. 
\label{tab:unit}
}
\end{table}
%%% 

In this paper we focus on three types of play hysteresis models under a common umbrella of {\kgeneral}: \general, \knonlinear, and regularized \kpreisach, which we analyze, compare, and parametrize.  The three types are interconnected. Generalized play can be approximated by \knonlinear. Furthermore, equivalent representation of Preisach model can be obtained as a superposition of an infinite number of unit hysterons of type \nonlinear; see \ourcite{Mayergoyz}{p.31 and Fig. 1.32}. 
These models are constructed with three steps which give output $w$ to input $u$ by adding $K$ unit hysterons $w_k$. 
Examples are given in Sec.~\ref{sec:illustrate}. 

\subsubsection*{{\bf (A)} Play models.} We build unit hysterons with initial value problems (IVP) for either generalized play \eqref{eq:odeg} or its special case, linear play \eqref{eq:odelp}, with  some given $\alpha\leq \beta$, and primary curves $\fr(u)=u-\beta,\fl(u)=u-\alpha$.
These models give output $v(t)$ which increases on $v(t)=\fr(u(t))$, decreases on $v(t)=\fl(u(t))$, and is constant between these bounding curves where it satisfies $\fr(u(t)) < v(t) < \fl(u(t))$. We consider $K$ such auxiliary functions $v_k(u)$, each corresponding to its own $\flk,\frk$. 

\subsubsection*{{\bf (B)} Truncation of play models.} The second step is to truncate each hysteron to limit the influence of the constraint. The output of the truncation is $b(v)$ as in Sec.~\ref{sec:b}, and is $b_k(v_k(u))$ for each $v_k$,
where each $b_k(\cdot)$ is some arbitrary monotone nondecreasing function, possibly different for each $k$. Models $u\to b_k(v_k(u))$ are called \nonlinear\ \cite{Visintin94}. The shape of each $b_k(v_k(u))$ follows the translates of $b(\cdot)$ (or $b_k$).

In particular we choose $b(\cdot)$ to be either $\bhinf=\id$ for the \linear, or  a truncation function with bounded range $[0,h]$ for \nonlinear.  Let $\eps>0,h>0$ and define the scaled ramp function  $\bh_{\eps}(x)=\tfrac{h}{\eps}(x_+-(x-\eps)_+)$, which has range $[0,h]$, slope $\tfrac{h}{\eps}$ on $(0,\eps)$, and equals $0$ for $x<0$ and $h$ for $x > \eps$.  The ramp function $\bh(x)=x_+-(x-h)_+$ is a particular case, with maximum slope $1$, and is the main building block in \knonlinear\ models, with unit hysterons of shape of truncated parallelograms. 

The scaled left continuous Heaviside function  $\bh_0$  equal to  $\bh_0(x)=h$ when $x>0$ and $\bh_0(x)=0$ when $x \leq 0$ has ``maximum slope'' equal to $\infty$, and produces discontinuous outputs. The output $(u,w)$ forms a hysteron  ``box'' of height $\mu h$ which is a building block of the discontinuous Preisach model, also called {\em basic relay model} \ourcite{Visintin94}{p.97}. These can be approximated by  their Yosida approximations $\bheps(\cdot)$ as $\eps\to 0$; we note $\bh=\bheps$ when $\eps=1$. Also, $\bh_0$ can be approximated by some $C^{\infty}$ smooth function $\bh_*(v)$ with range $[0,h]$; here we use the appropriately scaled {\tt erf} function
\ba
\label{eq:bsmooth}
\bh_*(s)= h/2( \mathrm{erf}(2s-1)+1/2)),
\ea
but other choices are possible. 

\begin{figure}[ht]
\begin{tabular}{ccccc}
\includegraphics
[height=20mm]{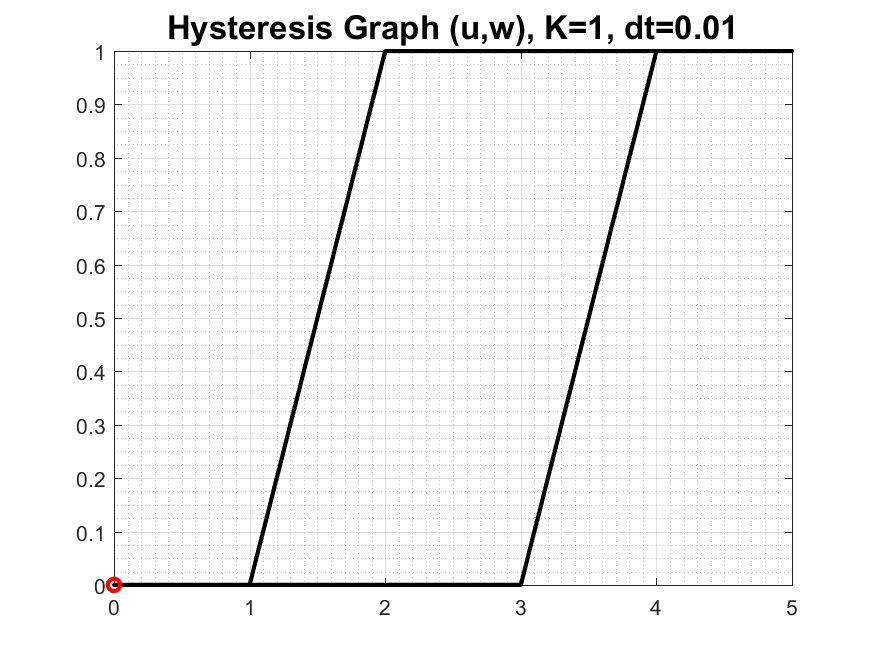}
&
\includegraphics
[height=20mm]{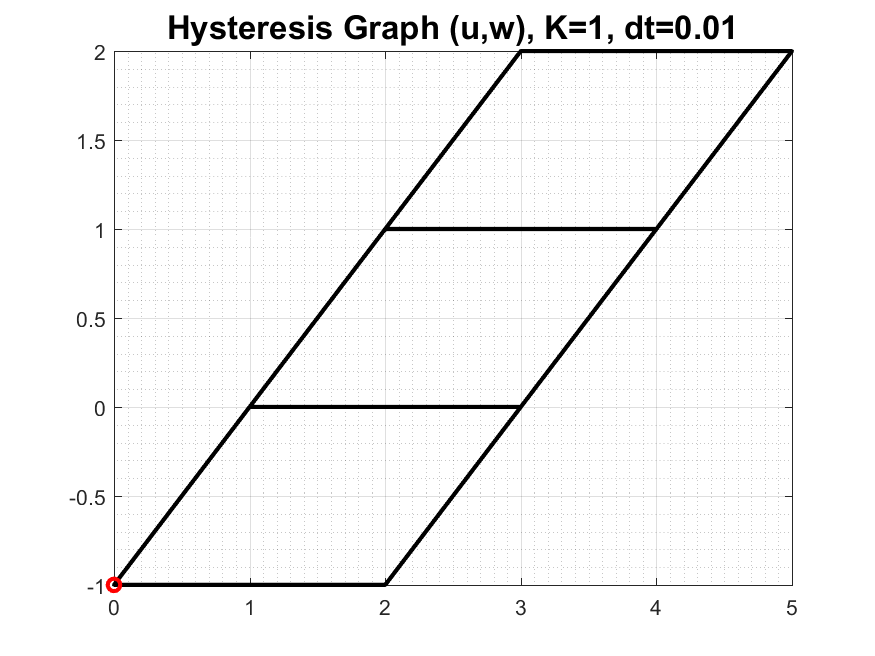}
&
\includegraphics
[height=20mm]{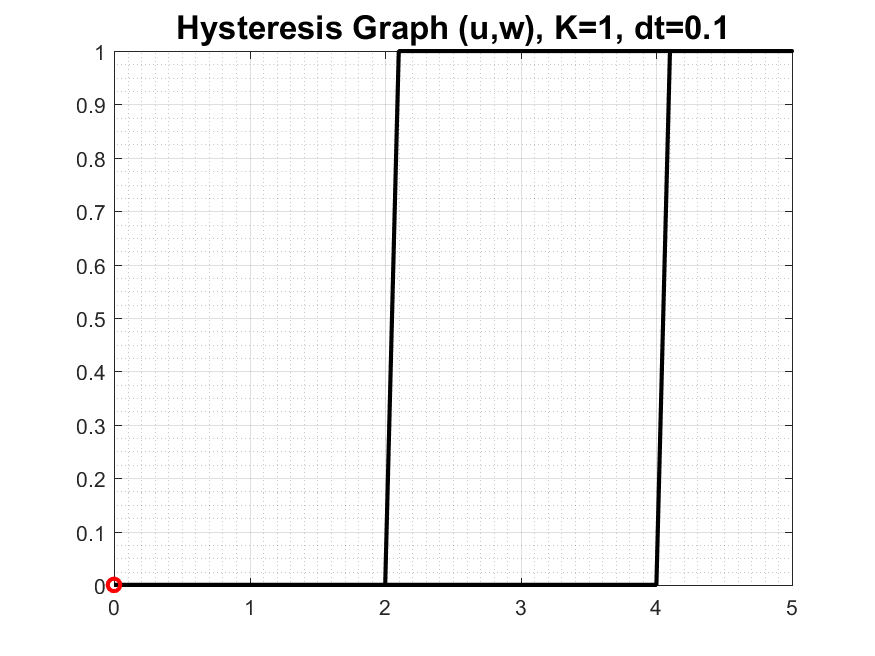}
&
\includegraphics
[height=20mm]{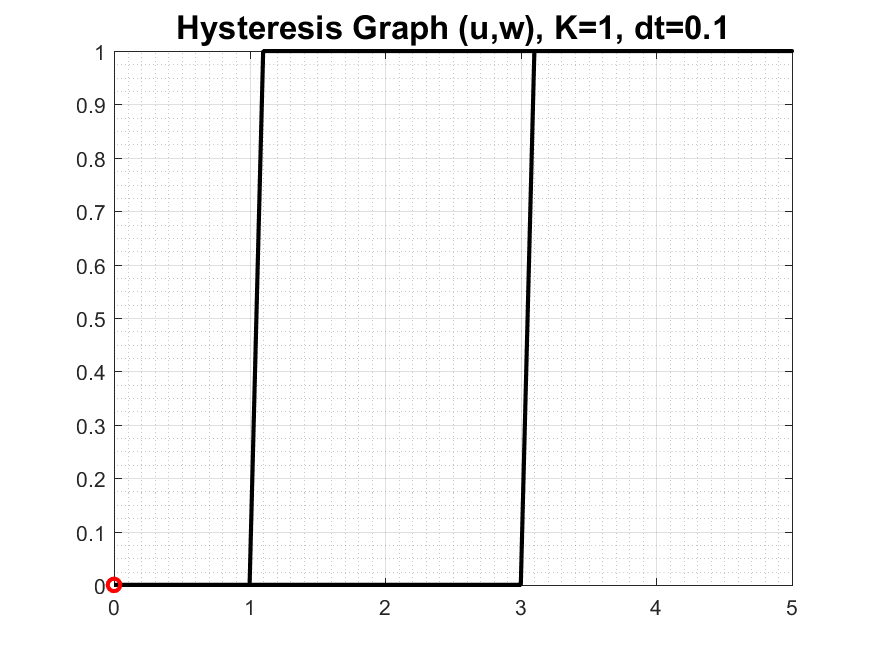}
&
\includegraphics
[height=20mm]{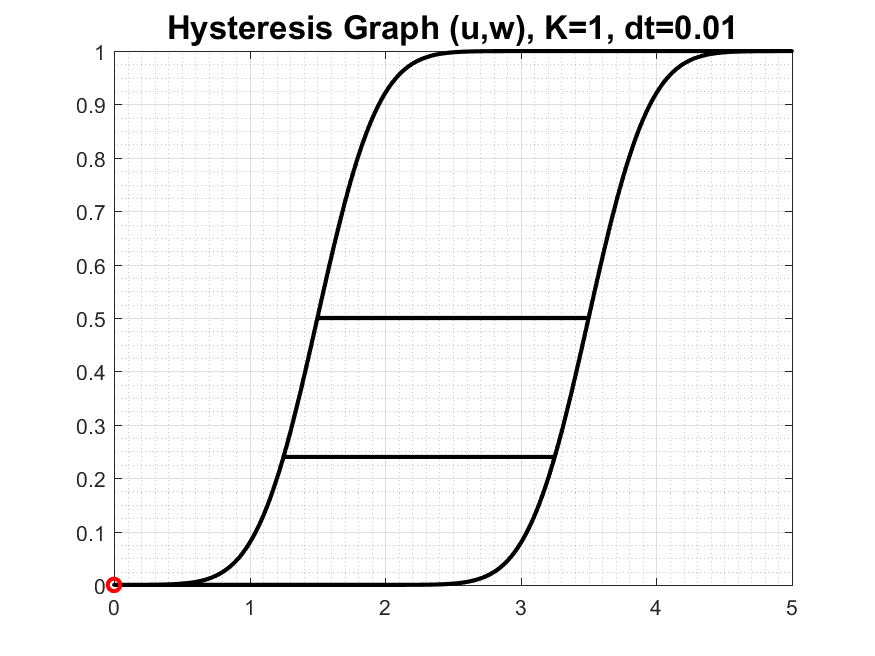}
\\
$\pP^{1,1}$=$[1,1,3,1]$&
$\pP^{1}_{\infty}$=$[1,1,3,\infty]$&
$\pP^{1,0}$=$[1,1,3,b^1_0]$&
$\pP^{1,\eps}$=$[100,1,3,1/100]$&
$\pP^{1,*}$=$[1,1,3,*]$
\\
(i)&(ii)&(iii)&(iv)&(v)\\
nonlinear&linear&Preisach&\multicolumn{2}{c}{(iv-v) regularized Preisach}\\
\end{tabular}
\caption{Examples of unit hysterons from Sec.~\ref{sec:model} with $\fr(u)=u-3$ and $\fl(u)=u-1$. Plotted are $H=\hH(\pP;u)$ with  $\pP$ as indicated and $u=\PL(0,5,0)$ in all cases except (ii) where $u=\PL(0,5,0,4,1,3)$, and (v) where $u=\PL(0,5,0,0,3.5,1.25)$. The interior loop for (ii) is the same as the truncated hysteron in (i). In all examples we set $\vz=\fl(u(0))$, and $\wz=0$ in all cases except (ii) when $\wz=\vz=\uz-1=-1$. The point $(\uz,\wz)$ is shown by a red dot.  
\label{fig:unit}
}
\end{figure}

The steps (A-B) produce unit hysterons shown in Fig.~\ref{fig:unit}. Note that a unit hysteron with $\bh_0$ has vertical sides and 
\eqref{assumption-b} does not hold, thus
\kpreisach\ model without regularization is not one of \kgeneral\ models. However, we include it in some illustrative examples.

\subsubsection*{{\bf (C)} Linear combinations of unit hysterons.} 
%%
%%%%%%%%%%%
\begin{definition}
\label{def:kgeneral}
The \kgeneral\ model is determined by a family of constraint curves $\{\frk, \flk\}$ and functions $b_k$, and scaling factors $\mu_k>0$. The output is  
\bsub
\label{eq:pk}
\ba 
\label{eq:bcollect}
w(t) = \hH(\pP;u(t);\vz) &=& \sum_{k=1}^K w_k = \sum_{k=1}^K \mu_k b_k(v_k(u(t))),
\\
\label{eq:kcom}
\tfrac{d}{d t} v_k(t)  + \cg(\frk(u),\flk(u);v_k(t))  &\ni& 0,\
\\
\nonumber
\quad v_k(0) &=& \vz_k \in [\frk(\uz),\flk(\uz)],
\quad 1 \le k \le K.
\ea
\esub
We collect the parameters in a K-tuple $\pP^K=(\pP_k)_{k=1}^K$. We also assume for the relevant data that
\bsub
\label{eq:pi}
\ba
\flk,\frk,b_k 
\mathrm{\ satisfy \eqref{assumption-a}, \eqref{assumption-b}, \eqref{eq:assumlr}}, 
\\
b_k \mathrm{\ are\ one\ of\ those\ listed\ in\ (B)},
\\
\alpha_k\leq \beta_k; h_k>0,\mu_k>0. 
\ea
\esub
\end{definition}

Each row $\pP_k$ of $\pP$ represents a {\em unit hysteron} identified by either some functions or numbers or special symbols, with interpretation clear from the context, as in an object-oriented software environment. 
For example, the numbers are interpreted as parameters of some fixed functions, and the symbol $\infty$ or $\ast$ have a special meaning. 
Table~\ref{tab:unit} summarizes the notation for unit hysterons, and Table~\ref{tab:kmodels} the properties and notation for the family of \kgeneral. For simplicity we consider only hysteresis operators made of   unit hysterons of the same type, even though  our theoretical results as well as algorithms apply to the more general case. 

We have the special cases of $[\mu_k,\flk(\cdot),\frk(\cdot),b_k(\cdot)]$ denoted by 
\bsub
\label{eq:pkall}
\ba
\label{eq:pkgen}
\pgen \ni \pP; \;\;
\pP_k&=&[1,\flk(\cdot),\frk(\cdot),\id], \;\; \mathrm{or}
\\
\label{eq:pknon}
\pknon \ni \pP; \;\;
\pP_k&=&[\mu_k,\alpha_k,\beta_k,h_k],
\;\; \mathrm{or}
\\
\label{eq:pklin}
\pklin \ni \pP;\;\;
\pP_k&=&[\mu_k,\alpha_k,\beta_k,\infty],
\;\; \mathrm{or}
\\
\label{eq:pkpreisach}
\ppreisach \cup \ppreisache \cup \ppreisachs \ni \pP;\;\;
\pP_k&=&[\mu_k,\alpha_k,\beta_k,\bkh_{s_k}(\cdot)]. 
\ea
\esub
The notation in \eqref{eq:pkall} is similar to MATLAB matrix notation. 

Additional remarks on \eqref{eq:pkall} are as follows. In \eqref{eq:pkgen}, from a modeling point of view,  it makes sense in for \general\ to subsume $\mu$ and $\bh$ in the definitions of $\fl,\fr$, since the alternative leads to cumbersome calibration. In addition, in practice \general\ model uses only one component, but the keyword \kgeneral\ is useful to denote the entire umbrella of models when discussing theory and implementation. 
For Preisach model in \eqref{eq:pkpreisach}, the maximum slope $s$ of $\bkh_{r_k}(v)$ is the maximum of $1/r_k$, and is the superscript in the parameter array $\ppreisache$, e.g., we have $s=\tfrac{1}{\epsilon}$  for $\bheps$. For the smooth choice $b_k=b^{h_k}_*$ we denote the variable slope by $*$ in $\ppreisachs$. Also, 
$r=0$ for $\bh_0$ with ``infinite'' slope $s=\infty$, thus $\bh_0$ is not Lipschitz and not part of \kgeneral\ family.
Finally, the special notation with symbols $\infty$ or $*$  need not to be interpreted literally in the formula \eqref{eq:bcollect}. 

%%%%%%%%%%%%%%%%%%%%%%%%%%
%%%%%%%%%%%%%%%

%%%%%%%%%%%%%%%
%\subsection{$K$-component family of hysterons for \knonlinear}
%\label{sec:mknon}
%%%%%%%%%%%%%
\begin{table}
\begin{tabular}{l||l|l|l|}
&\klinear&\knonlinear&{\kpreisach}\\
\hline
\hline
%%%
&$\pklin$&$\pknon$&$\ppreisach,\ppreisache,\ppreisachs$
\\
%%%
%%
row of $\pP^K$
&
$(\mu_k,\alpha_k,\beta_k,\infty)$
&
$(\mu_k,\alpha_k,\beta_k,h_k)$
&
$(\mu_k,\alpha_k,\beta_k,\bkh_{r_k})$
\\
$w =\hH(u)=$
&$\sum_k \mu_k v_k$
&
$\sum_k \mu_k \bkh(v_k)$
&
$\sum_k \mu_k \bkh_{r_k}(v_k)$
\\
\hline
primary
&
convex symmetric
&
unions of trapezoids
&
monotone
\\
&
piecewise linear
&
with monotone sides
&
steep stair-steps 
\\
\hline
secondary&
rich&
possibly rich&
possibly rich
\\
\hline\hline
\end{tabular}
\caption{\label{tab:kmodels}\knonlinear\ and \kpreisach\ hysteresis models $\hH(\pP;\cdot)$. Each component $v_k$ is the solution to \eqref{eq:kcom}, with finite difference solution $V_k^n$ given by \eqref{eq:resab}. The symbol $r_k$ as in Table~\ref{tab:unit} indicates  one of $0,\eps,*$; we recall that $\bh_0$ cannot be used in \kgeneral.}
\end{table}
%%% 

%%%%%%%%%%%%%
\subsection{Properties of $\hH(\pP;\cdot)$ with $\pP$ is as in Def.~\ref{def:kgeneral}} 

We provide brief remarks; see also  illustrations in Fig.~\ref{fig:stack} and Fig.~\ref{fig:stack-preisach}. 

For \klinear\ with $\pP=\pklin$, the finite sums of positive multiples of linear-play functionals over a collection of constraint intervals yield a linear play which contains {\em internal loops} consisting of a convex right-constraint for increasing values and a corresponding center-symmetric concave left-constraint for decreasing values.  The convex-concave character of the internal loops arises from the fact that the linear-play functionals are not truncated, so the {\em slope} of their sum is monotone with respect to the input. That is, once a constraint is active, it remains active until the input reverses direction; see \cite{PeszShow98}, \ourcite{Visintin94}{p. 84}; see also Fig.~\ref{fig:convex} from Sec.~\ref{sec:convex}. Interestingly, most of work on Preisach model features such symmetric convex-concave graphs. 

In contrast, in the \knonlinear\ model when $\pP=\pknon$, the variation of the $k$'th constraint is localized to the interval $[\alpha_k,\beta_k + h_k]$. The secondary curves depend significantly on the mutual arrangement of the parameters $\alpha_k,\beta_k$. Examples are shown in Fig.~\ref{fig:stack}. 
We come back to this impact on secondary scanning curves in Sec.~\ref{sec:secondary}. 

For the regularized \kpreisach\ model and $\pP=\ppreisache$ or $\ppreisachs$, the output $w$ is made of ``stair steps'' with steep slopes intermingled with some flat pieces. During calibration we actually set-up the \kpreisach\ model $\ppreisach$ which we later regularize  with  $\hH(\ppreisache;\cdot)$, but we do not attempt to eliminate the flat pieces unlike with $\pknon$; see comparison in Fig.~\ref{fig:stack-preisach}.

We mentioned earlier that we exclude the \kpreisach\ model with $\ppreisach$ from \kgeneral\ family. We recall that it is discrete, a sum of $K$ positive multiples of a family of delayed relay functionals 
$\pP^{1,0}$=$[1,\alpha_k,\beta_k,b^1_0],\ \alpha_k \le \beta_k$, and the output is discontinuous. Given $u$, we {\em can} produce (very rough) $w$ with $\ppreisach$. However, $\hH(\ppreisach;\cdot)$  is not tractable by a numerical solver when solving  for $u$ and $w$, e.g., in \eqref{eq:pde}. 

Preisach operator {\em can} produce smooth output \ourcite{Mayergoyz}{p.31},  \ourcite{Visintin94}{Chapter 4}
 if an uncountable collection of measures $\mu(\alpha,\beta,h)$ is given to create
$$\hH(u) = \iint_{\alpha < \beta} \mu(\alpha,\beta) b^1_0(v_{\alpha,\beta}(u)) d\alpha d\beta$$ where $v_{\alpha,\beta}(u)$ solves \eqref{eq:odelp}. Such an operator allows rich interior cycles, however the calibration necessary to obtain a particular model requires dense data $\dD$ 
\cite{HoffMeyer89,VerdiVis85,VerdiVis89,Krejci13}.
More generally, these are all examples of Prandtl-Ishlinskii play hysteresis; see \ourcite{Visintin94}{Ch III} for perspectives.

%%%%%
\subsection{Practical use of \kgeneral\ in numerical schemes}

The model \eqref{eq:pk} uses ODEs \eqref{eq:kcom} to define $v_k$ and the output $w(t)$ for input $u(t)$, $u \in C([0,T])$. In a numerical scheme, either $(U^n)_n$ are given as input, or they are themselves unknown. In the approximation scheme, we do not need actually to solve the ODEs \eqref{eq:kcom} for $V_k^n$. Rather, we have resolvent formulas \eqref{eq:vsol} which define the approximations $V_k^n$ and Lemma~\ref{lem:b} which defines $b_k(V_k^n)$. For concise notation, recalling the definition of $R$ in \eqref{eq:defr}, we adapt the formulas for the discrete version of \eqref{eq:pk}, and set $\resk$ to denote the appropriate resolvent for each component $\pP_k$
\ba
\label{eq:resk}
\resk(\vbar;U)=\resab{\frk(U)}{\flk(U)}{\vbar}.
\ea
For \knonlinear\ or \kpreisach\ models, the $\resk(\vbar;U)=\resab{U-\beta_k}{U-\alpha_k}{\vbar}$. 

%%%%%%%%%
\subsection{Illustration of \kgeneral\ models}
\label{sec:illustrate}
We now show examples of $H=\hH(\pP;u(t);\vz)$ parametrized with different $\pP$. 
A variety of unit hysteron shapes obtained with \eqref{eq:bcollect} is shown in Fig.~\ref{fig:unit}. Models with \general\ have already been shown; e.g., Fig.~\ref{fig:teaser}. We focus thus on $K>1$ and the $\pknon$ and $\ppreisach$ family. We make a uniform choice $\vz_k= \flk(u_0)$, compatible with \eqref{eq:kcom}, except as indicated. 

\myskip{
(i, ii) The simplest unit hysteron from \knonlinear\ family is given by 
\bas
\pknon&=&
\begin{tabular}{|l|l|ll|l|}
\hline
$k$&$\mu_k$&$\alpha_k$&$\beta_k$&$h_k$\\
\hline
1&1&
1&3&1
\\
\hline
\end{tabular}
= [1,1,3,1].
\eas
Its graph $\hH$ includes the parallelogram $H$ with sides of slope $1$ and bottom base $[1,3]$ where intersect the $v=0$ axis. The play output is truncated by the unit ramp function $\bhh(\cdot)$ with $h=1$. The corresponding \klinear\ model shown in (ii) is encoded with $h_1=\infty$ and we write
\bas
\pklin&=&
\begin{tabular}{|l|l|ll|l|}
\hline
$k$&$\mu_k$&$\alpha_k$&$\beta_k$&$h_k$\\
\hline
1&1&
1&3&$\infty$
\\
\hline
\end{tabular}
= [1,1,3,\infty].
\eas
Its graph $H$ is not truncated; it includes the segments lying on two lines $w=v=u-1$ and $w=v=u-3$.

(ii) Preisach hysteron features discontinuous output with ``infinite'' slopes and
\bas
\ppreisach&=&
\begin{tabular}{|l|l|ll|l|}
\hline
$k$&$\mu_k$&$\alpha_k$&$\beta_k$&$h_k$\\
\hline
1&$\infty$&
1&3&1
\\
\hline
\end{tabular}
= [\infty,1,3,1].
\eas
(iii) This hysteron can be approximated by 
\bas
\ppreisache&=&
\begin{tabular}{|l|l|ll|l|}
\hline
$k$&$\mu_k$&$\alpha_k$&$\beta_k$&$h_k$\\
\hline
1&$1/\eps$&
1&3&$\eps$
\\
\hline
\end{tabular}
= [1/\eps,1,3,\eps].
\eas
(v) One can also use  a smooth $\bh_*$; see Fig.~\ref{fig:unit}(v). 
}
%%%%%%%%%%%%%%%%%%%%%%%%%%%
\myskip{
\subsubsection{Examples of stacking hysterons with $K>1$}
\label{sec:kmodels}
}
\myskip{
For example, for $K=2$ we can write $\pP^2$ in two equivalent ways
\bas
\pP^2 =
\begin{tabular}{|l|l|ll|l|}
\hline
$k$&$\mu_k$&$\alpha_k$&$\beta_k$&$h_k$\\
\hline
1&$\mu_1$&
$\alpha_1$&$\beta_1$&$h_1$
\\
2&$\mu_2$&
$\alpha_2$&$\beta_2$&$h_2$
\\
\hline
\end{tabular}
= [\mu_1,\alpha_1,\beta_1,h_1;\mu_2,\alpha_2,\beta_2,h_2].
\eas
}

%%%%%%%%%%%%%%%%%%%%%%%
When $K>1$, the shape of primary scanning curves $u \to w$  in $H$ as well as of the secondary scanning curves  depends on the mutual arrangement of $\alpha_k$ and $\beta_k$, as well as on how the unit hysterons are stacked, truncated and scaled.   

\begin{figure}[ht]
\begin{tabular}{ccc}
\includegraphics
[height=30mm]{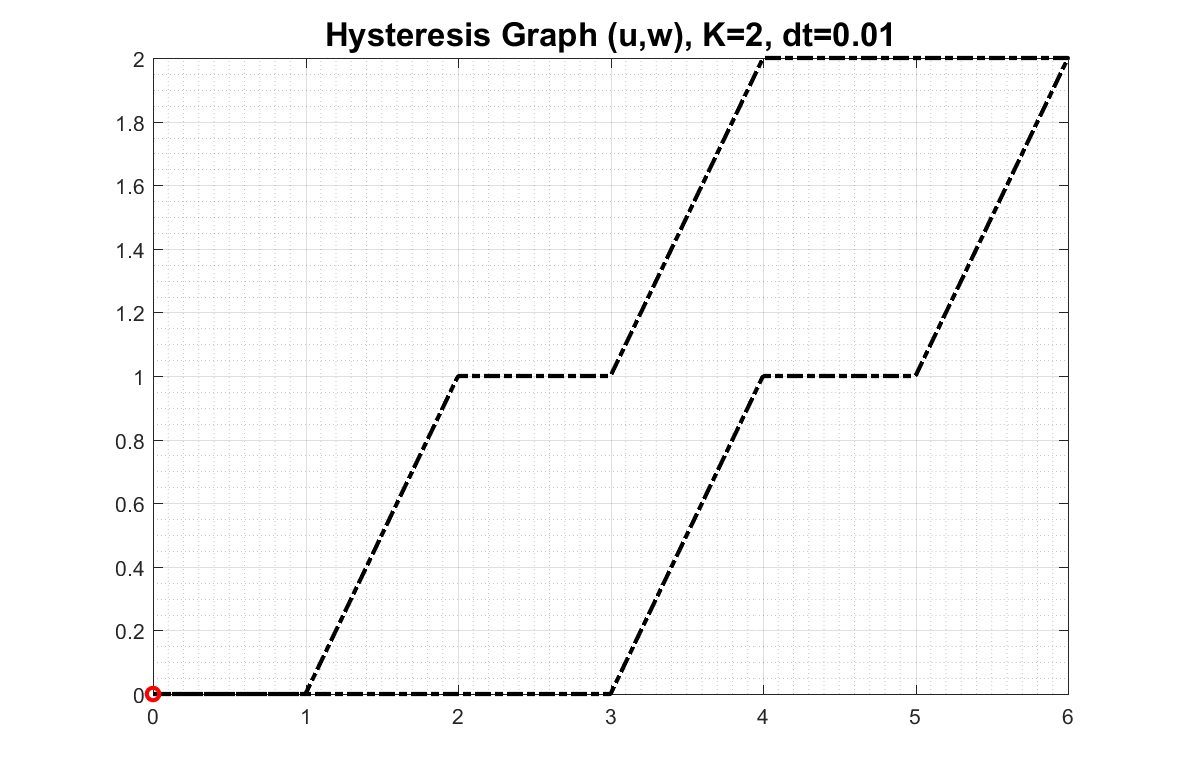}
&
\includegraphics
[height=30mm]{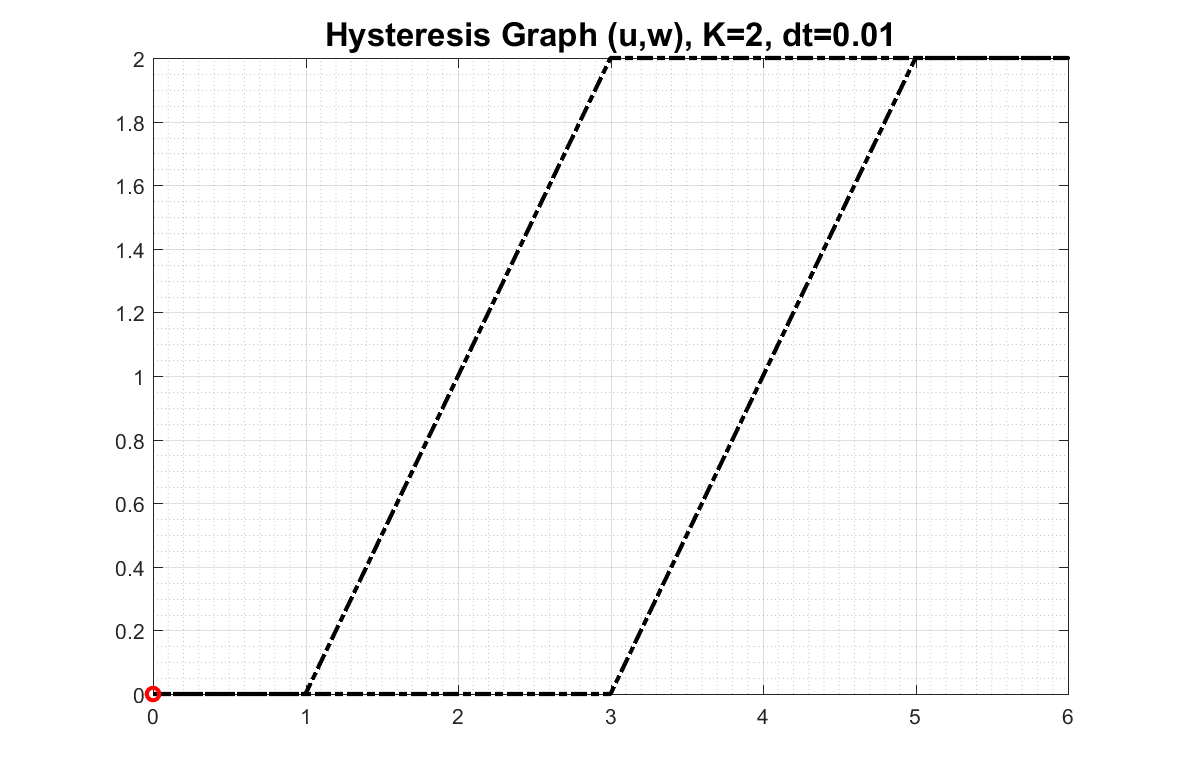}
&
\includegraphics
[height=30mm]{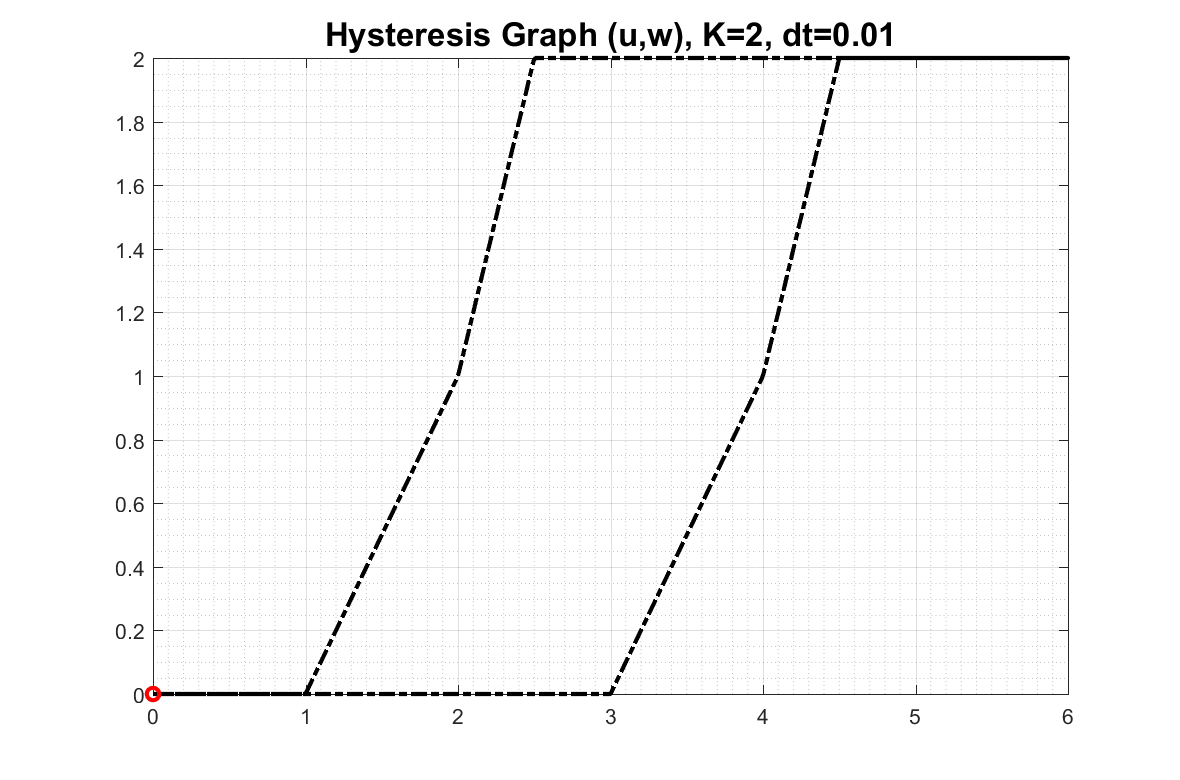}
\\
(a) $\pP%_{stackshifted} 
= [1,1,3,1;1,3,5,1]$
&
(b) $\pP%_{stacksum}
=[1,1,3,1;1,2,4,1]$
&
(c) $\pP%_{stack}
=[1,1,3,1;2,2,4,1/2]$
\\
\includegraphics
[height=30mm]{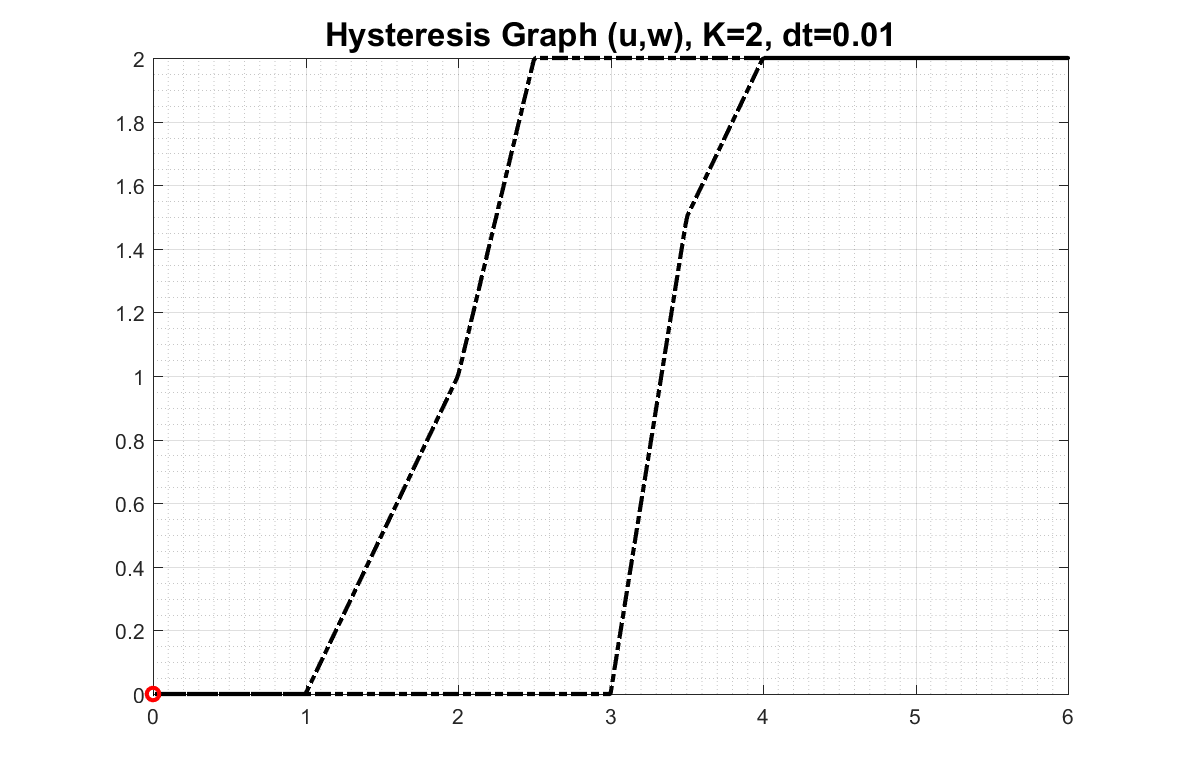}
&
\includegraphics
[height=30mm]{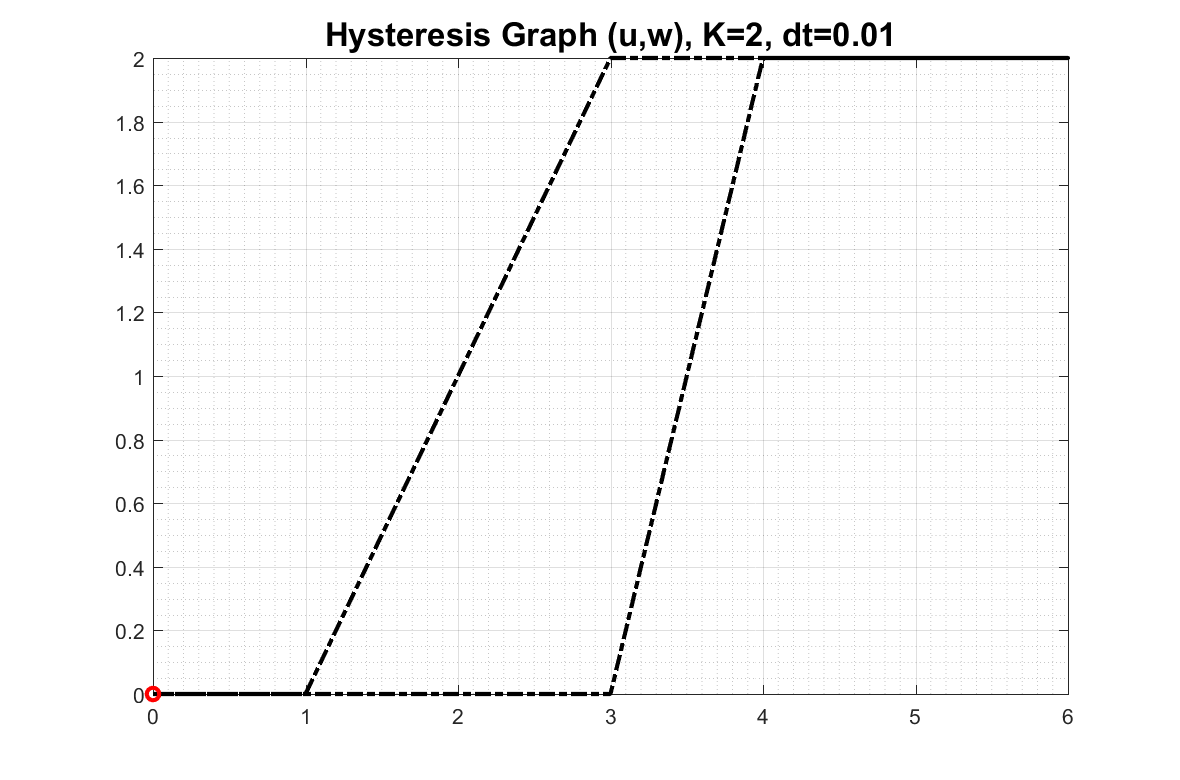}
&
\includegraphics
[height=30mm]{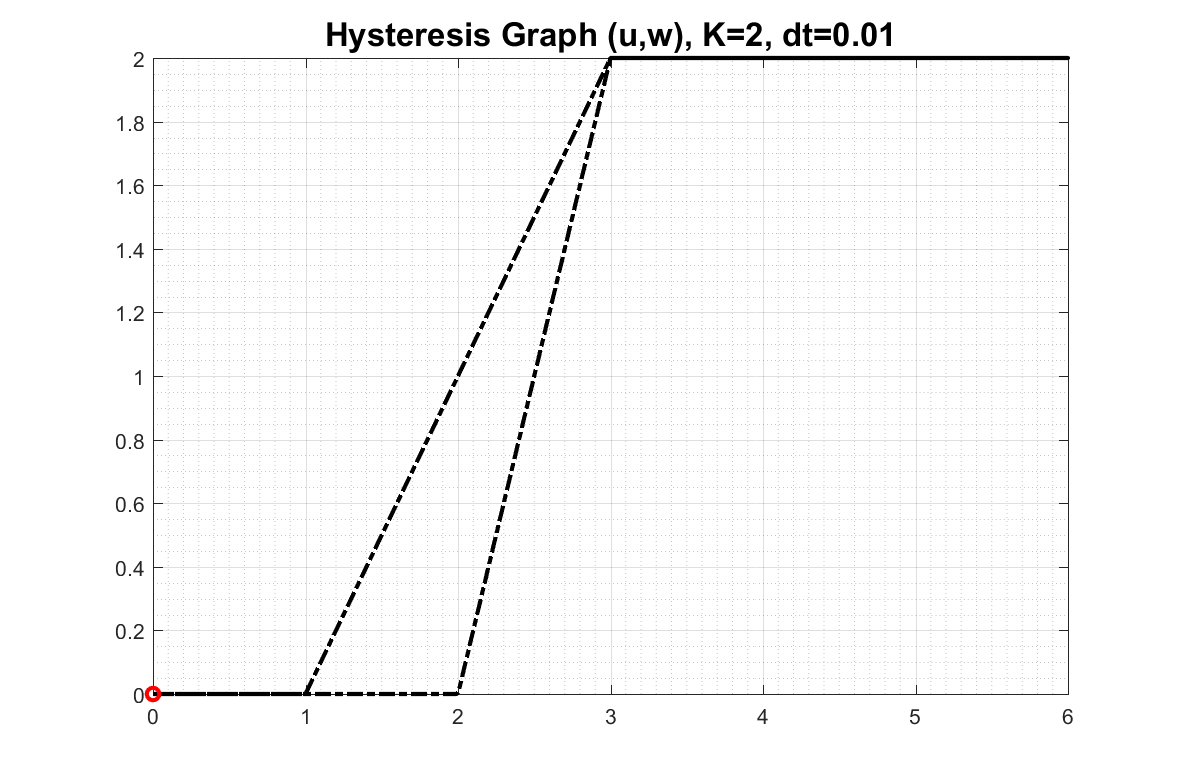}
\\
(d) $\pP%_{stacknonconvex} 
= [1,1,3,1;2,2,3,1/2]$
&
(e) $\pP%_{stacksquashed} 
= [1,1,3,1;1,2,3,1]$
&
(f) $\pP%_{stacktri}
= [1,1,2,1;1,2,2,1]$
\end{tabular}
\caption{Examples of stacking $K=2$ hysterons in \knonlinear\ model. Shown is $H=(u,w)=\hH(\pP;u)$ with  $\pP$ as indicated and $u(t)=\PL(0,6,0)$ in all cases.
\label{fig:stack}
}
\end{figure}

With $K=2$ one can easily write out the different possibilities; see  Fig.~\ref{fig:stack} for illustration. Recall $\alpha_k\leq \beta_k$, and denote
\bas 
A_k = \alpha_k + h_k,\ B_k = \beta_k + h_k.
\eas
We will say that two hysterons are adjacent on the left if $A_k=\alpha_{k+1}$, and on the right if  $B_k=\beta_{k+1}$,
They are coincident on the left if $\alpha_k=\alpha_{k+1}$, and on the right if $\beta_k=\beta_{k+1}$. 

We start by adding the hysterons 
$
\hH(\mu,\alpha_1,\beta_1,h_1;\cdot)+ \hH(\mu,\alpha_2,\beta_2,h_2;\cdot) 
$
when $\mu=1$. 

(a) If $A_1 < \alpha_2$, we obtain a flat section on the left bounding line. Likewise, if $B_1 < \beta_2$, there is a flat section on the right bounding line. See Fig.~\ref{fig:stack}~(a)
with $\pP=[1,1,3,1;1,3,5,1]$. 

(b) If $A_1 = \alpha_2$ and $B_1 = \beta_2$, then the flat sections are eliminated, i.e., the sections are { adjacent}, and this sum has the same bounding curves as the single hysteron $\hH(1,\alpha_1,\beta_1,h_1 + h_2;\cdot)$; see Fig.~\ref{fig:stack}~(b) with $\pP=[1,1,3,1;1,2,4,1]$. However, neither the operators nor the secondary curves match
\bas
\hH(1,\alpha_1,\beta_1,h_1;\cdot)+ \hH(1,A_1,B_1,h_2;\cdot) 
\neq
\hH(1,\alpha_1,\beta_1,h_1 + h_2;\cdot).
\eas

(c-d) Continuing with adjacent sections, the ranges match but the slopes do not for the two hysterons $\hH(1,2,4,1;\cdot) \neq \hH(2,2,4,1/2;\cdot)$; the bounding function on each side switches from slope $2$ to slope $1$ at the single node at which the sections are joined. See Fig.~\ref{fig:stack}~(c) with $\pP=[1,1,3,1;2,2,4,1/2]$. Another example is provided in (d) when $\pP=[1,1,3,1;2,2,3,1/2]$.

(e) Our most important example produces different slopes on the two sides of a single section, i.e., $H$ is a trapezoid. Towards this, we stack a pair of  hysterons for  which one end has adjacent sections but on the other end the sections are { coincident}. For example, we take 
$\alpha_1 + h_1 = \alpha_2$ and $\beta_1 = \beta_2,\ h_1 = h_2$. With $\mu_1=1,\mu_2=1$ we have slope 1 on the left and slope $2$ on the right. See Fig.~\ref{fig:stack}~(e)
with $\pP=[1,1,3,1;1,2,3,1]$.

\myskip{
\ba
\pP&=&
\begin{tabular}{|l|l|ll|l|}
\hline
$k$&$\mu_k$&$\alpha_k$&$\beta_k$&$h_k$\\
\hline
1&1&
1&3&1
\\
2&
1
&2&3&1
\\
\hline
\end{tabular}
= [1,1,3,1;1,2,3,1]
\ea
}

(f) Our final example shows what happens if in the case similar to (e) additionally, we have that $\beta_1=\alpha_1+h$. We obtain a degenerate trapezoid for which the two sides join at the top. 

%%%%%%%%%%%%%%%%%%%%%%%%%%
%%%%%%%%%%%%%%%

%%%
\subsubsection{Stacking $K$ hysterons for slopes of rational ratio}
\label{sec:ghysteron}
We continue case (e) from the previous example. Consider now some $m,n\in \mathbb{N}$ such that $r=\dfrac{m}{n}$ is irreducible. Consider 
$K=mn$ hysterons each with uniform $h_k=h,\mu_k=1$. Assume each of $\alpha_k$ and $\beta_k$ forms a non-decreasing set. The $\alpha_k$ are grouped as $m$ non-overlapping adjacent sets of $n$ points, and the $\beta$s are partitioned into $n$ non-overlapping adjacent sets of $m$ points:
\bsub
\label{eq:ghysteron}
\ba
\{\alpha_1=\dots=\alpha_n < \alpha_{n+1}=\dots=\alpha_{2n} < \dots < \alpha_{(m-1)n+1}=\dots=\alpha_{mn}\}
\\
\{\beta_1=\dots=\beta_m < \beta_{m+1}=\dots=\beta_{2m} < \dots < \beta_{(n-1)m +1}=\dots=\beta_{nm}\}
\ea
\esub
These give slope $n$ on $[\alpha_1,\alpha_1 + mh]$, and slope $m$ on $[\beta_1,\beta_1 + nh]$. The output $H$ has the shape of trapezoid, with the side slopes of ratio $r$. Multiplying all by some factor $\mu$ gives arbitrary slopes $s_l$ on left and $s_r$ on right, but their ratio $\dfrac{s_r}{s_l}=r$.  We call the resulting operator $\hH(\pP;u)$ a \ghysteron. 

\begin{definition}
\label{def:ghysteron}
$\hH(\pP;u)$ is called a {\ghysteron} if the primary curves in $H$ form a trapezoid whose top and bottom sides are parallel to the $u$-axis, and the left and right sides have positive slopes $s_l,s_r$ with a rational ratio 
\ba
\label{eq:ratass}
\frac{s_l}{s_r} \in \Q.
\ea
The parameters $\alpha_k,\beta_k$ satisfy \eqref{eq:ghysteron}, and $h_k,\mu_k$ are uniform. 
\end{definition} 

%%%%%%%%%
\subsubsection{Comparison of \kpreisach\ model with \knonlinear}

\begin{figure}[ht]
\begin{tabular}{ccc}
\includegraphics
[height=35mm]{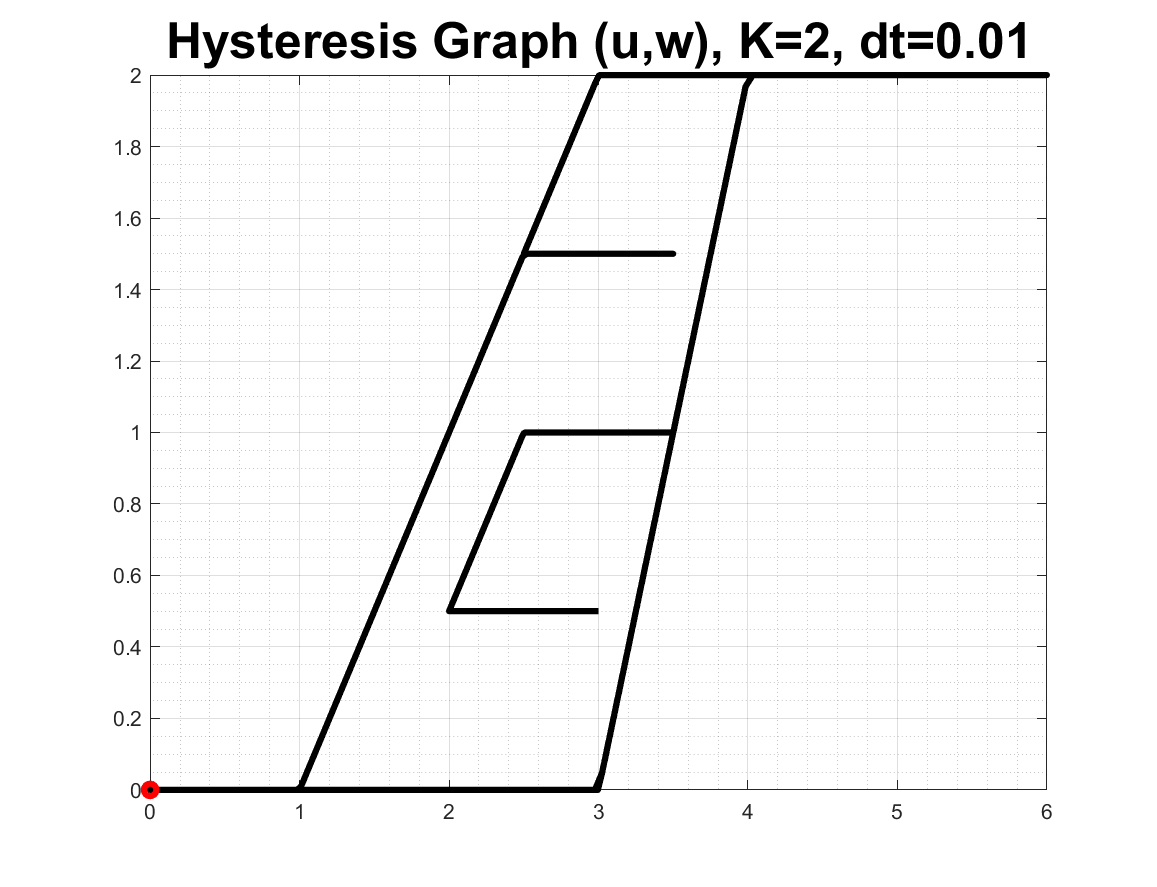}
&
\includegraphics
[height=35mm]{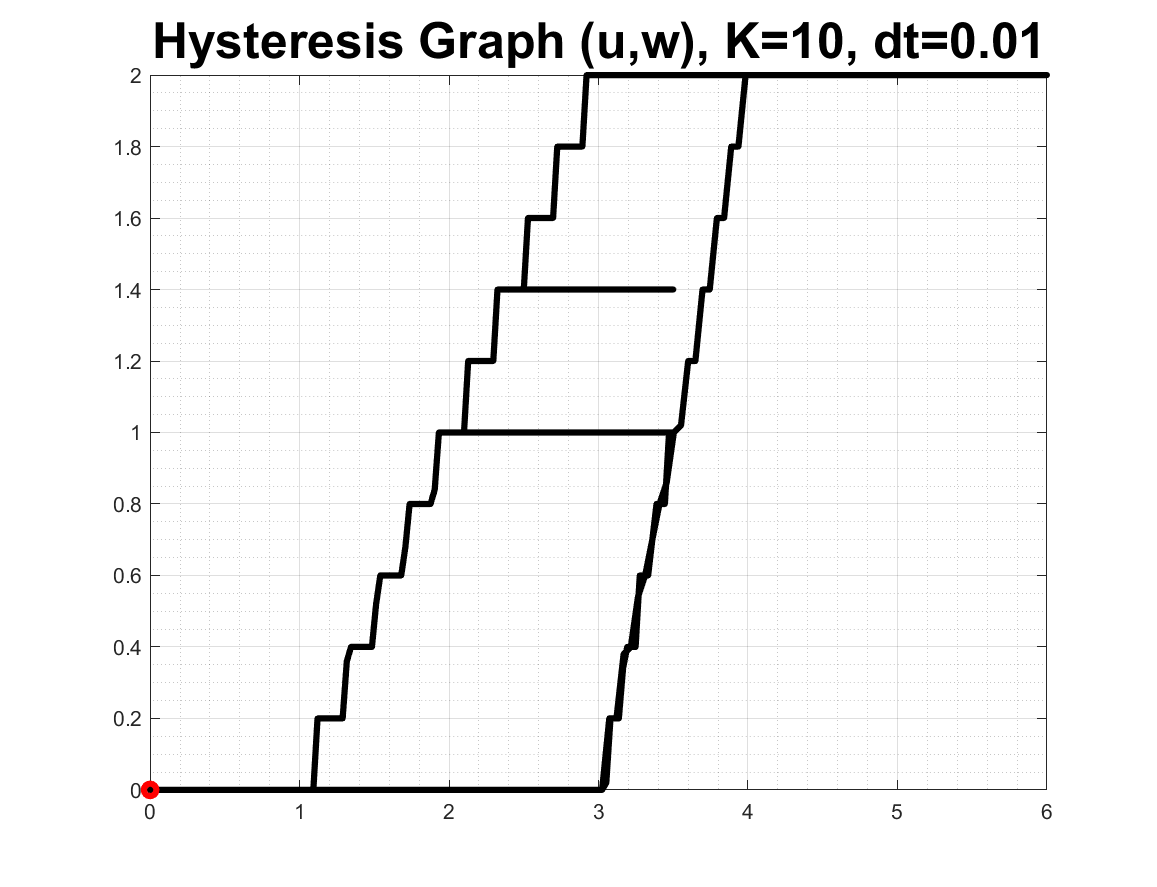}
&
\includegraphics
[height=35mm]{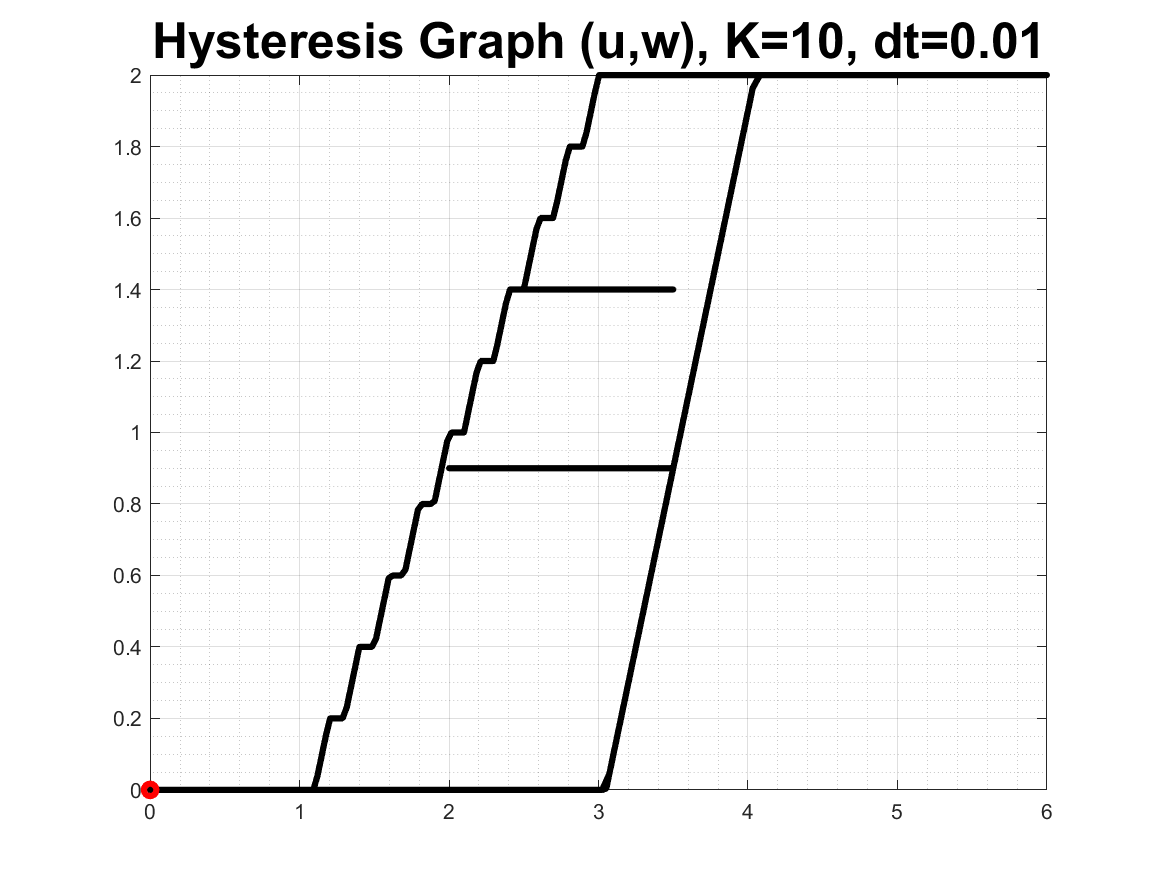}
%\includegraphics[height=30mm]{cimages/stacksquashed_preisach1_secondary.png}
%&
%\includegraphics[height=30mm]{cimages/stacksquashed_preisach_50_secondary.png}
\\
(a) $\pP^{2,1}$%_{stacksquashed} =
&
(b) $\pP^{10,0}$ 
&
(c) $\pP^{10,\eps}$ with $\eps=1/2$.
%&
%(d) $\pP^{10,\eps}$,  with $\eps=1$.
\end{tabular}
\caption{Comparison
of $\hH(\pP;u)$ with \knonlinear\ and  the $K$-Preisach model for the same \ghysteron\ $H$ as in Fig.~\ref{fig:stack}~(e). We use
$\uU=[0,6,2.5,3,3.5,0,3.5,2,3]$ which more than sweeps $H$. (a) \knonlinear with $\pP^{2,1}= [1,1,3,1;1,2,3,1]$. (b): $K$-Preisach approximation is rough with $K=10$, and (c) smoother if $\ppreisache$ is used. Note the richness of the secondary scanning curves in (a) and flatness in (b-c). 
\label{fig:stack-preisach}
}
\end{figure}

Finally we continue  Example~\ref{sec:ghysteron}(e), and  compare the resulting $\hH(\pknon;u)$ to that obtained with \kpreisach\ model; see Fig.~\ref{fig:stack-preisach}, and calibration with monotone $\{\alpha_k,\beta_k\}$. Since $K$ is finite, as expected, $\hH(\ppreisach;u)$  produces discontinuous output. With regularization with $\bheps$, the stair-step effect is diminished.  With $\eps=1$, the primary curves in $H$ for $\ppreisach$ are close to
those for $\pknon$, but the secondary curves differ substantially.

%%%%%%%%%%%%%

%\end{document}
%%%%%%%%%%%%%%%%%%%%%%%%%%%%%%%%%%%%%%%%%%%%%%%%%%%%%
\section{Calibration and approximation of hysteresis functionals}
\label{sec:inverse}
In this section we calibrate $\hH$, i.e., we consider the ``inverse problem'': Given some data $\dD$ including $H$, find $\pP^*$ so that %%
\ba
\label{eq:happrox}
H \approx H^* =\hH(\pP^*;u;\vz).
\ea
for any $u$ that sweeps $H$.
We recall that the calibration process provided in literature for the Preisach model involves the so-called Preisach plane and requires  $(u,w) \in \dD$ with $u$ from a family of inputs dense in $C^0$. See, e.g., the contributions in \cite{HSV88,Krejci13,Kadar88}. 

Finding $\pklin$ is easy but only possible if $H$ is symmetric as indicated in Tab.~\ref{tab:kmodels}. We show how to find $\pgen,\ppreisach,\pklin,\pknon$ for the same $H$. The algorithm for \knonlinear\ model is most involved, even if $H$ is a {\ghysteron}. 

We assume below that $H$ is {\adm}, i.e., that  $H$ is a boundary of some $\hH \subset \R^2=\{(u,w)\}$, a generalized trapezoid defined below.

%%%%%%%%%%%%%
\begin{definition}
\label{def:repr} $H$ is a generalized trapezoid $H$ if it is a boundary of closed  simply connected region $\hH=\{(u,w)\} \subset \R^2$ satisfying the following. $H$ must have top and bottom sides parallel to the $u$ axis, and  monotone lateral sides so that its boundary $H=\bd(\hH)= \hbot \cup \hleft \cup \htop \cup \hright$ of $\hH$ is made of, respectively, the bottom, left, top, and right bounding curves. We have
\bas
\hbot = [\alpha,\beta] \times \{w_{min}\}; \;\; \htop = [A,B] \times \{w_{max}\}
\eas
for given $\alpha,\beta,A,B,w_{min},w_{max} \in \R$ such that
\ba
\label{eq:graph}
0\leq \alpha \leq \beta, \;\; 0 \leq A\leq B, \; \alpha <A, \beta < B, \; 0\leq w_{min} < w_{max},\; 0< \frac{A-\alpha}{B-\beta}<\infty.
\ea
The left and right curves $\hleft=(u,\fl(u)),u\in [\alpha,A]$ and $\hright=(u,\fr(u)),u \in [\beta,B]$, are the graphs of functions which are continuous piecewise smooth increasing and injective on $(\alpha,B)$.
These functions either (i) coincide on all $[\alpha,B]$, or (ii) they satisfy
%%%%%
\ba
\label{eq:flu}
\fl(u) > \fr(u), \; u \in (\alpha,B); \; \mathrm{and}\;\; \fl(u)\geq \fr(u), \; u\in \{\alpha,B\}.
\ea
%%%
The curves $\hleft, \hright$ are the downward left curve and the upward right curve, respectively. 
For convenience we also list the vertices $\ver(H)$ of $H$ in the counter-clockwise order
\ba
\label{eq:V}
\ver(H) = \{(\alpha,w_{min}), \; (\beta, w_{min}), \; (B,w_{max}),\; (A,w_{max})\}.
\ea
Finally, the data on $H$ might be given from experiment, i.e., 
\ba
\label{eq:dd}
\dD=(u^l_k,w_k)_{k=1}^{K+1} \cup (u^r_k,w_k)_{k=1}^{K+1} \ea
chosen so that $w_k=\fr(u_k^r)=\fl(u_k^l)$, with $u_1^l=\alpha, u^l_{K+1}=A, u^r_1=\beta,u^r_{K+1}=B$. Assume that $(u^r_k)_k$ and $w_k$ are increasing sequences.
\end{definition}

We note that some graphs which are not \adm\ can be broken up into smaller pieces which are amenable to approximation. Further, if $H$ is non-hysteretic, i.e., $\alpha=\beta$, and $A=B$, it can be parametrized by a single unit hysteron. Lastly, if instead of \eqref{eq:dd}, the data $
\dD=\dD^l \cup \dD^r = (u^l_k,w^l_k)_{k=1}^{K_l+1} \cup (u^r_k,w^r_k)_{k=1}^{K_r+1} 
$ with $K_l \neq K_r$ or $w^l_k\neq w^r_k$, then one must pre-process $\dD$, e.g., by taking an intersection of $\dD^l$ and $\dD^r$ and interpolating. 

\medskip
Now we comment on the inputs $u(t)$.
We say that the input $u(t)$ {sweeps} $H$ if $u:[0,T]=[0,T^+]\cup [T^+,T]\to [\alpha,B]$ is an absolutely continuous function $u \in W^{1,1}[0,T)$ such that 
\bsub
\label{eq:usweep}
\ba
u(0)=u(T)=\alpha; u(T^+)=B; 
\\
u'(t) \geq 0, \; \mathrm{a.e.}\;\; t \in [0,T^+]; \;\;  u'(t) \leq 0, \; \mathrm{a.e.} \; t \in [T^+,T].
\ea
\esub
We denote the set of such sweeping functions by $u \in \sS(H)$. Typically we choose $u\in \PL(\uU)$ for some $\uU$ that includes $\{\alpha,B,\alpha\}$.  

%%%%%%%%%%%%%%%%%%
\subsection{Finding $\pgen$ for \general} 
Calibration of $\pgen$ does not take any effort. We take as $\fl$ and $\fr$ the functions whose graphs form $\hleft$ and $\hright$.  We record $\pgen=[1,\fl(\cdot),\fr(\cdot),\id]$. %
If discrete experimental data is used, we can set $\fl(\cdot), \fr(\cdot)$ to be, e.g., piecewise linear interpolants of the data on $\hleft,\hright$, respectively. 

\subsection{Finding $\ppreisach, \ppreisache,\ppreisachs$  for \kpreisach\ models}
\label{sec:kpreisach-i}
Given $K$,  we partition the range $[w_{min},w_{max}]$ into $K$ intervals, and form $K$ rectangles, each of height $h_k$ so that $\sum_k h_k=w_{max}-w_{min}$. It is easiest to choose uniform $h_k=h= (w_{max}-w_{min})/K$. Then we set $w_k=w_{min}+\sum_{j=1}^{k-1}h_j, k=1,\ldots K+1$. For each $w_k$ we find $u^l_k=\fl^{-1}(w_k)$ and $u^r_k=\fr^{-1}(w_k)$.
We see that $u^l_1=\alpha$, and $u^l_{K+1}=A$, while 
 $u^r_1=\beta$, and $u^r_{K+1}=B$. 
Finally we set rectangles, each of height $h_k$, with left corner $\alpha_k=(u^l_k+u^l_{k+1})/2$ and $\beta_k=(u^r_k+u^r_{k+1})/2$. We record each $\pP_k=[\infty,\alpha_k,\beta_k,h_k]$, and collect in $\ppreisach$. The output $u \to w=\hH(\ppreisach;u)$  will be discontinuous.

We can now choose some $\eps$ for  
$\ppreisache$, with each $\pP_k$ replaced by $\pP^{\eps}_k=[1/\eps,\alpha_k,\beta_k,\eps h_k]$, with some small $\eps$. The output $u \to w=\hH(\ppreisache;u)$  will be continuous with intermittent flat pieces. 

%%%%
\subsection{Finding $\pklin$ for \klinear}
For a \klinear\ operator $u \to \hH(\pP;u)$ finding $\pP$ requires only the knowledge of one of $\hright$ or $\hleft$, since $\hleft$ must be a symmetric reflection of $\hright$ with respect to the midpoint of $H$.  For meaningful calibration we assume that the top and bottom parts of $H$ are single points where $\hleft$ and $\hright$ intersect. 
To calibrate, wlog, we take $H_r$.  For every interval $(u^r_k,u^r_{k+1})$ we approximate the slope of $\fr(\cdot)\vert_{(u^r_k,u^r_{k+1})}$ with finite differences $s_k = \frac{w_{k+1}-w_k}{u^r_{k+1}-u^r_k}$. Next we set simply
\bas
\alpha_k = u_1; \forall k;\;
\beta_k = u_k; \forall k; \;
\mu_1=s_1;\; \mu_k=s_k-\sum_{m=1}^{k-1}\mu_m; k>1. 
\eas
Since the weights $\mu_k$ approximate the second derivatives of $\fr(\cdot)$, which is convex, we get $\mu_k\geq 0$; see 
\cite{PeszShow98} for more. 
The secondary curves of \klinear\ $\hH$ are not horizontal; rather, they have shape similar to that of translates of  $\fl,\fr$, with a rich structure completely determined by $\pklin$. This will be evident in examples in Sec.~\ref{sec:convex}. 

%% Use code such as HysteresisOutput(10,0.01,0,construct_klinear(uu,wu),[uu(1),uu(end),uu(1),-22,-45,-25]);

%HysteresisOutput(10,0.01,0,construct_klinear(uu,wu),[uu(1),uu(end),uu(1),4,2,4,4.5,2,4]);

%HysteresisGraph;HysteresisOutput(10,0.01,0,construct_klinear,[1,4,1,2,3.5,1.7,3]);

%uu=linspace(1,5,50);wu=(uu-1).^2+1;

%%%%%%%%%
\subsection{Finding $\pknon$ for \knonlinear}
\label{sec:knon} 
We start in Sec.~\ref{sec:trapezoid} by describing how to parametrize $H$ when $H$ is a \ghysteron\  as in  Def.~\ref{def:ghysteron}. We follow in Sec.~\ref{sec:multilevel} by a hierarchical algorithm which approximates any \adm\ $H$ (as in Def.~\ref{def:repr}) with  curvilinear sides, by a sum of {\ghysteron}s.

%%%%%%%%%%%%
%%%%%%%%%
\subsubsection{Algorithm $H \approx H^* \to \pP$ for \ghysteron\ $H$ with linear sides} \label{sec:trapezoid}
Assume a trapezoid $H$  as in Def.~\ref{def:ghysteron}, with vertices
$
\ver(H) = \{(\alpha,\underline{w}), \; (\beta, \underline{w}), \; (B,\overline{w}),\; (A,\overline{w})\}.
$
We calculate
\bas
h=\overline{w}-\underline{w}; \; \fl'=s_l=\frac{h}{A-\alpha}; \fr'=\frac{h}{B-\beta}; \;
r=\frac{s_r}{s_l}=\frac{A-\alpha}{B-\beta}.
\eas
We can set $K=mn$ and  calculate $\pP=(\pP_k)_{k=1}^K \in \R^{K \times 4}$ with simple formulas given below.
However, we can also approximate $\Q \ni r^*=\tfrac{m^*}{n^*} \approx r$, with a new $K^*=m^*n^* \ll K$. This is useful when $r \not \in Q$, i..e, \eqref{eq:ratass} does not hold, or when  $K=mn$ is impractically large.  The approximation $r \approx r^*$ now gives a new $H^* \approx H$ with $\ver(H) \approx \ver(H^*)$. The parameters $\pP^*$ for $H^*$ are found with the calculations below.
 
\noindent
(STEP (a)) Approximate $r$ by an irreducible fraction $r \approx r^*=\frac{m^*}{n^*} \in \Q. \; \mathrm{Set\ } K^*=m^*n^*$.
%%
%\bsub

\noindent 
(STEP (b)) Choose the new slopes $s_l^* \approx s_l,s_r^* \approx s_r$ of the sides of $H^*$. The choice of $s_l^*$ and $s_r^*$ is not unique, but we require $r^*=\frac{s_r^*}{s_l^*}$.  Once $s_l^*$ and $s_r^*$ are set,  calculate the points $A^*,B^*$, the scaling factor $\mu^*$ and the subinterval length $h^*$, 
$
A^*=\alpha+\frac{h}{s_l^*};  B^*=\beta+\frac{h}{s_r^*}; 
\mu^* =s_l^* m^*=s_r^* n^*;  h^*=\frac{h}{\mu^*}.
$
For example, we can set $s_l=s_l^*$, and use $s_r^*=s_l r^*$,  $\mu^*=s_l{m^*}$, and $A^*=A$ and $B^*=\beta+ \frac{h}{s_r^*}$. Alternatively, we can set $s_r^*=s_r$, $s_l^*=\frac{s_r^*}{r^*}$, with
$\mu^*={s_r}{n^*}$, $B^*=B$ and $A^*=\alpha+ \frac{h}{s_l^*}$.  Other options are possible. 

\noindent
(STEP (c)) With $\mu^*$, and $h^*$ known, set
$
\mu_k=\frac{\mu^*}{K^*}; \;\; 
h_k=h^*; \;\; k=1,\ldots K^*.
$

\noindent
(STEP (d))
Define $\alpha_k, \beta_k$, enumerating $1 \leq k \leq K^*$ as $k=(j-1)n+l$ or $k=(l-1)m+j$.  
\bas
\alpha_{(j-1)n+l}=\alpha+(j-1)h^*; \;
\beta_{(l-1)m+j}=\beta+(l-1)h^*; \; 1 \leq j \leq m,\; \; 1\leq l \leq n.
\eas
%%
%%%

\medskip

In the end we collect $\pP^*$ with $H^*=\hH(\pP^*;u), u \in \sS(H)$. In particular, we require 
\ba
\label{eq:Vs}
\ver(\hH^*) = \{(\alpha,\underline{w}), \; (\beta, \underline{w}), \; (B^*,\overline{w}),\; (A^*,\overline{w})\},
\ea
satisfy \eqref{eq:graph}, and that $\abs{A-A^*}+\abs{B-B^*} \approx 0$ and $\hleft \approx \hleft^*$ and $\hright \approx \hright^*$.  

\begin{remark}
\label{rem:rat}
Given $r \in \R$, the choice of irreducible fraction $r^*\approx r$ is not unique. To maintain accuracy, i.e., to minimize $\abs{r-r^*}$, we use {\em diophantine} approximations, e.g., with the MATLB function {\tt{rat}} with desired accuracy \cite{DIOPH}. To constrain the magnitude of $K^* = m^*n^*$, we consider some $K^{max}$ and search for $m^*,n^*$ 
\ba
(m^*,n^*)=\mathrm{argmin\ }_{(m,n): mn \leq K^{max}}  \left|\frac{m}{n}-r\right|^2.
\ea
In practice, it suffices to search in $[1,\sqrt{K^{max}}]^2$, with $K^{max}=O(100)$.  
\end{remark}

%%%%%%%%%%%%
\begin{figure}[ht]
\begin{tabular}{ccc}
\includegraphics
[height=40mm]{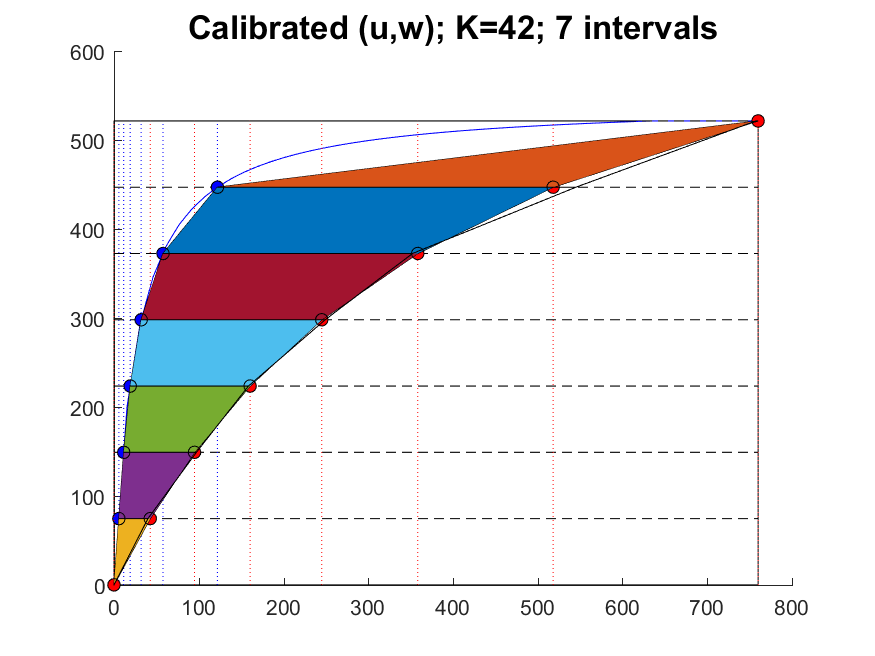}
&
\includegraphics
[height=40mm]{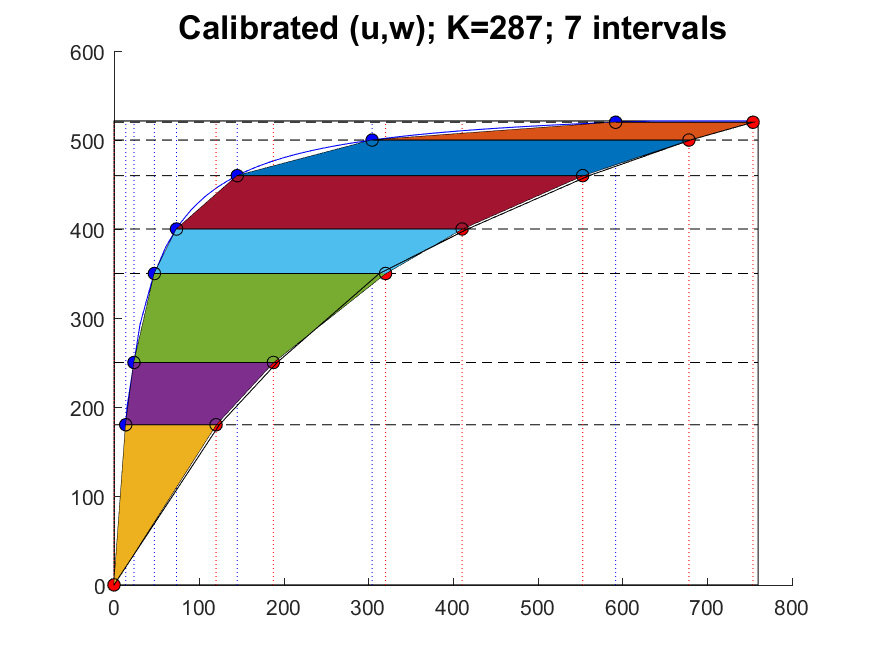}
&
\includegraphics
[height=40mm]{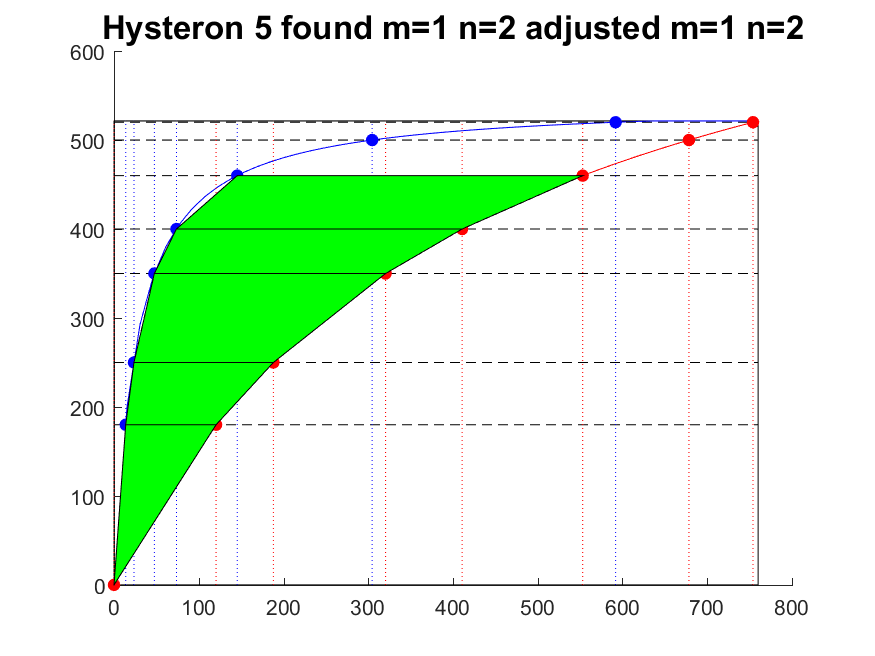}
\\
(a)&(b)&(c)
\end{tabular}
\caption{\label{fig:calibrate}
Calibration for adsorption hysteresis graph, with $H$ made of $\fl$ curve (in blue) and $\fr$ (in red), with original $H$ as shown in Fig.~\ref{fig:teaser} (a).  Partition of the range into $I=7$ intervals is the first step in the calibration with $\pknon$, with the corresponding {\ghysteron}s filled in. (a) With a uniform partition, some portions of $\hH$ are inadequately covered. (b) An adaptively chosen partition gives a better approximation $H^* \approx H$ and better coverage of $\hH$. The \knonlinear\ graph is found for (a) with $K=42$ and (b) $K=287$, respectively. (c) For the adaptive vertical partition, the individual {\ghysteron}s are found step by step. The step $i=6$ (with $K_6=2$) is shown, with the {\ghysteron}s $\hH_1,\hH_2\ldots \hH_{i-1}$ filled with green.}
\end{figure}

%%%%%%%%%%%
\subsubsection{Hierarchical approach for  $H$ with general shape in Def.~\ref{def:repr}} 
\label{sec:multilevel}
We approximate $H\approx H^*$ and cover $\hH$ by a union of $I$ {\ghysteron}s with a process illustrated in Fig.~\ref{fig:calibrate}. 
The key is the partition  
$
w_{min}=w_0 < w_1 < \ldots w_I=w_{max}.
$
of the range $[w_{min},w_{max}]$ of $\hH$ into $I$ subintervals which can be, in principle, arbitrary.  We isolate the partitions 
\bas
\hH_i = \hH \cap \{(u,w) \in \hH: w_{i-1} \leq w \leq w_{i} \}
\eas
of $\hH$. The boundary $H_i$ of each graph $\hH_i$  has the top and bottom sides parallel to the $u$-axis, with the height at most $h_I=\max_i (w_i-w_{i-1})$. The curvilinear left and right sides of $H_i$ follow the curves $(u,\fl(u))$ and $(u,\fr(u))$, respectively, and the vertices of $H_i$ 
\bas
%\label{eq:vertices}
\ver(H_i) = \{(\alpha_i,w_{i-1}),(\beta_i,w_{i-1}),(B_i,w_i),(A_i,w_i)\}; \; 
\alpha_{i}=(\fl)^{-1}(w_{i-1}), 
\beta_{i}=(\fr)^{-1}(w_{i-1})
\eas
are well defined. The continuity of $\fr,\fl$ requires
\ba
\label{eq:coincide}
A_i=\alpha_{i+1}, \;\;  B_i=\beta_{i+1}, \;\; 0\leq i \leq I-1.
\ea
If the $i'$th portion of the graph is not hysteretic, we have $\alpha_i=\beta_i$, and $A_i=B_i$. 
If the sides of each $H_i$ are linear and satisfy \eqref{eq:ratass}, we find some $\pP_i=(\pP_{k_i})_{k_i=1}^{K_i}$ for each $H_i$, and collect $\pP=((\pP_{k_i})_{k_i=1}^{K_i})_{i=1}^I$ renumbered appropriately. 
However, if the sides are not linear,  or if \eqref{eq:ratass} does not hold, or if $\sum_i K_i$ is too large, we proceed by iteration to satisfy the accuracy and efficiency needs, while we maintain continuity as in \eqref{eq:coincide}. An example of such iterative algorithm is given in the Appendix~\ref{sec:iteration}.

%

%%%%%%%%%%%%
\subsection{Quality of approximation $H \approx H^*$}
The efforts to approximate a given $H\approx H^*$ are not much different from those of piecewise interpolation of the sides $\hleft,\hright$ of $H$. Once $\pP^*$ is calibrated, the model $u \to w \in \hH(\pP^*;u)$  can be coupled to some external dynamics, e.g., to ODE or PDE. However, then an additional modeling error arises since the actual accuracy of output $(u(t),w(t)) \in \hH$ depends on the accuracy of approximation of $u(t)$ and $w(t)$ by some $(U_{\tau},W_{\tau})$ in these other equations coupled to $u \to \hH(\pP;u)$. 

With a large number $K$ of components $v_k$, the error in $u$ and $w$ would seem to accumulate in $w-W_{\tau}$ from that for the individual components $v_k-V_{k,\tau}$, and would affect $u-U_{\tau}$. However, for all models except \klinear, the effect of the accumulation seems insignificant in practice, and a large $K$ is not an issue for accuracy of time-stepping, but may be desired to reduce model error $H-H^*$. At the same time, large $K$ requires more computational time. 

Finally when discussing $H \approx H^*$, we must realize that in practice we encounter ($U^*_{\tau}, W^*_{\tau}) \in \hH_{\tau}(\pP^*;u^*_{\tau})$. 
 The modeling error $H\approx H_{\tau}^*$ contributes  to the global approximation error  $u-U_{\tau}^*$, and $w-W_{\tau}^*$.
%%%%%%%%%

\subsection{Examples of \knonlinear, \kpreisach\ and \klinear\ models} 

%%%%%%%%%%
\subsubsection{Trapezoid with curvilinear sides}
\label{ex:bottom}
Consider first $H^{lin}$ with $\ver(H^{lin}) =\{(3,0),(9,0),(11,5),(4,5)\}$ and linear sides. The slopes of these sides are $s_l=5$, $s_r=\frac{5}{2}$, with $h=5$, and $r=\frac{s_l}{s_r}=\frac{1}{2}$, We accept $H^*=H^{lin}$, and set $m^*=m=1,n^*=n=2$, with $K^*=2$, and calculate $h^*=\frac{A-\alpha}{m^*}=1$, $\mu^*=\frac{h}{h^*}=5$. Each $\mu_k=\frac{\mu^*}{K^*}=\frac{5}{2}$, and $h_k=h^*=1$. We summarize $\hH(\pP^{lin};\cdot)$ with
$
\pP^{lin}=[\frac{5}{2},3,9,1;\frac{5}{2},3,10,1]$.

%\label{ex:curvy}
Consider now curvilinear hysteresis graph $\hH$ with  $\ver(H)=\{(3,0),(9,0),(B,5),(4.1,5)\}$, and $B=11+\tfrac{\pi}{10}$. The left side of $\hH$ is given by the quadratic polynomial
$\fl(u)=\frac{50}{11}(x-3)+\frac{1}{2}[(x-3.55)^2-0.55^2]$, while $\fr(u)$ is the piecewise linear function which connects the vertices $(9,0)$ and $(B,5)$.   To find an approximation $\hH^*$, we proceed by iteration. We set $q=0$ and
$
\ver(H^{(0)})=\{(3,0),(9,0),(B,5),(4.1,5)\}.
$
Here the linear left side of $H^{(0)}$  with the slope $s_l=\frac{50}{11}$ does not match very well the primary scanning curve in $H$. We also see that $s_r$ in exact precision is not rational, thus we continue. 

In iteration $q=1$ we consider a {\ghysteron} $\hH^{(1)}$  associated with
\bas
\ver(H^{(1)})=\{(3,0),(9,0),(B^{(1)},5),(4.1,5)\},
\eas
with the double prevision decimal approximation $B^{(1)} =11.314159265358979 \approx B$. 
We have $r^{(1)} \approx \tfrac{11}{23.14159265358979}$. We can find $m^{(1)},n^{(1)}$ so that $r^{(1)} \in \Q$. However, the corresponding $K^{(1)}=O(10^{16})$, very large and impractical.  
We try next $\ver(H^{(2)})$, setting $B^{(2)}=11.3$, with
%%%%%%
$
\ver(H^{(2)})=\{(3,0),(9,0),(11.3,5),(4.1,5)\}.
$
%%%%
Now  $r^{(2)}=\frac{11}{23}$.
We set $m^{(2)}=11$, $n^{(2)}=23$, and parametrize $H^{(2)}$ with $K^{(2)}=mn=253$ {\uhysteron}s. In particular, we have $h^*=\frac{1.1}{11}=\frac{1}{10}$, and $\mu^*=50$, while each $\mu_k=\frac{50}{11 \cdot 23}$, and each $h_k=h^*$. We number the {\uhysteron}s with $k=(j-1)n+l$ or $k=(l-1)m+j$ when $1\leq j\leq 11,1\leq l\leq 23$. In particular, $\alpha_1=\ldots \alpha_{23} =\alpha=3$, but $\alpha_{24}=3+h^*$, and $\alpha_{253}=\alpha + 10h^*=4$.  In turn, $\beta_1=\beta_{24}=\beta=9$ but $\beta_2=\beta+h^*$, and $\beta_{253}=\beta + 22h^*=9+2.2=11.2$.  
We get
\ba
\pP^{(2)}=
\begin{tabular}{|l|ll|ll|}
\hline
$k$&$\alpha_k$&$\beta_k$&$h_k$&$\mu_k$\\
\hline
1&3&9&1&$\tfrac{50}{253}$
\\
2&3&9.1&1&$\tfrac{50}{253}$
\\
\ldots&&&&\\
24&3.1&9&1&$\tfrac{50}{253}$
\\
\ldots&&&&\\
253&4&11.2&1&$\tfrac{50}{253}$
\\
\hline
\end{tabular}
\ea
However, $K^{(2)}$ may be still too large to be practical, and we try again. Setting $B^{(3)}=11.2$ we
have $\frac{s_r}{s_l}=\frac{1}{2}$ and we can set $m^*=1,n^*=2$. The parametrization of 
$
\ver(H^{(3)})=\{(3,0),(9,0),(11.2,5),(4.1,5)\}.
$
is similar to that for $H^{lin}$.  The difference is the new scaling factor $\mu^*=\frac{50}{11}$, and $h^*=\frac{11}{10}$. The individual scaling factor for each {\uhysteron} is now $\mu_k=\frac{25}{11}$. We set $\hH^*(u)=\hH(\pP^{(3)};u)$, with 
\ba
\pP^{(3)}=
\begin{tabular}{|l|ll|ll|}
\hline
$k$&$\alpha_k$&$\beta_k$&$h_k$&$\mu_k$\\
\hline
1&3&9&1.1&$\frac{25}{11}$
\\
2&3&10.1&1.1&$\frac{25}{11}$
\\
\hline
\end{tabular}
\ea

%%%%%%%%%%%
\myskip{
\subsubsection{Stack of two}
Now consider $\hH$ corresponding to the two {\ghysteron}s with bounding curves $H_{bottom}$ and $H_{top}$ stacked on top of one another
\myskip{
\bas
\ver(H_{bottom})&=&\{(3,0),(9,0),(11,5),(4,5)\},
\\
\ver(H_{top})&=&\{(4,5),(11,5),(15,10),(12,10)\}.
\eas
}
We recognize $H_{bottom}$ as that from Example~\ref{ex:bottom}. A quick check on $H_{top}$ reveals that it corresponds to a {\ghysteron}, with the slopes $s_l=\tfrac{5}{8}$ and $s_r=\tfrac{5}{4}$. Thus $m=2,n=1$ works naturally, and we amend the parametrization in Ex.~\ref{ex:bottom} as follows  for $\hH((\pP;\cdot)$ with 
%%%%%%%%%%%%
%\bsub
\ba
\myskip{= \frac{5}{2}\left(
b^1(\gamma_{3,9}(u))+b^1(\gamma_{3,10}(u))
\right)
+ \frac{5}{8} \left(
b^4(\gamma_{4,11}(u))+b^4(\gamma_{8,11}(u))
\right)
.
\ea
%%%
Enumerating the {\uhysteron}s together, we have 
\myskip{
\begin{center}
\begin{tabular}{|l|ll|ll|}
\hline
$k$&$\alpha_k$&$\beta_k$&$h_k$&$\mu_k$\\
\hline
1&3&9&1&$\frac{5}{2}$
\\
2&3&10&1&$\frac{5}{2}$
\\
3&4&11&4&$\frac{5}{8}$
\\
4&8&11&4&$\frac{5}{8}$
\\
\hline
\end{tabular}
\end{center}
}
\ba
}
\pP=
\begin{tabular}{|l|ll|ll||l|}
\hline
$k$&$\alpha_k$&$\beta_k$&$h_k$&$\mu_k$
&
Vertices $\ver(H_i)$ for $i=bottom,top$\\
\hline
1&3&9&1&$\frac{5}{2}$
&
$\ver(H_{bottom})=\{(3,0),(9,0),(11,5),(4,5)\}$
\\
2&3&10&1&$\frac{5}{2}$
&
\\
\hline
3&4&11&4&$\frac{5}{8}$
&
$\ver(H_{top})=\{(4,5),(11,5),(15,10),(12,10)\}$
\\
4&8&11&4&$\frac{5}{8}$
&
\\
\hline
\end{tabular}
\ea
%\esub
}
%%%%%%%%%%%
\myskip{
\subsubsection{Triangles and non-hysteretic curves}

Next three examples show how use our algorithm for the portions of the graph $\hH$ which coincide at some point. For simplicity, they have linear sides. 
\myskip{
\bas
\ver(H_{line})&=&\{(0,0),(0,0),(1,2),(1,2)\},
\\
\ver(H_{tri-down})&=&\{(0,0),(0,0),(9,1),(3,1)\},
\\
\ver(H_{tri-up})&=&\{(12,10),(15,10),(18,13),(18,13)\}.
\eas
Consider first $\ver(H_{line})$. We see that $h=2$, $s_l=2=s_r$, thus $m=1=n=K$, $\mu=s_l m=2$, $h^*=\frac{h}{2}=1=h_1$, and $\mu_1=2$, with $\alpha_1=0=\beta_1$. 

Consider next $\ver(H_{tri-down})$. With $h=1$, we have $r=\tfrac{1}{3}$, and so $m=1,n=3=K$. Immediately $\mu^*=\tfrac{1}{3}$, $h^*=3=h_k$ for all $k$, and $\mu_k=\tfrac{1}{9}$. 
}
\begin{center}
\begin{tabular}{|l|ll|ll||l|}
\hline
$k$&$\alpha_k$&$\beta_k$&$h_k$&$\mu_k$
&
Vertices $\ver(H_i)$ for each portion $\hH_i$\\
\hline
1&0&0&2&$\frac{5}{2}$
&
$\ver(H_{line})=\{(0,0),(0,0),(1,2),(1,2)\}$
\\
\hline\hline
1&0&0&3&$\frac{1}{9}$
&
$\ver(H_{tri-down})=\{(0,0),(0,0),(9,1),(3,1)\}$
\\
2&0&3&3&$\frac{1}{9}$
&
\\
3&0&6&3&$\frac{1}{9}$
&
\\
\hline\hline
1&12&15&3&$\frac{1}{4}$
&
$\ver(H_{tri-up})=\{(12,10),(15,10),(18,13),(18,13)\}$
\\
2&15&15&3&$\frac{1}{4}$
&
\\
\hline
\end{tabular}
\end{center}
}
%%%%%%%%%%%%%
\subsubsection{Same trapezoid, two ways.}
A given {\ghysteron} can be parametrized in more than one way. 
Let $\ver(H)=\{(4,0),(8,0),(8,1),(10,1)\}$ and $H$ have linear sides, with the corresponding $\frac{s_r}{s_l}=2$. We can use $I=1$ or $I=4$ and the parametrizations 
%%%%%%%%
\ba
\label{eq:trap1}
\pP^2&=&
\begin{tabular}{|l|ll|ll|}
\hline
$k$&$\alpha_k$&$\beta_k$&$h_k$&$\mu_k$\\
\hline
1&4&8&2&$\frac{1}{4}$
\\
2&6&8&2&$\frac{1}{4}$
\\
\hline
\end{tabular},
\;\;
\pP^8=
\begin{tabular}{|l|ll|ll|}
\hline
$k$&$\alpha_k$&$\beta_k$&$h_k$&$\mu_k$\\
\hline
1&4&8&0.5&$\frac{1}{4}$
\\
2&4.5&8&0.5&$\frac{1}{4}$
\\
3&5&8.5&0.5&$\frac{1}{4}$
\\
4&5.5&8.5&0.5&$\frac{1}{4}$
\\
5&6&9&0.5&$\frac{1}{4}$
\\
6&6.5&9&0.5&$\frac{1}{4}$
\\
7&7&9.5&0.5&$\frac{1}{4}$
\\
8&7.5&9.5&0.5&$\frac{1}{4}$
\\
\hline
\end{tabular}
\ea
Here $\pP^8$ is obtained with a uniform partition of $[0,1]$ to $\{w_0,w_1,w_2,w_3,w_4\}=\{0,0.25,0.5,0.75,1\}$, and Algorithm in Sec.~\ref{sec:multilevel} yields $K=8$ {\uhysteron}s, each with $h^*=0.5$. Denote by $[r]$ the integer part of $r$. We have
$
{\hH^*}(\pP^8;\cdot) = 
\sum_{k=1}^8 
\hH(\tfrac{1}{4},{4+\tfrac{k-1}{2},8+\tfrac{[k-1/2]}{2}},h^*;\cdot)
$. 
In other words, both $\pP^2$ and $\pP^8$ yield the same bounding curves in $H$. Since it is more efficient to have $K=2$ than $K=8$, model $\pP^2$ is preferable over $\pP^8$. We note that even though $H(\pP^2)=H(\pP^8)$, the operator $\hH(\pP^2)\neq \hH(\pP^8)$ since the images $\hH(\pP^2;u)\neq \hH(\pP^8;u)$ when $u \not \in \sS(H)$.

%%%%%%%%%%%%%%%  kite example: skip
\myskip{
\subsubsection{Comparison of {\klinear} vs {\knonlinear} model}
\label{sec:klin}

%%%%%%%%%%
\begin{figure}[ht]
\begin{tabular}{cc}
\includegraphics
[height=30mm]{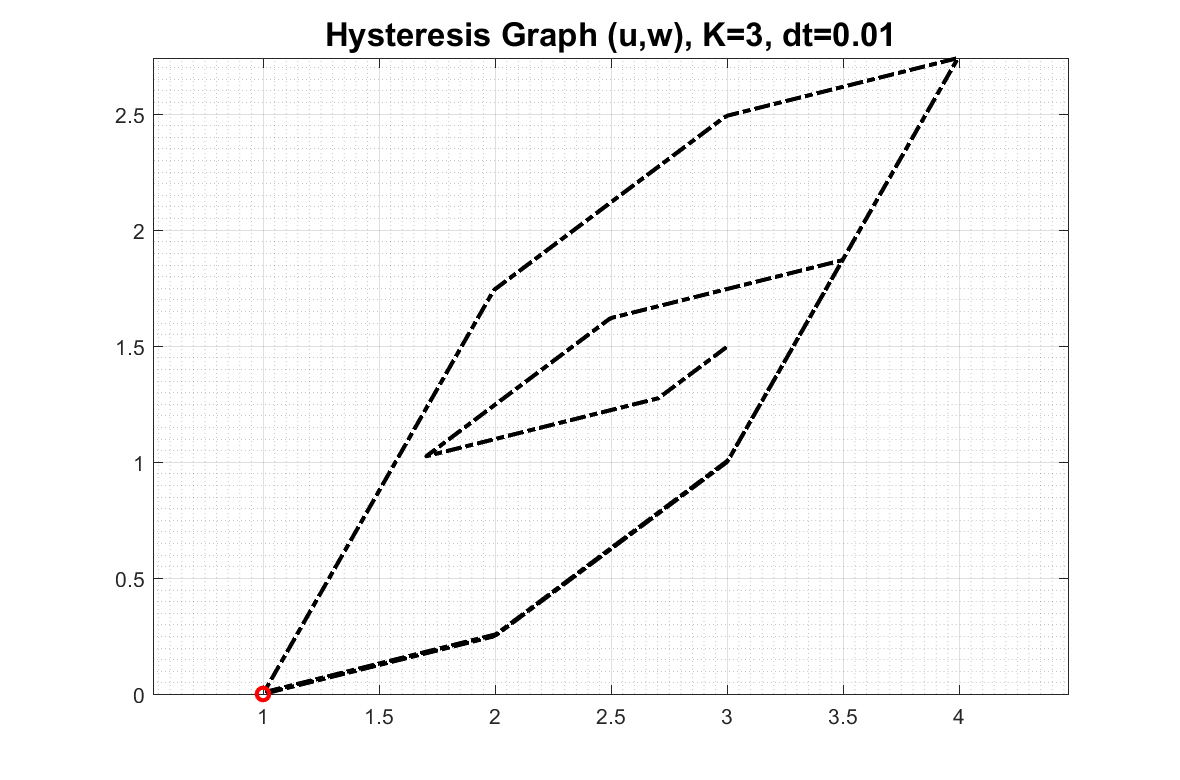}
&
\includegraphics
[height=30mm]{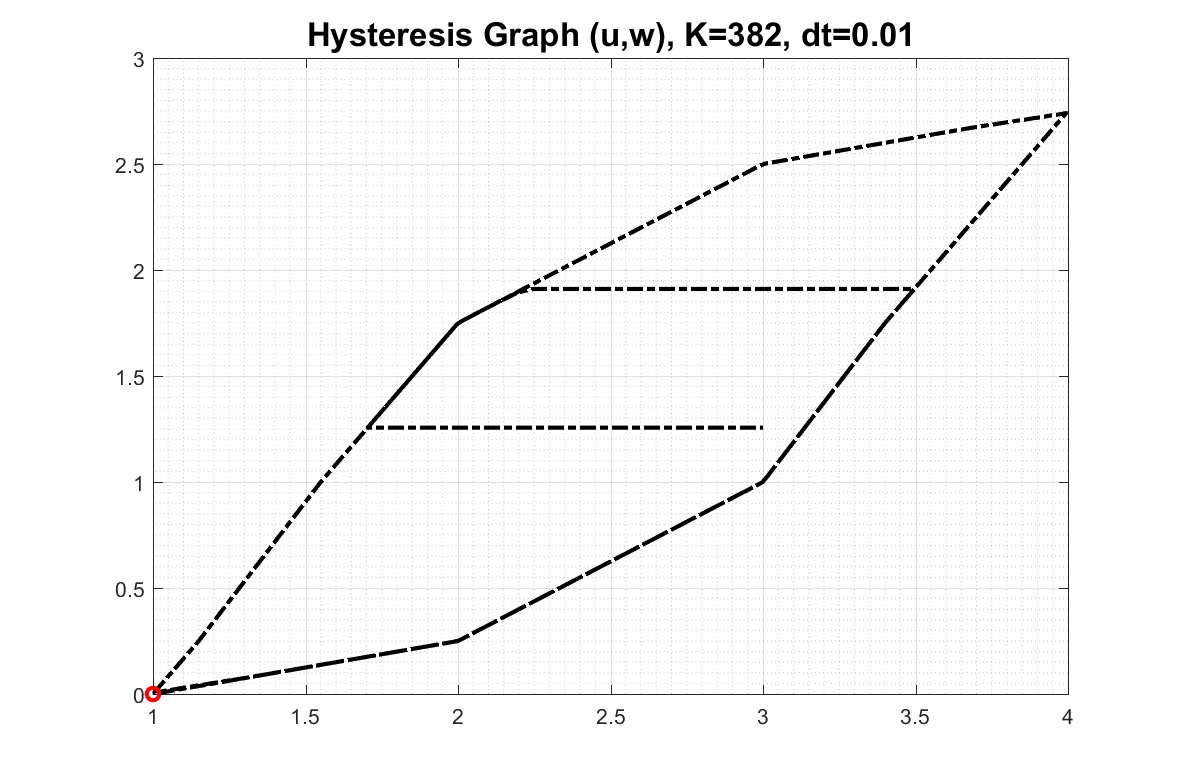}
\\
(a) $\pklin%_{stacksquashed} 
= [1/4,1,1,\infty;1/2,1,2,\infty;1,1,3,\infty]$
&
(b) $\pknon$ with $K=382$
\end{tabular}
\caption{Hysteresis graph with \klinear\ for ``kite'' shape (left) and  \knonlinear\ for the same shape; here $\uU=[1,4,1,2,3.5,1.7,3]$. 
\label{fig:klinear}
}
\end{figure}

We consider a simple example with $K=3$ which produces $H$ of convex ``kite'' shape 
\ba
\label{eq:three}
\hH(\pP^{3,\infty};u) = \hH(\tfrac{1}{4},1,\infty;u) + \hH(\tfrac{1}{2},0,1,\infty;u)+\hH(1,0,2,\infty;u).
\ea
See Fig.~\ref{fig:klinear}(a). The graph $H$ is convex-concave, as some portions of its sides are parallel to the sides on opposite ends. The largest slope of the sides equals $\tfrac{1}{4}+\tfrac{1}{2}+1=\tfrac{7}{4}$. 

Next we identify a \knonlinear\ model for the same $H$. With a process described in Sec.~\ref{sec:knon}, we identify the corners of $H$, and draw it as a union of $I=5$ trapezoids, each parametrized with \knonlinear\ model. The result is the graph shown in Fig.~\ref{fig:klinear}~(b). We see a big difference in the secondary curves between \klinear\ and \knonlinear. 
}

%%%%%
\subsubsection{Smooth graph with \klinear\ model, \kpreisach\ model, and \knonlinear\ model}
\label{sec:convex}
Now let $H$ be symmetric convex-concave, with $\fr(u)=(u-1)^2+1/3(u-1)$. The left curve is obtained on $[u_{min},u_{max}]$ by symmetric reflection with $\fl(u)=\fr(u_{max})-\fr(u_{min}+u_{max}-u)$. 

\begin{figure}[ht]
\begin{tabular}{cccc}
\includegraphics
[height=30mm]{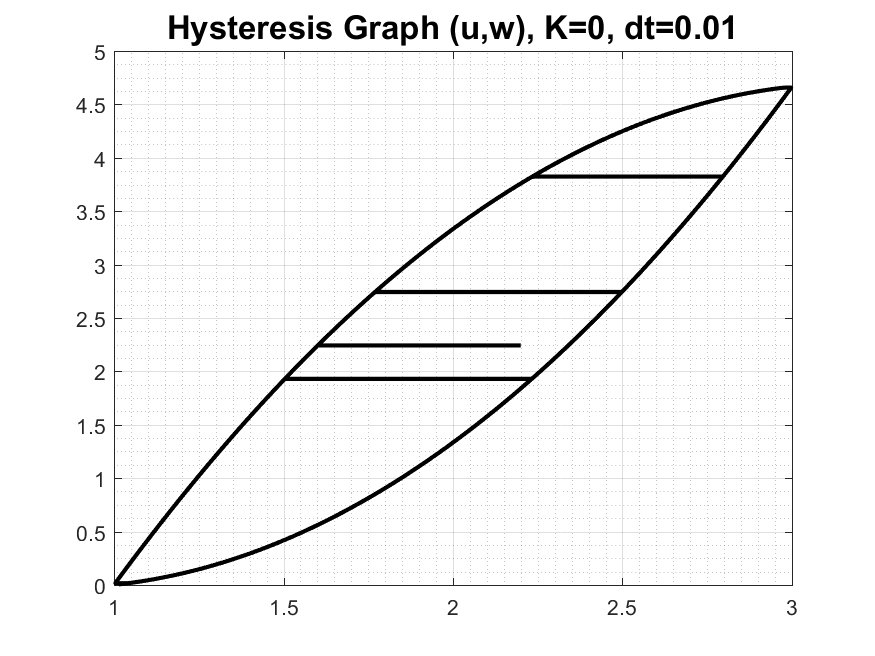}
&
\includegraphics
[height=30mm]{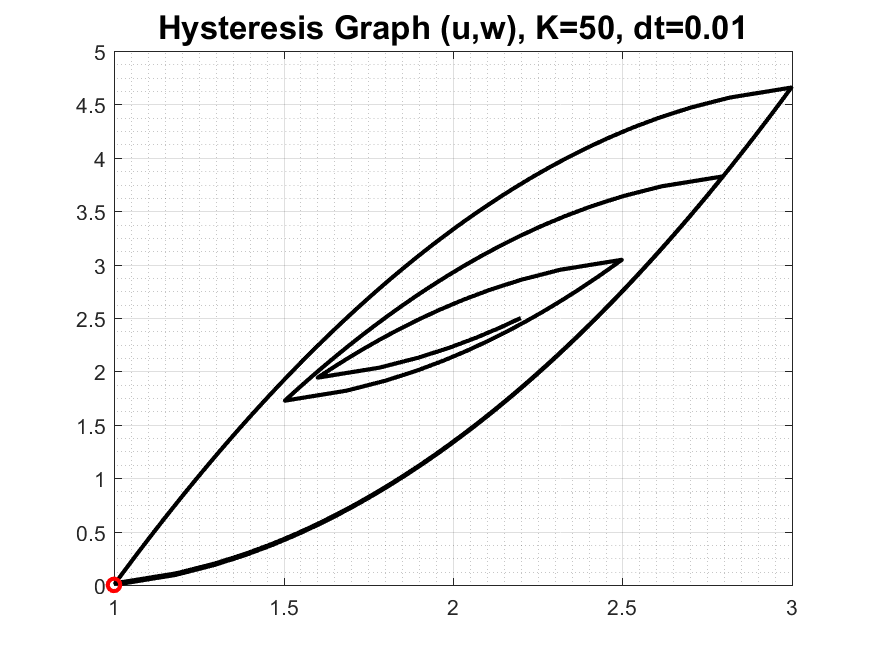}
&
\includegraphics
[height=30mm]{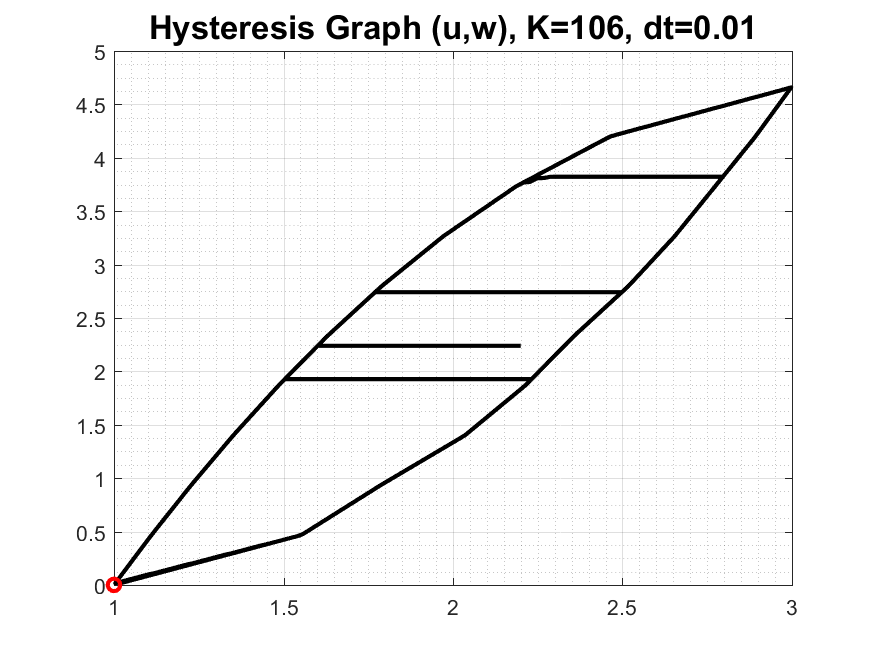}
&
\includegraphics
[height=30mm]{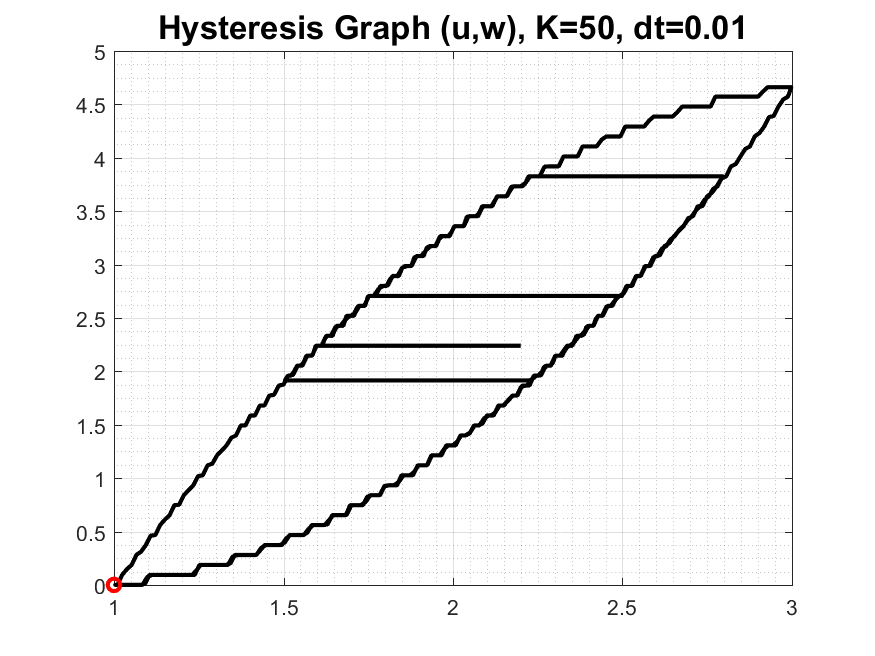}
\\
generalized&\klinear&\knonlinear&\kpreisach
\\
$\pgen$
&
$\pklin$, $K$=50
&
$\pknon$, $K$=106
&
$\ppreisache$, $K$=50, $\eps=1/10$
\end{tabular}
\caption{
Hysteresis graph with \klinear\ for a smooth curve $\fr(u)=(u-1)^2+1/3(u-1)$. Note the richness of interior loops with \klinear\ model with $\pklin$. Here
$\uU=[1,3,1,2.8,1.5,2.5,1.6,2.2]$. 
\label{fig:convex}
}
\end{figure}

We compare the output for generalized play, with $\pgen$, the \klinear\ model, \knonlinear, and \kpreisach. Fig.~\ref{fig:convex} illustrates $\hH(\pP;u)$ obtained for a particular $u\in \PL(\uU)$ with $\uU$ designed to sweep $H$ as well as to show a few secondary loops. 

%%%%%%%%%%

%%%%%%%
\subsubsection{Adsorption-desorption hysteresis graph from experimental data}
\label{sec:adsorption}
Now we calibrate hysteresis functional $\hH(\pP;\cdot)$ for realistic experimental adsorption--desorption data for methane CH4 from \cite{JessenK}. This data is fit in \cite{JessenK} to a particular algebraic model called  Langmuir isotherm $\gamma(u)=V\tfrac{Bu}{1+Bu}$ for each of the desorption $\fl(u)$ and adsorption $\fr(u)$ curves, respectively.  The data is given in Tab.~\ref{tab:adsorption}.

%%%%%%%%%%%%%%
\begin{table}
\centering
\begin{tabular}{l|ll}
&V&B\\
\hline
adsorption, $\fr(\cdot)$&
811&0.00237\\
desorption, $\fl(\cdot)$&
543&0.0382\\
\hline
\\
\end{tabular}
\caption{\label{tab:adsorption}Data for the Langmuir fit of the adsorption-desorption curves in Sec.~\ref{sec:adsorption}.}
\end{table}
%%%%%%%%%%%%

\begin{figure}[ht]
\begin{center}
\begin{tabular}{ccccc}
\includegraphics
[height=30mm]{cimages/ch4_generalized_secondary.png}
&
 \includegraphics
[height=30mm]{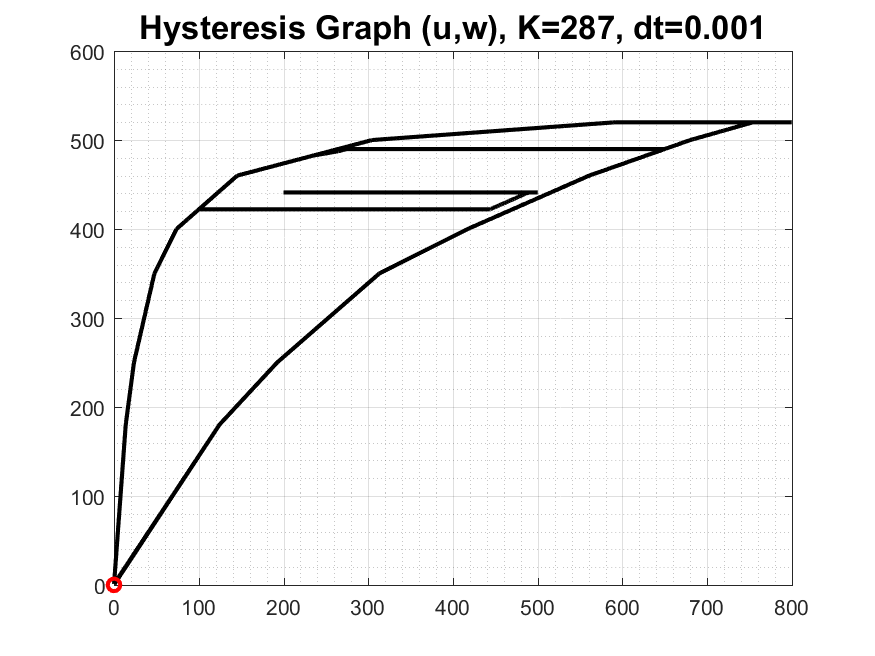}
&
\includegraphics
[height=30mm]{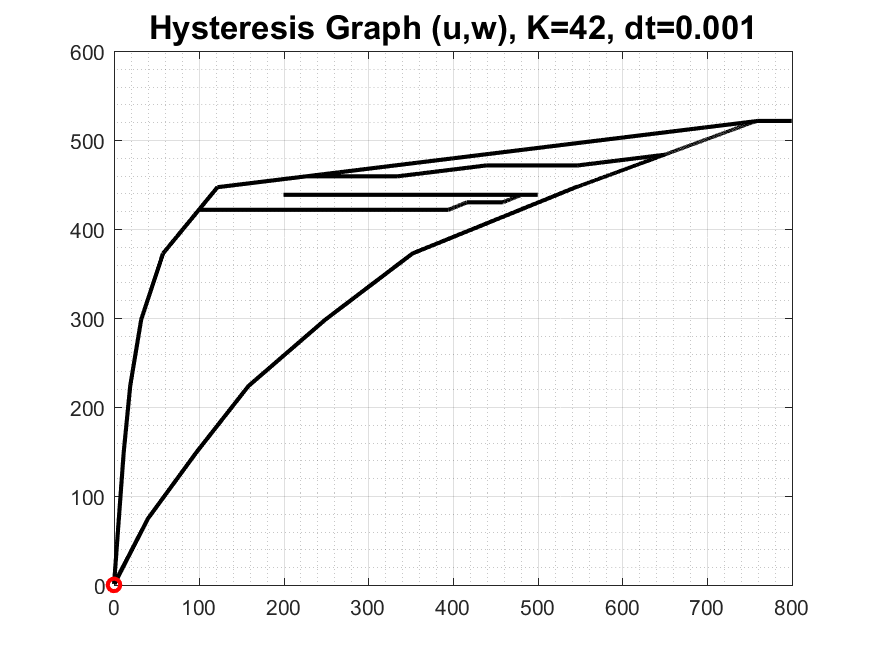}
&
\includegraphics
[height=30mm]{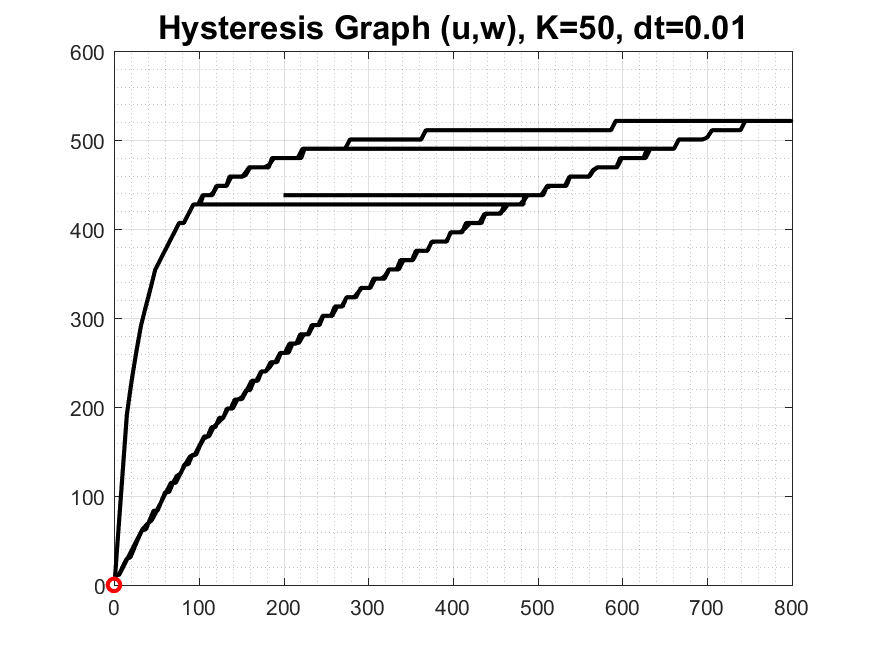}
\\ 
generalized&\knonlinear&\knonlinear&\kpreisach
\\
$\pgen$
&
$\pknon$, adaptive, $K=287$
&
$\pknon$, uniform, $K=42$
&
$\ppreisache$, $K=50$,$\eps=1/10$
\end{tabular}
\caption{Comparison of hysteresis graphs for the CH4 adsorption-desorption curves 
from \cite{JessenK}
$\hH$ plotted with $\uU= 
[0,800,0,650,100,350,500,200]$. The \klinear\ model is not available since $H$ is not symmetric.  The graphs can be compared to those in Fig.~\ref{fig:teaser}.
Note that the graph for $\pgen$ is annoted specially with $K=0$.
\label{fig:compare-adsorption}
}
\end{center}
\end{figure}

We use this data to produce $H$ and calibrate $\hH(\pP;)$ with generalized play, \knonlinear\ and \kpreisach\ models. For \knonlinear, we consider $I=7$ {\ghysteron}s with a uniform partition or an adaptive partition, both shown in Fig.~\ref{fig:calibrate}, with the corresponding $K=42$ and $K=287$.  The latter adheres more closely to $H$, and seems a better model in spite of a large $K$. 
See Fig.~\ref{fig:compare-adsorption} for illustration. 

%%%%%%%%%%%

%\end{document}
%%%%%%%%%%%%
\section{Analysis of evolution with \kgeneral\ hysteresis}
\label{sec:analysis}
To discuss the well-posedness of \eqref{eq:pde} with hysteresis, and to prove stability of the numerical scheme, we first formulate some auxiliary results for the stationary problem for the related system
\bsub 
 \label{eq:coupODEab}
\ba
 \label{eq:coupODEa}
\tfrac{d}{dt}(a(u) + b(v)) &=& f,
\\
\tfrac{d}{dt}b(v) + \cg(\fr(u),\fl(u); v) &\ni& 0. \label{eq:coupODEb}
\ea
\esub
The system \eqref{eq:coupODEab} corresponds to \eqref{eq:pde} with $A = 0$ and generalized play $\hH(u) = \hH(\pgen;u)$ obtained from the solution to the initial-value problem \eqref{eq:gplay}. Our results extend those we proved in \cite{PeszShow20} for \knonlinear. 

We start by writing \eqref{eq:coupODEab} in the form
\bsub 
 \label{eq:coupSYS}
\ba
 \label{eq:coupSYSa}
\tfrac{d}{dt}a(u) - \xi &=& f,\quad \xi \in \cg(\fr(u),\fl(u); v),
\\
\tfrac{d}{dt}b(v) + \xi  &\ni& 0. \label{eq:coupSYSb}
\ea
\esub
%
%%%%%%%%%%%%%%%%%%%%%%%%%%%%%%%%%%%%%%%%%%%%%%%%%%%%%%%%%%%%%%%%%
We prove properties of \eqref{eq:coupSYS} with {\kgeneral} in Sec.~\ref{sec:ode}, which we later use in Sec.~\ref{sec:odefd} for the numerical schemes for \eqref{eq:coupODEab}.

Note that for each $u \in \R$ the graph $\cg(\fr(u),\fl(u);\cdot)$ is maximal monotone, and for each $v \in \R$ the graph $- \cg(\fr(\cdot),\fl(\cdot);v)$ is maximal monotone. 
This is the form which occurs in the coupling with PDEs to give well-posed initial-value problems in the space $L^1 \times L^1$ for systems such as \eqref{eq:coupODEab}. 
We will add to \eqref{eq:coupODEa} an appropriate operator $A$ in $L^1$ for the PDE, and in Sec.~\ref{sec:pde} we use these towards the statement on well-posedness of \eqref{eq:pde}; later we use these to show the stability of an implicit-explicit numerical scheme for \eqref{eq:pde} when $A$ is an advection{-diffusion} operator {in $L^1$}.  
Both the well-posedness results for these systems as well as estimates for their discrete approximations depend on these results.  

%%%%%%%%%%%%%%%%%%%%%%%%%%%%%%%
\subsection{Estimates for \eqref{eq:coupSYS}}
\label{sec:ode}
Implicit-difference approximations of \eqref{eq:coupSYS} lead to consideration of the following systems.

\begin{lemma} 
\label{Cis2accretive}
If assumptions \eqref{assumptions} hold,
then  solutions $(u,v)$ and $(\bar u, \bar v)$ of 
\bsub \label{eq:Cis2accretive}
\ba
a(u) - \xi = f,\ b(v) + \xi = b(g), \ \xi \in \cg(\fr(u),\fl(u);v),
\\
a(\bar u) - \bar \xi = \bar f,\ b(\bar v) + \bar\xi = b(\bar g), \ \bar\xi \in \cg(\fr(\bar u),\fl(\bar u);\bar v),
\ea
\esub
satisfy the estimates
\begin{equation}  \label{eq:order-accretive}
(a(u) - a(\bar u))^+ + (b(v) - b(\bar v))^+ 
\le (f - \bar f)^+ + (b(g) - b(\bar g))^+,
\end{equation}
and
\begin{equation}  \label{eq:accretive}
|a(u) - a(\bar u)| + |b(v) - b(\bar v)| 
\le |f - \bar f| + |b(g) - b(\bar g)|.
\end{equation}
\end{lemma}
\begin{proof}
Subtract the respective equations,
$$
(a(u) - a(\bar u)) - (\xi - \bar \xi) = f - \bar f,\ (b(v) - b(\bar v))   
+ (\xi - \bar \xi) = b(g) - b(\bar g),
$$
multiply by $\sgn_0^+(u-\bar u)$ and then $\sgn_0^+(v-\bar v)$ and add to obtain
\begin{equation*}
(a(u) - a(\bar u))^+ + (b(v) - b(\bar v))^+ + (\xi - \bar \xi)
\big(\sgn_0^+(v - \bar v) - \sgn_0^+(u - \bar u)\big) \le 
(f - \bar f)^+ + (b(g) - b(\bar g))^+.
\end{equation*}
Below we verify the third term is non-negative, so we obtain
\eqref{eq:order-accretive}.
The corresponding estimates hold  for the negative parts 
$( \cdot )^-$ and then for their sum,
\eqref{eq:accretive}.
This shows the solution is order-preserving.

Finally, we check that $(\xi - \bar \xi)
\big(\sgn_0^+(v - \bar v) - \sgn_0^+(u - \bar u)\big) \ge 0$. If $u = \bar u$ or $v = \bar v$ it is $\ge 0$ by the monotonicity of $\cg$ in each variable. Otherwise consider the displacement $(u - \bar u,v - \bar v)$ in each quadrant of $\R \times \R$ with $(u,v)$ on the constraint and $(\bar u,\bar v)$ in the interior. In the first quadrant this term is $(\xi - \bar \xi)(1 - 1)$ and in the third it is $(\xi - \bar \xi)(0 - 0)$. In the second it is $(\xi - 0)(1 - 0) \ge 0$ and in the fourth quadrant it is $(\xi - 0)(0 - 1) \ge 0$.
\end{proof}
%%%%%%%%%%%

Corresponding results hold as well for K-generalized play as given by \eqref{eq:pk}. 

\begin{proposition}     \label{CisKaccretive}
Assume in \eqref{eq:pk} that for each $1 \le k \le K$, the functions $a(\cdot), b_k(\cdot), \flk,\frk$ satisfy \eqref{assumptions}.
Then the solutions $(u,v),\ v = (v_k)$ and $(\bar u, \bar v),\ \bar v = (\bar v_k)$ of the respective systems
%%%%%%%%%%
\bsub       \label{eq:KcoupSYS}
\ba     \label{eq:KcoupSYSa}
a(u) - \sum_{k=1}^K \xi_k = f, \ b_k(v_k) + \xi_k = b_k(g_k), \ \xi_k \in \cg(\gamma_{r k}(u),\gamma_{\ell k}(u);v_k(t)), \ 1 \le k \le K,
\\          \label{eq:KcoupSYSb}
a(\bar u) - \sum_{k=1}^K \bar \xi_k = \bar f,   
\ b_k(\bar v_k) + \bar \xi_k = b_k(\bar g_k), 
\ 
\bar \xi_k \in \cg(\gamma_{r k}(\bar u),\gamma_{\ell k}(\bar u);\bar v_k(t)), \ 1 \le k \le K,
\ea
\esub
with data $(f,g),\ g = (g_k) \in \R^K$ and $(\bar f,\bar g),\ \bar g = (\bar g_k)$
satisfy the estimates
\begin{equation}  \label{eq:order-Kaccretive}
(a(u) - a(\bar u))^+ + \sum_{k=1}^K(b_k(v_k) - b_k(\bar v_k))^+ 
\le (f - \bar f)^+ + \sum_{k=1}^K (b_k(g_k) - b_k(\bar g_k))^+,
\end{equation}
and
\begin{equation}  \label{eq:Kaccretive}
|a(u) - a(\bar u)| + \sum_{k=1}^K |b_k(v_k) - b_k(\bar v_k)| 
\le |f - \bar f| + \sum_{k=1}^K|b_k(g_k) - b_k(\bar g_k)|.
\end{equation}
\end{proposition}
This follows from the same proof as in Lemma~\ref{Cis2accretive} with the corresponding estimates $(\xi_k - \bar \xi_k) (\sgn_0^+(v_k - \bar v_k) - \sgn_0^+(u - \bar u)) \ge 0$ for $1 \le k \le K$. 

The results hold also when $\flk, \frk$, are permitted to be maximal monotone relations without common points of multiple-values.  See \ourcite{Visintin94}{VIII.2}.

%%%%%%%%%%%%%%%%%%%%% move to end of section 2.4. 
%%%%%%%%%%%%%%
\subsection{PDE coupled to Hysteresis}
\label{sec:pde}
Now we consider the evolution equation \eqref{eq:pde} in which $\hH(\cdot)$ is a generalized play given by \eqref{eq:pk} and $A$ is a PDE on a domain $\Omega$ for advection or diffusion with appropriate boundary conditions. For the case of generalized play, the equation takes the form of the system \eqref{eq:coupODEab} with $A$ added to the first equation. The estimates in  Lemma~\ref{Cis2accretive} show that the operator $\cg$ on $\R \times \R$ given by $(a(u), b(v)) \to (- \xi, \xi)$ with $\xi \in \cg(\fr(u),\fl(u);v)$ is {\em accretive}. This is the analogue of monotone on a Banach space $B$, and with the range condition  $\Rg(I + \cg) = B$ it is called {\em m-accretive}.  
If $\cg$ is m-accretive, $\vz \in \overline{\Dom(\cg)}$, and $f \in L^1(0,T;B)$, then the initial-value problem \eqref{eq:ode}
has a unique {\em integral solution}. This is a limit in $C([0,T];B)$ of implicit-difference approximations \eqref{eq:fd}. For Banach spaces which possess the Radon-Nikodym property, in particular, for finite-dimensional spaces, if $f$ has bounded total variation, the integral solution is a strong solution  $v \in W^{1,\infty}(0,T;B)$; see \ourcite{Visintin94}{XII.4}.

Corresponding estimates for the system in the product space $B = L^1(\Omega) \times L^1(\Omega)$ show the initial-value problem for \eqref{eq:coupSYS} is well posed, that is, it has a unique {integral solution} in $C([0,T];B)$. The same is true with $A(u)$ added to the first equation if $A$ is accretive on $L^1(\Omega)$. Moreover the same argument with Proposition~\ref{CisKaccretive} extends to the other forms of play hysteresis and we obtain the following general result.
See \ourcite{Visintin94}{VIII.3} and \cite{Visintin93, LittShow94, PeszShow20}.

\begin{proposition}
\label{prop:pde}
Assume the conditions \eqref{assumptions} and that $A$ is m-acretive on the Banach space $L^1(\Omega)$.
Then the initial-value problem for \eqref{eq:pde} with $\hH$ as in \eqref{eq:pk}, under \eqref{eq:pi}, is well posed. In particular, it holds when $\hH(\cdot)$ is either \general\ with $\pgen$ or \knonlinear\ with $\pknon$ or regularized {\kpreisach} (with $\ppreisache$ or $\ppreisachs$).
\end{proposition}

%%%%%%%%%%%%%%%%%%%%%%%%%%%%%%%%%%%%%%%%%%%%%%%%%%%%%%%%%%

%%%%%%%%%%%%%%%%%%%%%%%
\section{Numerical scheme for an ODE  and PDE with hysteresis} 
%%%%%%%%%%%%%%%%%%%%
%%%%%%%%%%%%%%%%%%%%%%%%%%%%%%%%%%%%%%%%%%%%%%%%%%%
%%%%%%%%%%
\label{sec:numerical}
Now we discuss some practical challenges when solving numerically a dynamical problem involving $\hH(\cdot)$. 
Numerical models with hysteresis were discussed and analyzed in \cite{VerdiVis85,VerdiVis89,HoffMeyer89}, all for Preisach type models
and parabolic PDEs. Here we use different techniques motivated by explicit schemes for transport. 

We define first an implicit numerical scheme for the ODE \eqref{eq:coupODEab} 
and discuss its solvability, solver, and convergence rate.
Next we consider \eqref{eq:pde} with an advective transport $A(u)$. We apply an Ex-Im scheme (explicit-implicit) in which the PDE term $A(u)$ is treated first explicitly in time, followed by an implicit solver for  the hysteresis accumulation term $u+w$ with the scheme for the ODE \eqref{eq:coupODEab}.  

%%%%
\subsection{Numerical scheme for ODE \eqref{eq:coupODEab}}
 \label{sec:odefd}
 We approximate $(u,w)$ solving \eqref{eq:coupODEab}  with $w \in \hH(\pP^K;u)$. We seek $(U^n,W^n)$ defined by 
\bsub
\label{eq:coupFD}
\ba
\label{eq:coupFDu}
a(U^n)-a(U^{n-1})+W^n-W^{n-1}&=&\tau F^n,
\ea
Here we set $F^n=f(t_n)$ for smooth $f$, or $F^n=\tfrac{1}{\tau}\int_{t_{n-1}}^{t_n} f(s)ds$ otherwise. 
Next we need $W^n$ as a function of $U^n$; we also seek $V^n=(V^n_k)_k$ which approximate $v_k$ solving \eqref{eq:odelp}, and $w=b(v)$ solving \eqref{eq:odeb}. 
For these we have a unified formula \eqref{eq:resk}, and we have
\ba
\label{eq:coupFDw}
W^n= \ggres(\pP;V^{n-1};U^n)=\sum_{k=1}^K \mu_k b_k(\resk(V^{n-1}_k;U^n)).
\ea
\esub
The function $\ggres(\pP;\vbar;U)$ has properties which extend from these of each $\resk$ and of $b_k$. In particular, these are the monotonicity properties expressed in Lemmas~\ref{lem:gres} and \ref{lem:b}. From these the next result follows.

\begin{lemma}
\label{lem:solver}
The solution $(U^n,W^n)$ to \eqref{eq:coupFD} exists and is unique. 
\end{lemma}

\begin{proof}
We substitute $W^n$ from \eqref{eq:coupFDw} into \eqref{eq:coupFDu} to obtain the stationary problem for $U^n$ 
\ba
\label{eq:statgenehys}
a(U^n) + \ggres(\pP;V^{n-1};U^n) = a(U^{n-1})+W^{n-1}+ \tau F^n.
\ea
Since $a(\cdot)+ \ggres(\pP;\vbar;\cdot)$ is monotone strictly increasing, with  {range $\R$}, this equation has a unique solution $U^n$. Once $U^n$ is known we calculate $(V^n_k)_k$ and $W^n$ from \eqref{eq:coupFDw}.
\end{proof}
We also prove additional properties. They are analogues of what was proven in \cite{PeszShow20} for \knonlinear. We only need to prove the counterparts of \ourcite{PeszShow20}{Lemma 4.4} and \ourcite{PeszShow20}{Prop.4.5}. We also need  an additional result similar to that in \ourcite{PeszShow20}{Prop.4.6}.
\begin{lemma}
\label{lem:order}
The solutions to \eqref{eq:coupFD} satisfy
\ba
\label{eq:specialg}
\abs{a(U^n)-a(U^{n-1}}+\sum_k \mu_k\abs{b_k(V_k^n)-b_k(V_k^{n-1})}=\tau \abs{F^n}.
\ea
\end{lemma}
\begin{proof}
The lemma was proved for \knonlinear\ model in \cite{PeszShow20}. It remains to verify \eqref{eq:specialg} for $\pP=\pgen$ and $\hH(\pP;\cdot)$ represents the \kgeneral\ model. 
We recall Lemma~\ref{lem:gres}. From this monotonicity result we can write $W_k^n-W_k^{n-1}=\mu_k b_k(V_k^n)-\mu_k b_k(V_k^{n-1})=\psi_k(U^n-U^{n-1})$ with some 
$\psi_k\geq 0$; this is similar to the use of mean value theorem, but we need not identify $\psi_k$ with a derivative. 
Similarly, we also have $a(u)-a(\ubar)=\chi(u-\ubar)$, with positive $\chi\geq0$. Thus \eqref{eq:coupFDu} is equivalent to
\bas
(\chi+\sum_k \psi_k)(U^n-U^{n-1})=\tau F^n,
\eas
Next, $F\geq 0$ implies $U^n-U^{n-1}\geq 0$, and $F\leq 0$ implies $U^n-U^{n-1}\leq 0$, thus we obtain the desired result. \end{proof}

%%%%%%
\subsection{Implementation, solver, and convergence of discrete scheme \eqref{eq:coupFD}}
\label{sec:solver}

%%%%%%%%%%%
\begin{figure}[ht]
\centering
\begin{tabular}{cc}
\includegraphics
[height=40mm]{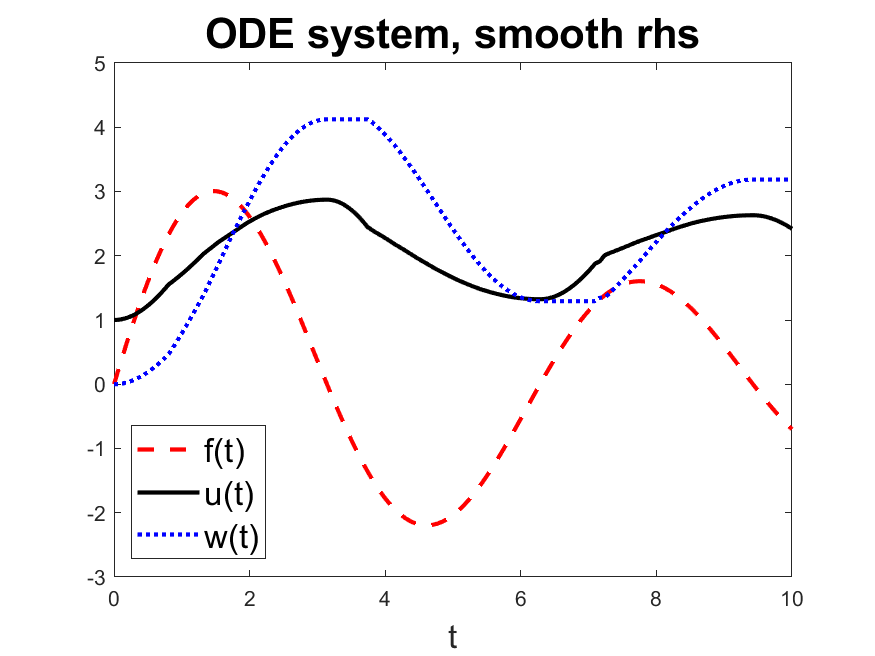}
&
\includegraphics
[height=40mm]{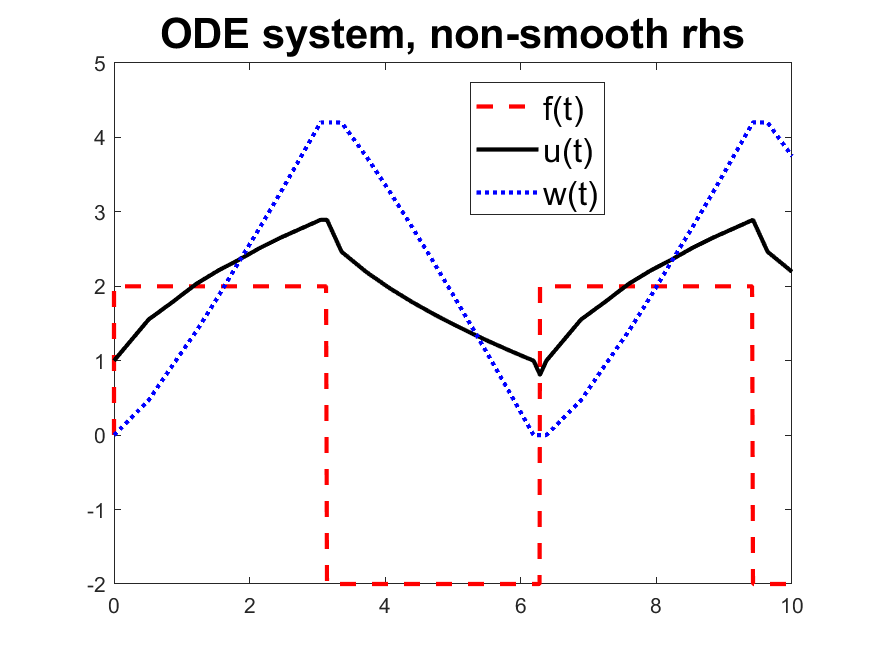}
\\
(a)&(b)
\\
%%%%
\includegraphics
[height=40mm]{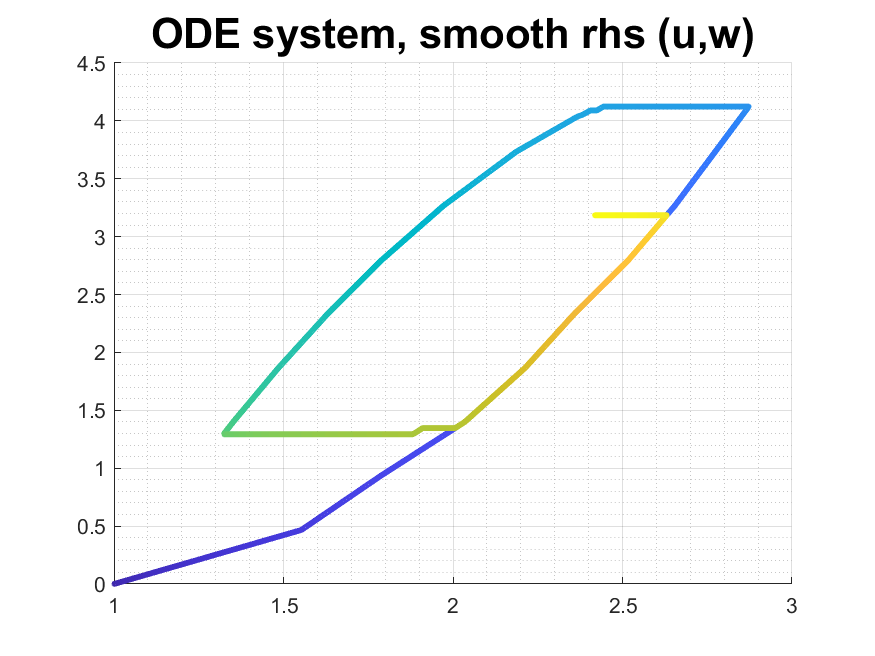}
&
\includegraphics
[height=40mm]{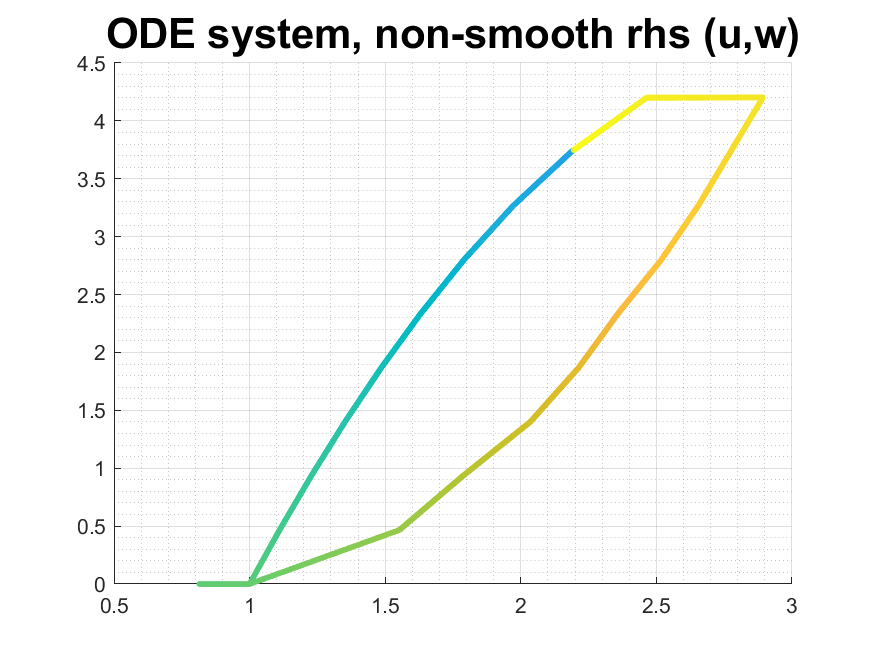}
\\
\\
(c)&(d)
\end{tabular}
\caption{Illustration of ODE with hysteresis with  $\pknon$ from Sec.~\ref{sec:convex}. Top: $f(t),U_{\tau},W_{\tau}$ for (a) $f=f^{cont}$ and (b) $f=f^{disc}$, and $\tau=10^{-4}$. Bottom: (c) the trace $(U^n,W^n)_{n}$ of $\hH(\pknon;)$ for 
(c) $f=f^{cont}$ and (d) $f=f^{disc}$. 
\label{fig:ode}}
\end{figure}

We discuss how the solution to \eqref{eq:coupFD} is found in practice. We set-up an example for \eqref{eq:coupODEab} with $\hH(\pP;\cdot)$  for the convex--concave graph $H$ from Sec.~\ref{sec:convex} and $\pP$ as one of $\pgen,\pklin,\pknon,\ppreisache$. We
use $\uz=1$, a compatible $\vz$, and the source function as one of the two
\ba
T=10; \;\; f^{cont}(t)= 3.5 \mathrm{sin}(t) e^{-0.1 t};\; f^{disc}(t)=\mathrm{sign}(f^{cont}(t)).
\ea
The experiment is designed so that $(U_{\tau},W_{\tau})$ for both $f^{cont}$ and $f^{disc}$ loop almost all over $H$; see Fig.~\ref{fig:ode}. Here $U_{\tau},W_{\tau}$ are the piecewise linear interpolants of $(U_n)_n$ and $(W_n)_n$, respectively. 
 
%%%%%%%%%%
\subsubsection{Solver for \eqref{eq:coupFDu}}
Since \eqref{eq:statgenehys} is a scalar root-solving problem on $\R$, there are many excellent solvers. Newton's method uses derivative information and converges fast close to the root \cite{Kelley} whenever the residual is a smooth function of $u$. However, Newton's method may have occasional difficulty handling only piecewise differentiable functions such as $\ggres(u)$. In this (semismooth) case \cite{Ulbrich} the solver still may converge, but with occasional failures.  A robust alternative is a method we call \Root\ which brackets the root and then uses a secant method \cite{Brent}; this method is implemented in {\tt fzero} in MATLAB. The \Root\ solver is robust but also much slower than Newton iteration. Both require an implementation of the calculation for $\ggres(U^n)$; Newton solver also requires its derivative.  

We test solver performance. Of interest is the average $\Nave$ number of iterations as well as the computational time shown in Table~\ref{tab:solver}. For perspective we show the case without hysteresis, which still requires a nonlinear solver at every time step $n$. 

%%%%%%
\begin{table}
%%%%%%%%% 
%%% NO HYSTERESIS  
% 0.1,1,0,generalized,'Root');
\begin{tabular}{|l|lll|lll|}
\hline
%(Root), 
Solver
&
\multicolumn{3}{c|}{Root}
&
\multicolumn{3}{c}{Newton}
\\
\hline
case/$\tau$
&
0.1&0.01&0.001&
0.1&0.01&0.001
\\
\hline
%%%%%%%%%%%%%%%%%%%%%%%%%%%%%%%%%%%%%%%%%%
$\fl$=$\fr$
%&0.1 
&9.31 
% max_iters=13 maxres=1
(0.171)
%
%&0.01 
&7.33 
% max_iters=9 maxres=1
(1.01)&
% 
%&0.001
7.1782 
% max_iters=10 maxres=1
(9.44)&
%% 0.001,1,0,generalized,'Newton');
%(Newton), 
%&0.1 
3.68 
%max_iters=4 maxres=1.94786e-07
(0.117)&
%
%&0.01 
3 
%max_iters=3 maxres=4.25181e-09
(0.583)&
%&0.001 
3 
% max_iters=3 maxres=4.36763e-13
(5.29)
\\
\hline
%%%%%%%%%  smooth rhs 
%tic;HysteresisSolver(10,0.1,1,0,generalized,'Root');toc,
%Root
$\pgen$
&&&&&&
\\
$f^{cont}$
%&0.1 &
&9.08 
%max_iters=15 maxres=1
(0.218)&
%%
%&0.01 &
8.15
%max_iters=14 maxres=1
(1.057)&
%&0.001 &
8.602 
%max_iters=16 maxres=1
(10.49)
%%
%%(Newton), 
%&0.1 &
&3.29 
%max_iters=5 maxres=1.73712e-07
(0.113)
%(Newton), 
 %&0.01 &
&2.79 
 %max_iters=4 maxres=1.72131e-09
(0.529)&
%
%(Newton), &0.001 &
 2.78
 %max_iters=4 maxres=2.33862e-11
(4.776)
\\
%%%%%%%%%  non smooth; generalized
$f^{disc}$
&
%0.1 &
10.27 
 %max_iters=20 maxres=1
(0.174)&
%0.01 &
8.97 
 %max_iters=12 maxres=1
(1.135)&
%0.001 
9.50 
 %max_iters=13 maxres=1
(11.21)
%%%%%
%(Newton), 
%0.1
&3.73 
%max_iters=5 maxres=1.87262e-07
(0.124)& 
%(Newton), &0.01 &
2.98 
%max_iters=4 maxres=1.87262e-08
(0.590)&
%(Newton), &0.001 &
2.97 
 %max_iters=4 maxres=2.83507e-12
(5.059)
\\
\hline
%%%%%%%%%%%  nonsmooth  KNONLINEAR, #106 components
%% 0.1,1,0,knon_symmetric,'Root');
$\pknon$
&&&&&&
\\
%(Root), &0.1 &
$K$=$106$&
7.98 
%max_iters=14 maxres=1
(1.066)&
%(Root), &0.01 &
7.14 
%max_iters=16 maxres=1
(8.954)&
%(Root), &0.001 &
7.94 
%max_iters=17 maxres=1
(98.61)&
%(Newton), &0.1 &
2.42 
%max_iters=5 maxres=2.498e-15
(0.56)&
%(Newton), &0.01 &
2.07 
%max_iters=3 maxres=4.90233e-15
(4.58)&
%(Newton), &0.001 &
2.01 
%max_iters=3 maxres=2.105e-13
(47.40)
\\
%(Root), &0.1 & build with 30 
$K$=$921$&
9.01 
%max_iters=14 maxres=1
(8.93)&
%(Root), &0.01 &
7.7 
%max_iters=16 maxres=1
(77.14)&
%(Root), &0.001 &
8.46 
%max_iters=17 maxres=1
(861.2)&
%(Newton), &0.1 &
 2.98
%max_iters=5 maxres=2.498e-15
(5.12)&
%(Newton), &0.01 &
 2.09
%max_iters=3 maxres=4.90233e-15
(36.31)&
%(Newton), &0.001 &
2.01 
%max_iters=3 maxres=2.105e-13
(380.8)
\\
\hline
%%%% now 0.1,1,0,klin_symmetric,'Root');
$\pklin$
&&&&&&
\\
$K$=$50$
&
%(Root), &0.1 &
8.88 (0.455)
&
7.11 (3.54)
&
7.731 (39.65)
%%%(Newton), &0.1 &
&
3.07 (0.34)
&
1.15 (1.87)&
2.02(19.07)
\\
%%% change number of components 
%0.1,1,0,klin_symmetric,'Root');
%(Root), &0.1 &
K=200&
9.67 (1.63)
& 
7.58 (15.26)
&
8.08 (167.4)
&
3.08 (1.12)
&
2.47 (8.90)
&
2.04 (78.19)
\\
\hline
%%%%%%%%%%%%%%%%%%% Preisach
% 0.01,1,0,symmetric100,'Root');
%% HysteresisSolver(10,0.1,1,0,symmetric100,'Root');
%&0.1 &
$\ppreisache$
&&&&&&
\\
$K$=$100$
&
14.35 
 %max_iters=26 maxres=1
(1.56)&
%&0.01 &
11.92 
%max_iters=21 maxres=1
(13.07)&
%,&0.001 &
11.48
%max_iters=22 maxres=1
(132.92)
%%% WITHOUT REGULARIZAITON
% 0.01,1,0,symmetric100,'Newton');
%res=0.0466667 iter=11 step=3 t=0.03
%Error using HysteresisSolver (line 105)
%no convergence
&
--
&
--
&
--
\\
$\bh_*$& 
%%% with regulARIZATION
%0.1,1,0,symmetric100,'Newton');%0508763
% (Newton), &0.1 &
10.08 (1.23)&
8.63 (10.19)&
9.39 (113.04)& 
3.78 (0.84)&
3.014 (6.54)&
2.98 (67.09)
\\
\hline

\end{tabular}
%%%
\caption{Solver performance for 
scheme \eqref{eq:coupFD} and
experiments in Sec.~\ref{sec:solver}. The numbers tabulated are the number of iterations and the computational time (in parentheses). The first row presents the case without hysteresis. The second two rows compare the use of $\pgen$ for $f^{cont}$ and $f^{disc}$. Examples with $\pknon,\pklin,\ppreisache$ are reported for the more challenging $f^{disc}$. Newton solver does not converge for some time steps when $\pP=\ppreisache$, unless regularization with $\bh_*$ is used. 
\label{tab:solver}
}
\end{table}

The results in the Table show that the $\ppreisache$ model is the most sensitive to the solver choice, while $\pgen$ and $\pklin$ require the least amount of computational effort. As usual, when Newton converges, it converges faster than the Root solver. We use similar convergence criteria for both, with a combination of absolute ($10^{-14}$) and relative ($10^{-6}$) tolerance. As usual, smaller $\tau$ decreases somewhat the number of iterations, but overall the computational effort scales roughly linearly with the number of time steps. Newton solver performs poorly for the \kpreisach\ model, even though it can be made to work upon regularization.  

 The computational time for \kgeneral\ model is not significantly higher than that for case without hysteresis, and so is the number of iterations. As expected, higher $K$ requires more computational time than lower $K$, and the effort scales about linearly with $K$; compare, e.g., the $\pklin$ case with $K=50$ to that with $K=200$. 

%%%%%%%%%%%%%%%%%%%%%%%%%%%%%%%%%
\begin{table}
\centering
%%%%%%%%%%%%%%%%%%%
\begin{tabular}{|l|ll|ll|l|l|}
\hline
\hline
&\multicolumn{4}{c|}{$\pgen$ }
&\multicolumn{1}{c|}{$\pknon$}
&\multicolumn{1}{c|}{$\pknon,\bh_*$}
\\
\hline
&\multicolumn{2}{c|}{$f^{cont}$}
&\multicolumn{2}{c|}{$f^{disc}$}
&\multicolumn{1}{c|}{$f^{cont}$}
&\multicolumn{1}{c|}{$f^{disc}$}
\\
\hline
$\tau$&
$E_u$&$E_w$&
$E_u$&$E_w$
&$E_u$
&$E_u$
\\
\hline
0.1&0.0669978&0.108278
&0.111336 &0.134983
&
0.0592
&
0.10087
\\
0.01&0.0068853&0.011001
&0.010153 &0.008205
&
0.00710
%knon reg 0.0100    
&
0.014807    
%0.0079
\\
0.001&0.0006917&0.001092
&0.002179 &0.003043
&
0.00065
&0.001316
\\
\hline
\end{tabular}
\caption{
\label{tab:error}Error $E_u=\norm{u-U_{\tau}}{\infty}$ and $E_w=\norm{w-W_{\tau}}{\infty}$ for the experiments in Sec.~\ref{sec:solver} with $\pP$ from Sec.~\ref{sec:convex}. We use the fine grid solution $(U_{\tau_{fine}},W_{\tau_{fine}})$ with $\tau_{\infty}=0.0001$ as a proxy for the true solution, and run the experiments for several $\tau$; some are shown here. We seek the order $p$: $E_u = O(\tau^p)$.    For the smooth source $f^{cont}$, we get $p\approx 1$ for $\pknon$ or $\pklin$ ($p=0.98$), but less dependably so when $f^{disc}$ is used. In fact, for $\pknon$ it is necessary to use additional smoothing $\bh_*$ or of $f^{disc}$ to achieve $O(\tau)$ convergence, at least within reasonable range of $\tau$. }
\end{table}

%%%%%%%%%%%%%%%%
\subsubsection{Convergence rate for 
$\norm{u-U_{\tau}}{\infty} \to 0$ found in
\eqref{eq:coupFD}}
\label{sec:err}

Table~\ref{tab:error} is devoted to the rate.  Generally, if $u \in W^{1,1}$ then $v\in W^{1,1}$ as well, and so is $w=b(v)$.  However, it is not clear what to expect for $(u,w)$ solving \eqref{eq:coupODEab}. For nonsmooth $f$ such as $f=f^{disc}$, we expect $u+w\in W^{1,1}$, but the regularity of each $u(t)$ and $w(t)$ remains unclear. Since the analytical solution is not known, we use the fine grid solution. 

While generally we expect about linear rate of convergence $O(\tau)$ in $u$,  this expectation is not easy to confirm for graphs with large $K$, or even for $\pgen$, within the range of practical time steps.  As concerns $W$, the error seems stable, but we do not expect or always observe convergence, except for very simple $\pP$. 

%%%%%%%%%%%%
\subsection{Numerical scheme and analysis for a PDE}
\label{sec:pdenum}

We now consider the homogeneous IVP which specializes \eqref{eq:pde} to when $A=\partial_x (\alpha(\cdot))$. 
\begin{subequations}    \label{eq:ivp}
\ba
\Dt \left(a(u) + w\right)
+ \Dx \alpha(u) &=& 0, \;\;
u(x,0) = u_{init}(x)\,; v(x,0)=\vz(x)\;\;
\\
%\label{eq:uw}
w(x,t) &=& \hH(u(x,t))
\ea
\end{subequations}
The scheme we use is Ex-Im: explicit in the transport, with upwind treatment of advection, and implicit in the resolution of the nonlinearity under $\Dt$. For a spatial grid parameter $h$, we define the gridpoints $x_j=jh$ on the support of $\uz$, we set $U^0_j=\uz(x_j)$, $V^0_{j,k}=\vz_k(x_j)$, and at every time step $n>0$ we solve for $U_j^n \approx u(x_j,t^n)$ and 
$W_j^n \approx w(x_j,t^n)$ as follows:
\bsub
\label{eq:PDEFD}
\ba
\tfrac{1}{\tau} \left(a(U_j^n)-a(U_j^{n-1})+W_j^n-W_j^{n-1}\right)
+ \tfrac{1}{h}(\alpha(U_j^{n-1})-\alpha(U_{j-1}^{n-1}))=0; 
\\
W_j^n=
\sum_k \mu_k b_k(\resk(V_{k,j}^{n-1};U_j^n)).
\ea
\esub
Rearranging, we see that at every $n$ and $j$, one has to solve  for $(U_j^n,W_j^n)$ a problem analogous to \eqref{eq:coupFDu} with $F^n_j=-\tfrac{\tau}{h}(\alpha(U_j^{n-1})-\alpha(U_{j-1}^{n-1}))$.

We have the following result on weak $TV_T$ stability of \eqref{eq:PDEFD}. We state it  with a brief proof. 

\begin{proposition}
\label{prop:TVT}
Assume the CFL condition  $0\leq \tfrac{\tau}{h}\max_u \tfrac{d\alpha}{du}\leq 1$. The scheme \eqref{eq:PDEFD} is uniquely solvable, and is weakly $TV_T$ stable in the $L^1 \times (L^1)^K$ space. Namely, 
$TV_T(a(U),(W_k)_k) \leq C(T)$ where the constant $C(T)$ depends on time $T$ but not on $\tau,h$. 
\end{proposition}
\begin{proof}
The proof for $\pP=\pknon$ is given in \cite{PeszShow20} and relies on the properties of solvability (\cite{PeszShow20}, Prop.4.5), comparison principle (\cite{PeszShow20}, Lemma 4.4), and ordering (\cite{PeszShow20}, Prop.4.6). The proof for general play requires that we verify properties listed above. These are, respectively, Lemma~\ref{lem:solver}, Lemma~\ref{Cis2accretive}, and Lemma~\ref{lem:order} proven in this paper. 
\end{proof}

The stability result combined with the usual truncation error analysis \cite{leveque-red,Toro-book} which is fairly simple for this upwind scheme suggest that the convergence rate $\norm{u-U_{h,\tau}}{1}$ is first order in $h$ or less, if the solutions are not smooth. We confirm this below. The error in $w$ is a different matter, since $w$ is only evaluated as a function of $u$. In fact, this evaluation may incur a local $O(\tau)$ modeling error since the manner of scanning the actual $\hH(\pP;\cdot)$  is dependent on the discretization error. So while 
$\norm{w-W_{h,\tau}}{1}$ remains stably bounded, in our experiments it is $O(1)$. 
The issue is similar to that discussed in Sec.~\ref{sec:err}.

%%%%%%%%%%%
\subsection{Numerical examples for transport with hysteresis}

We illustrate and compare the different hysteresis models for \eqref{eq:pde} with the scheme from Sec~\ref{sec:pdenum}. The operator $\hH(\pP;\cdot)$ is constructed from either the symmetric concave-convex graph $H$ from Sec.~\ref{sec:convex}, or from the adsorption--desorption graph from Sec.~\ref{sec:adsorption}. 
We set $a(u)=u$ and $\alpha(u)=u$ in all examples, and focus on the hysteresis models alone; there is no substantial difficulty, but the  exposition for more general functions $a(\cdot)$ and $\alpha(\cdot)$ takes more time. 

We recall that the convergence rate was demonstrated to be $O(\sqrt{h})$ for \knonlinear\ case in \cite{PeszShow20}, and essentially first order for \klinear\ on examples similar to those presented below. We present here the case with \general\ shown in Table~\ref{tab:errPDE} which reveals that the order in $u$ is $O(h)$ for \general\ model. 
However, even though the solutions $W_{h,\tau}$ appear to visually converge, the error in $w$ scales like $O(1)$.  

%\end{document}
%%%%%%%%%
\begin{table}
\centering
\begin{tabular}{|c|c|llll|c|}
\hline
&h&0.01&0.005&0.001&0.0005&
\\
\hline
$\pgen$, $\uz$ from Fig.~\ref{fig:ch4box}&$E_h$&19.3351&10.0101&2.6699&1.32687&$p=0.89$
\\
\hline
&h&0.05&0.01&0.005&0.001&
\\
\hline
$\pknon$, $\uz$ from Fig.~\ref{fig:ch4box}&$E_h$&63.0521&17.7914&8.4411&1.3292&$p=0.98$
\\
\hline
$\pknon$, $\uz$ from Fig.~\ref{fig:ch4linear}&$E_h$&28.70&2.96&2.72&0.08&
$p=1.4$
\\
\hline
\end{tabular}
\caption{
\label{tab:errPDE}
Error $E_h=\norm{u-U_{h,\tau}}{1}$ for \eqref{eq:PDEFD}, with $\pgen$ and $\pknon$ for the adsorption model from Sec.~\ref{sec:adsorption}, and initial conditions as indicated We use $\tau=0.9h$. The solution for $h_{fine}=0.0001$ and $h_{fine}=0.0005$ is a proxy for $u(x,t)$ for $\pgen$ and $\pknon$ with $K=287$.
The error for $w$ which is stable but appears $O(1)$is not shown. The error in $w$ is as shown, with roughly first order for all cases.}
\end{table}

\myskip{some values when lambda = 0.9
Error for tauc=0.01 is uerr=19.3351 werr=373.737
Error for tauc=0.005 is uerr=10.0101 werr=374.695
Error for tauc=0.001 is uerr=2.66996 werr=375.459
Error for tauc=0.0005 is uerr=1.32687 werr=375.447
}

%%%%%%%%%
\subsubsection{Transport with adsorption--desorption hysteresis from Sec.~\ref{sec:adsorption}}
\label{sec:ch4-transport}
\begin{figure}[ht]
\begin{tabular}{cc}
\includegraphics
[height=50mm]{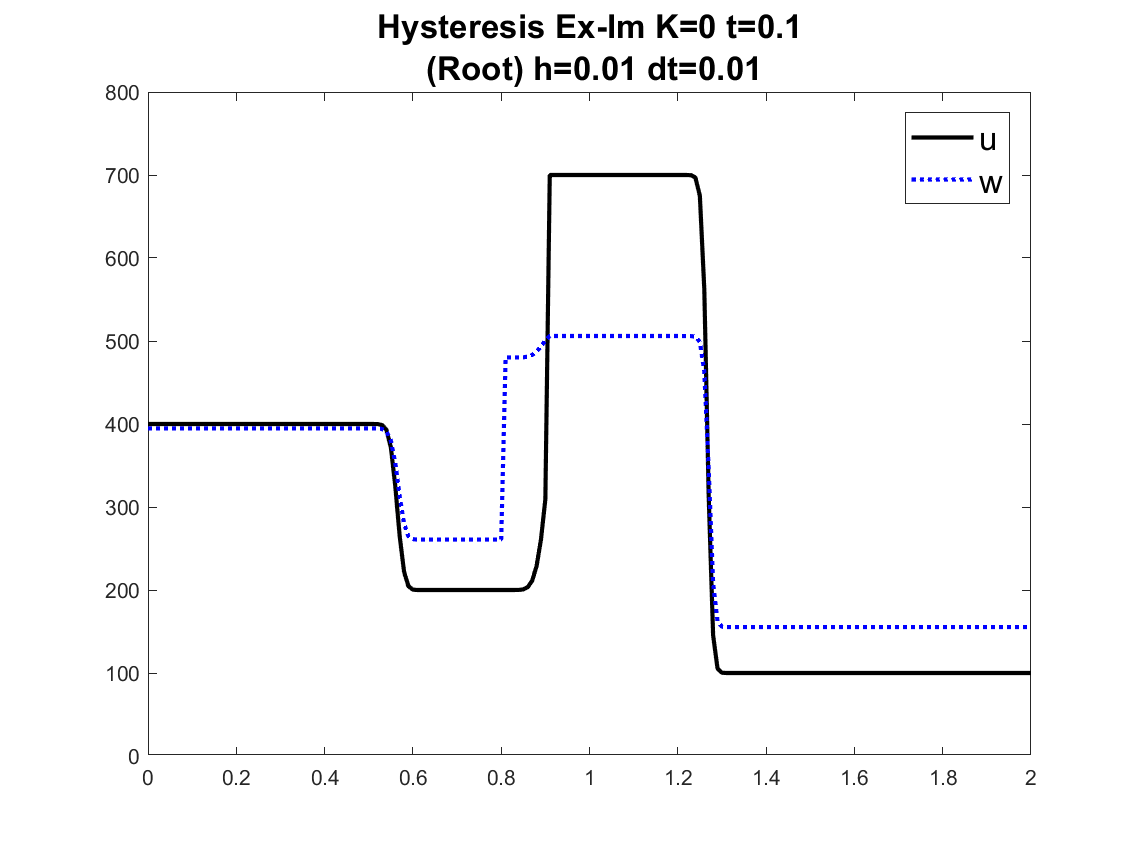}
&
\includegraphics
[height=50mm]{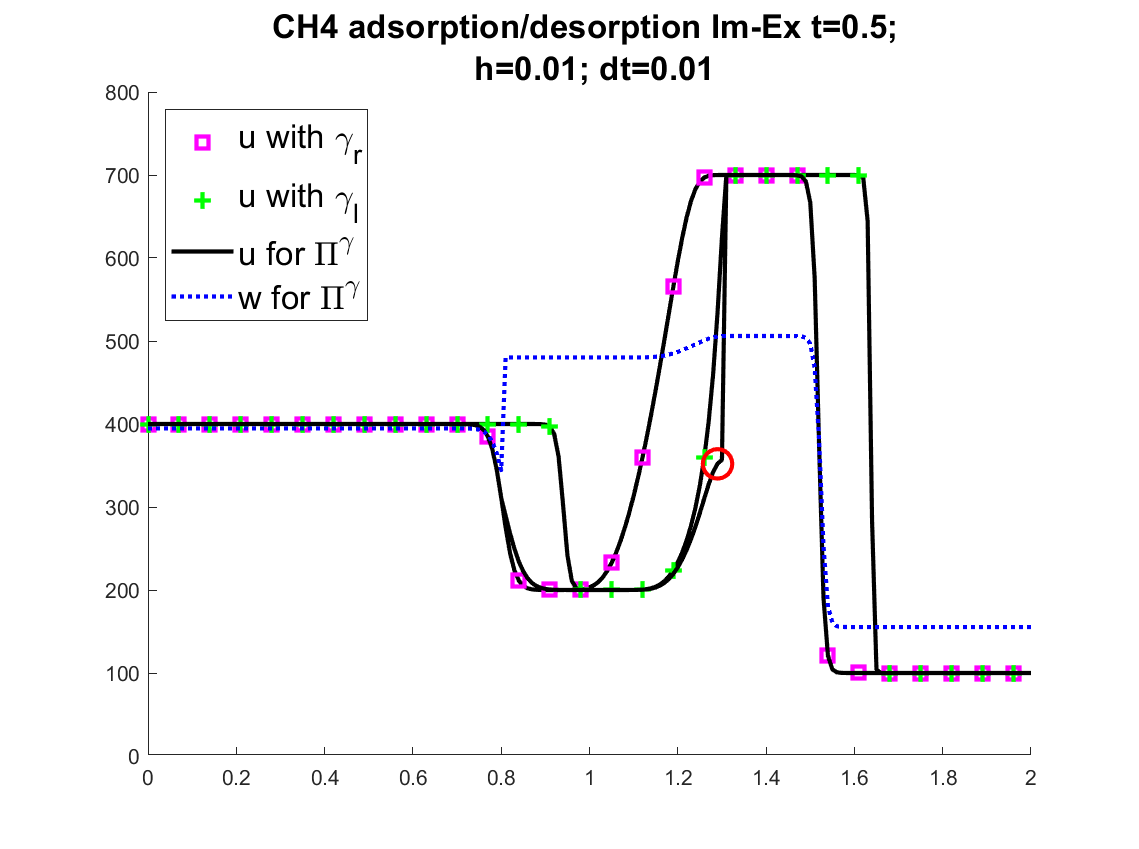}
%\includegraphics
%[height=40mm]{cimages/ch4_finalBOX.png}
%&
\\
(a)&(b)
\\
%%%%
\includegraphics
[height=50mm]{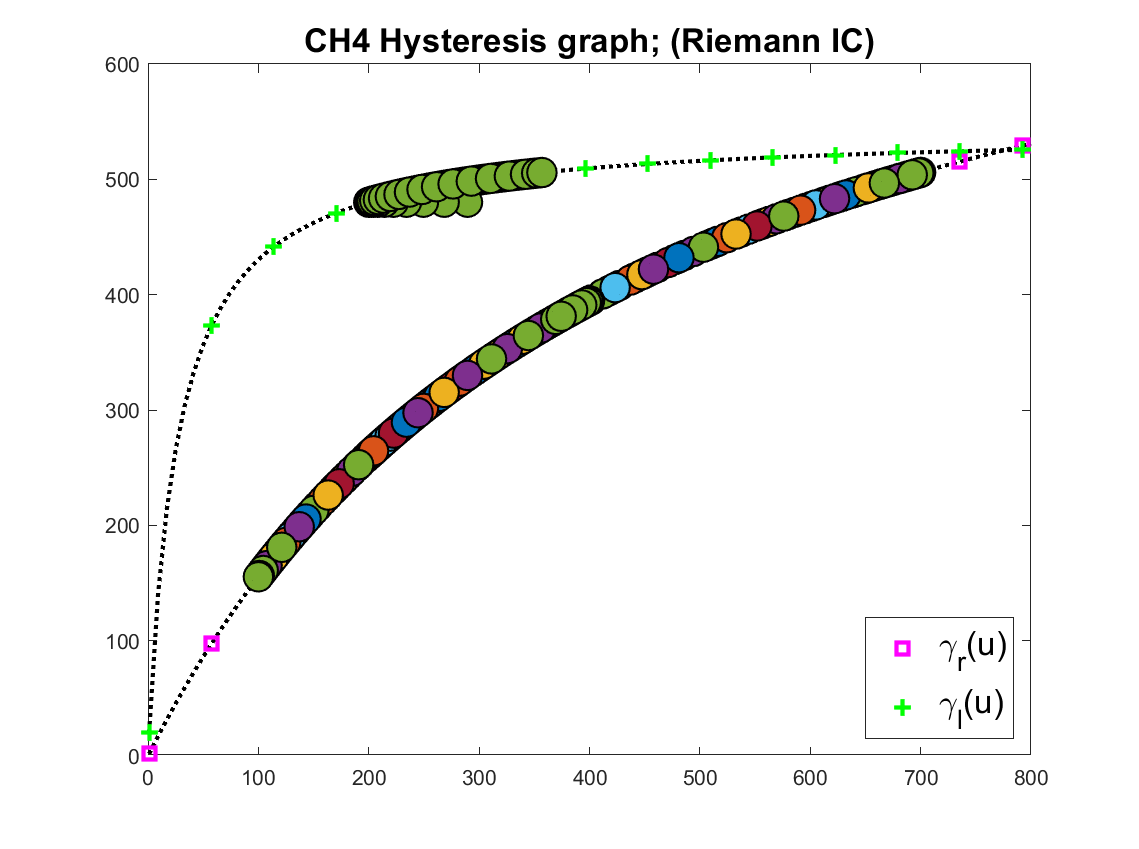}
&
\includegraphics
[height=50mm]{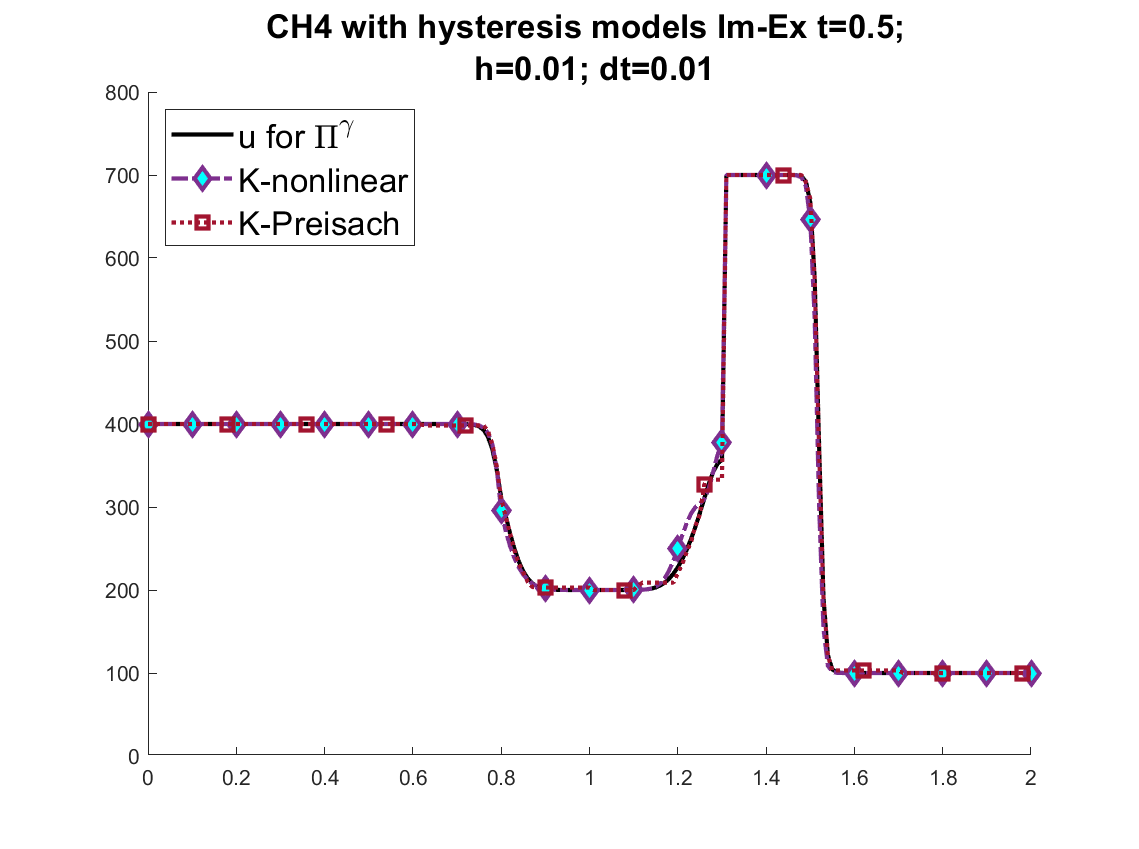}
\\
(c)&(d)
\end{tabular}
\caption{Illustration of transport with adsorption hysteresis models  for the CH4 adsorption-desorption curves from Sec.~\ref{sec:ch4-transport}. Top: (a) initial and (b) final results at $t=0.5$. 
Bottom: (c) the trace $(U_j^n,W_j^n)_{j,n}$ of $\hH(\pP;)$ for $\pP=\pgen$, and (d) comparison of $U_{h,\tau}(\pP)$ for different $\pP$.
\label{fig:ch4box}}
\end{figure}

\begin{figure}[ht]
\begin{tabular}{ccc}
\includegraphics
[height=40mm]{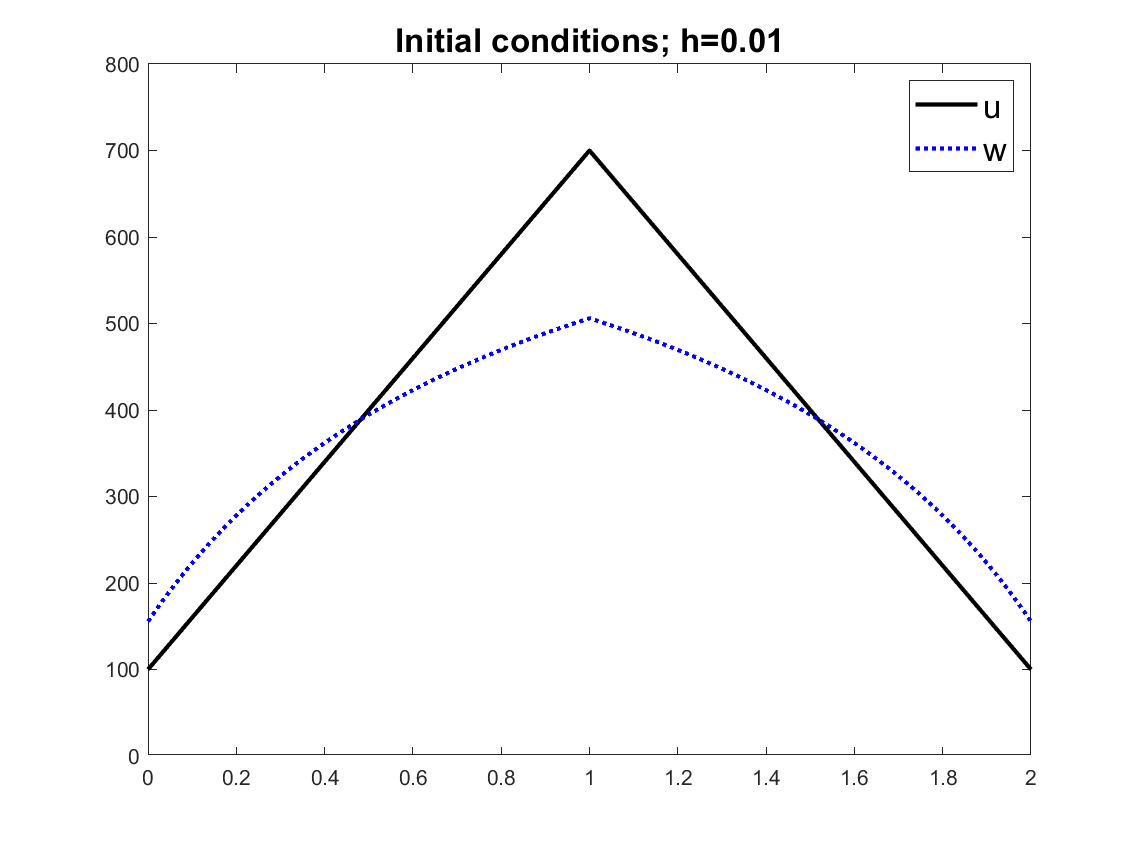}
&
\includegraphics
[height=40mm]{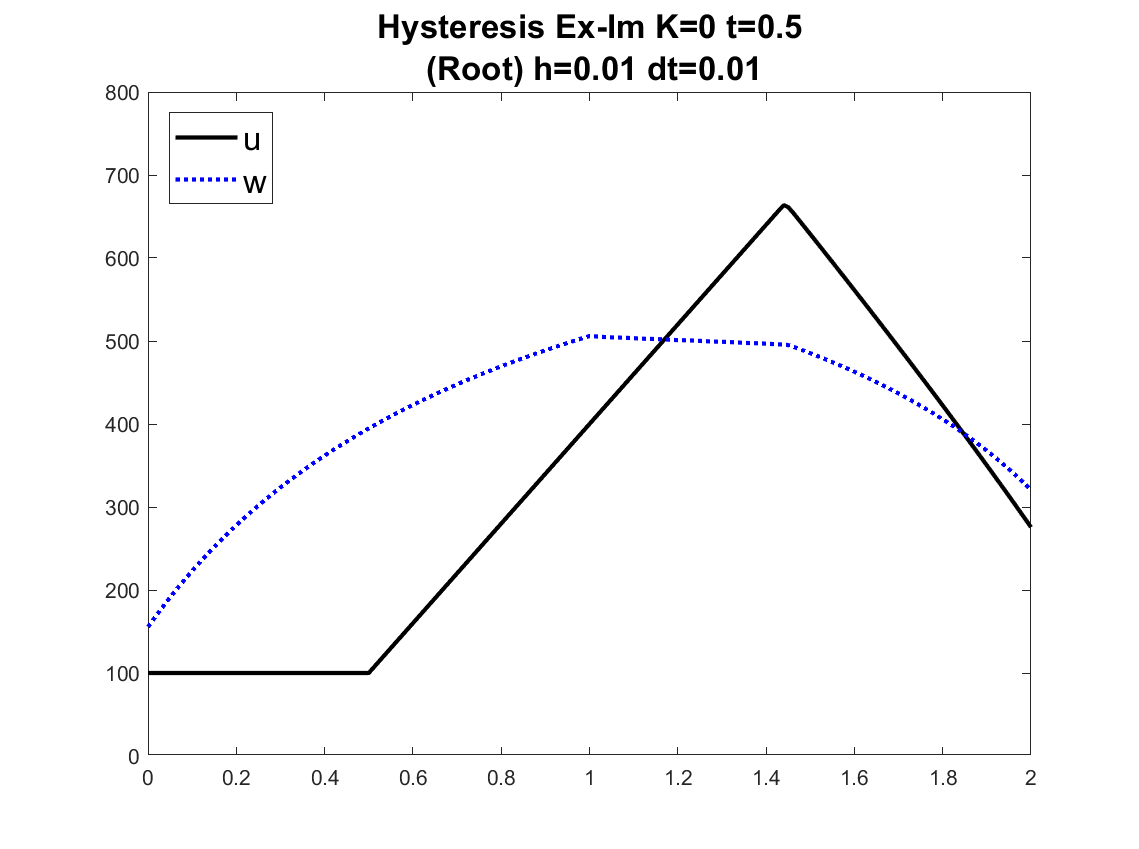}
&
\includegraphics
[height=40mm]{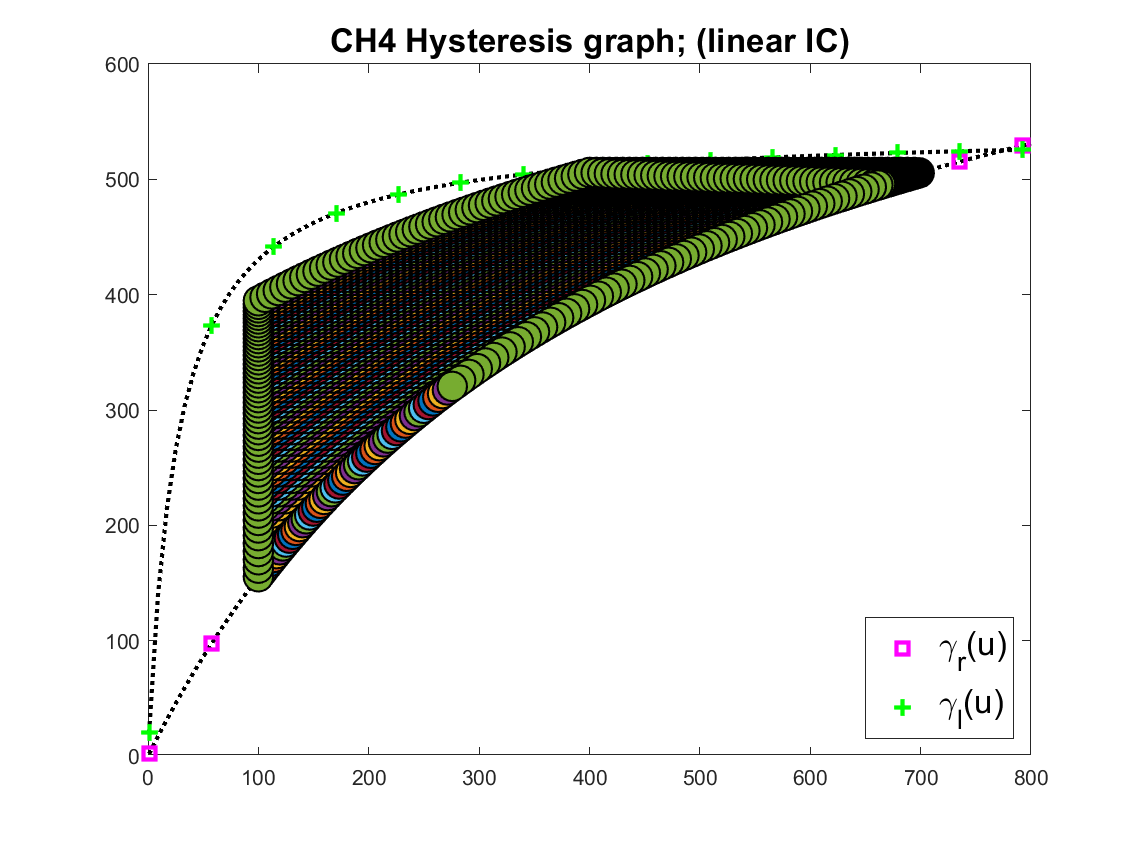}
\\
\end{tabular}
\caption{Simulation of transport with adsorption hysteresis from Sec.~\ref{sec:adsorption}. Left and middle: initial and final $(U_{h,\tau},W_{h,\tau})$ with $\pgen$. Right:  the trace $(U_j^n,W_j^n)_{j,n})$ of $\hH(\pgen,:)$. 
\label{fig:ch4linear}}
\end{figure}

We consider \eqref{eq:ivp} with two different $\vz(x)$ as shown in Fig.~\ref{fig:ch4box} and Fig.~\ref{fig:ch4linear}. The first is made of a superposition of Riemann problems, the second is a simple piecewise linear. 

The first case in Fig.~\ref{fig:ch4box} is designed to show how the solution $U_{h,\tau} \approx u(x,t)$ develops several interesting features. Since $\tau=h$ exactly, there is little diffusion and the fronts are not overly smeared. First, we observe the shock and rarefaction waves forming for the no-hysteresis case when suing $\fl(u)$ and $\fr(u)$. Since both of these are concave, the resulting function is convex, thus the fronts form in the front and rarefaction in the back, with the profile for $\fl(\cdot)$ with steeper slopes travelling faster than that for $\fr(\cdot)$. Then we turn to study the profile of $(U_{h,\tau},W_{h,\tau})$ corresponding to $\pP=\pgen$. Because the back of the first ``box''drops almost discontinuously from $u=700$ down to $u=350$, the resulting $w$ follows a secondary scanning curve with slope $0$, as expected for the generalized hysteresis. Thus the resulting flux function is linear, and thus this part of the graph does not develop a smooth rarefaction, unlike for $\fl$ or $\fr$. For the input between $200\leq \leq 300$, $w$ follows again $\fl(w)$ and the rarefaction develops again.   The ``trace'' of the graph $(U_j^n,W_j^n)_{j,n})_{j,n}$ over all $j$ and $n$ from the entire simulation shows rather sparse sampling of $H$ due to discontinuous profile of the solution. 

The second case in Fig.~\ref{fig:ch4linear} leads to a much richer trace of $\hH(\pP)$  which results from a large collection of intermediate values $U_j^n$ in the front and back of the wave. The front travels with velocities found from $\fr$, and the back with velocities from $\fl$. This profile eventually steepens and becomes a shock followed by rarefaction (not shown).  

%%%%%%%%%%%

%%%%%%%%%%%

%%%%%%%%%%
\subsubsection{Transport with convex-concave graph $\hH$ from Sec.~\ref{sec:convex}.}

We set-up the initial condition to be the ``trough'' (``well''), and simulate with the different hysteresis models. We also compare the simulation with hysteresis to that without. The latter examples show the expected behavior of the sides of the ``well'. If only the (convex) $\fr$ is used, the flux function, the inverse of $u+\fr(u)$, is concave, thus we expect a sharp front on the increasing right hand side of the well. The opposite happens when $\fl(\cdot)$ is used. 

When the hysteresis model is used with one of $\pklin, \pknon,\ppreisache$, both sides of the graph show behavior typical of rarefaction which arises because of concavity of $\fr$ on the increasing side, and the convexity of $\fl$ on the decreasing side. 

The results for all models are qualitatively consistent with this description and with each other, and the models \knonlinear, \kpreisach, and \kgeneral\ all give very similar resuts, with the $u(\ppreisach; x,t)$ corresponding to the \kpreisach\ being the most rough. The graphs $\hH(\pP;\cdot)$ are very similar to each other. 

The results with \klinear\ model have the most rich secondary curves, thus the values of $u(x,t)$ corresponding to the bottom of the ``well'' travel with a larger variety of velocities than those for $\pknon$.

%%%%%%%%%%%
\begin{figure}[ht]
\centering
\begin{tabular}{cc}
\includegraphics
[height=40mm]{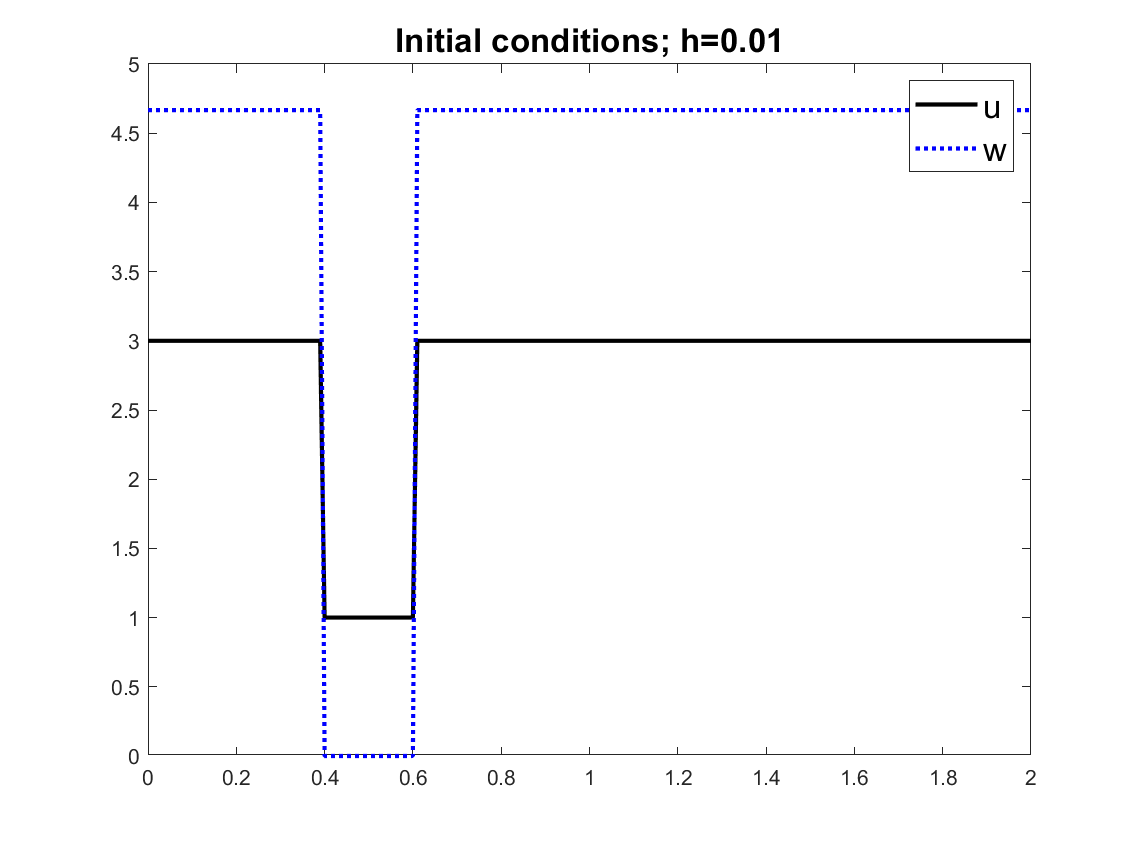}
&
\includegraphics
[height=40mm]{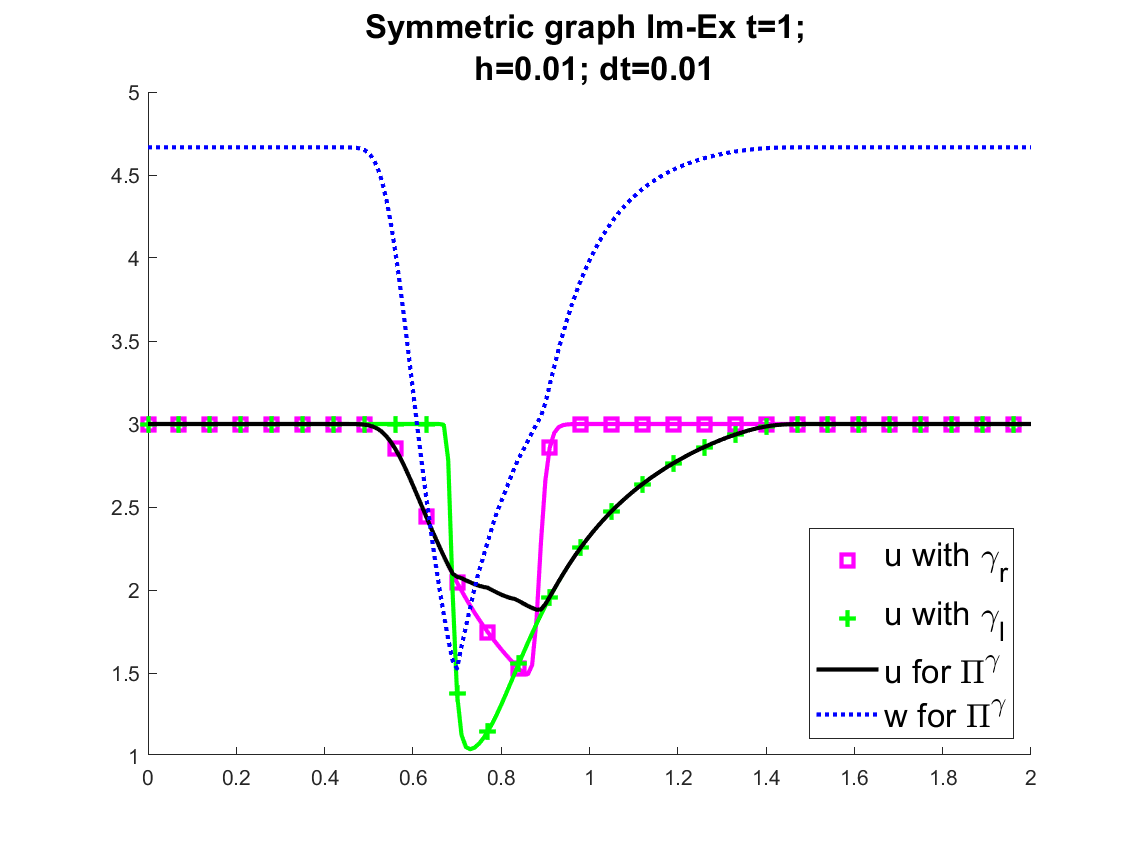}
%%%%
\\
(a)&(b)
\\
\includegraphics
[height=40mm]{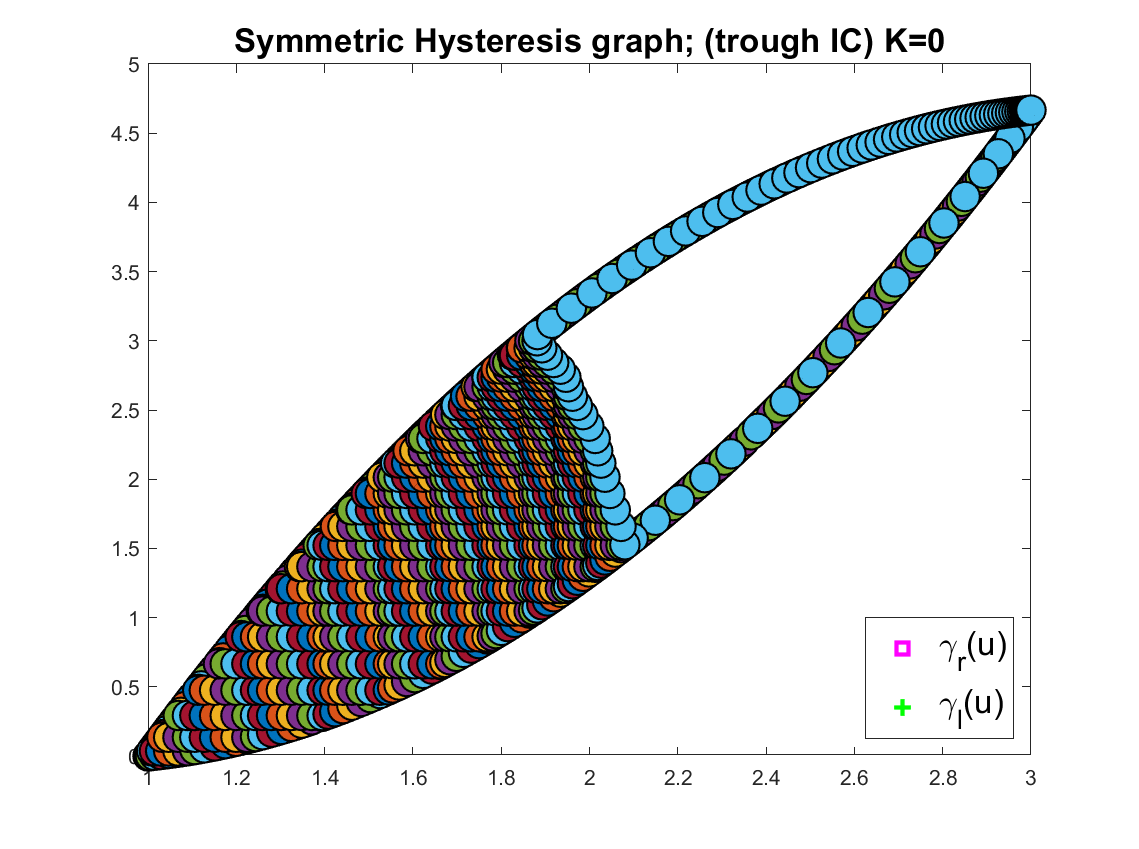}
&
\includegraphics
[height=40mm]{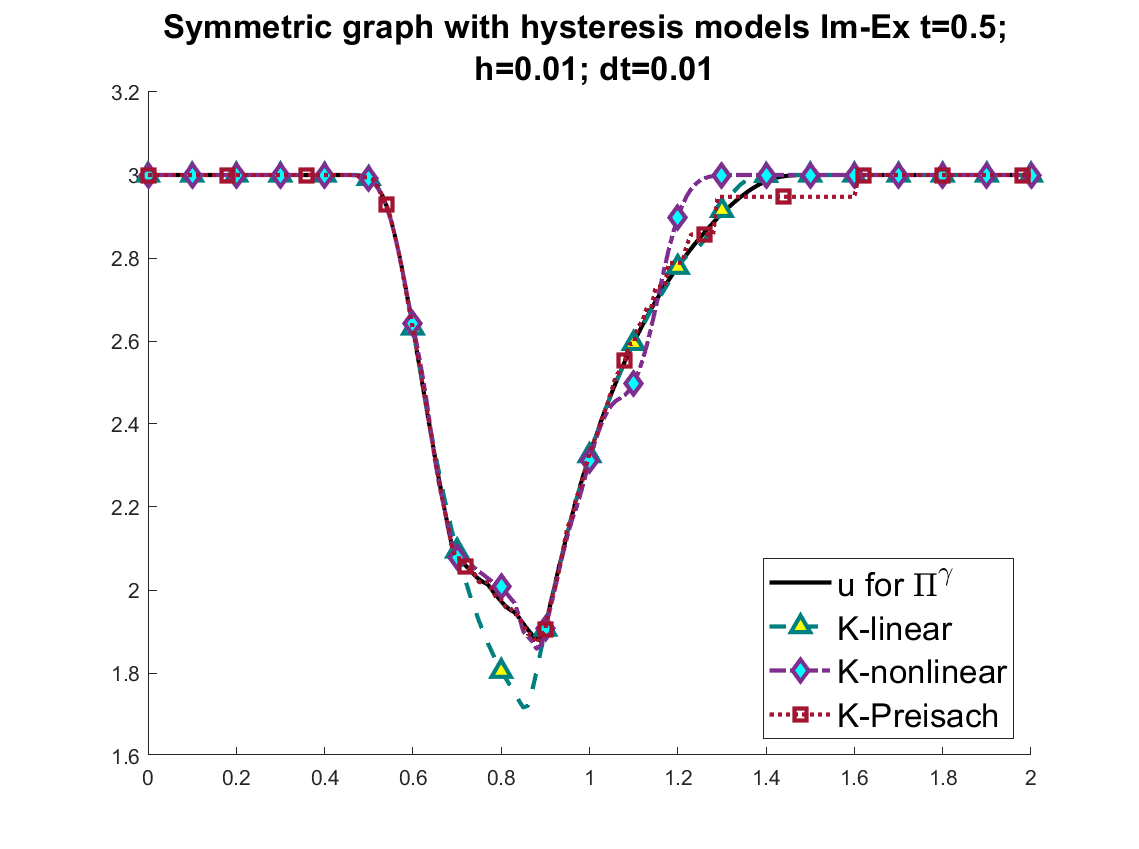}
\\
9c)&(d)\\
\multicolumn{2}{c}{
\includegraphics
[height=30mm]{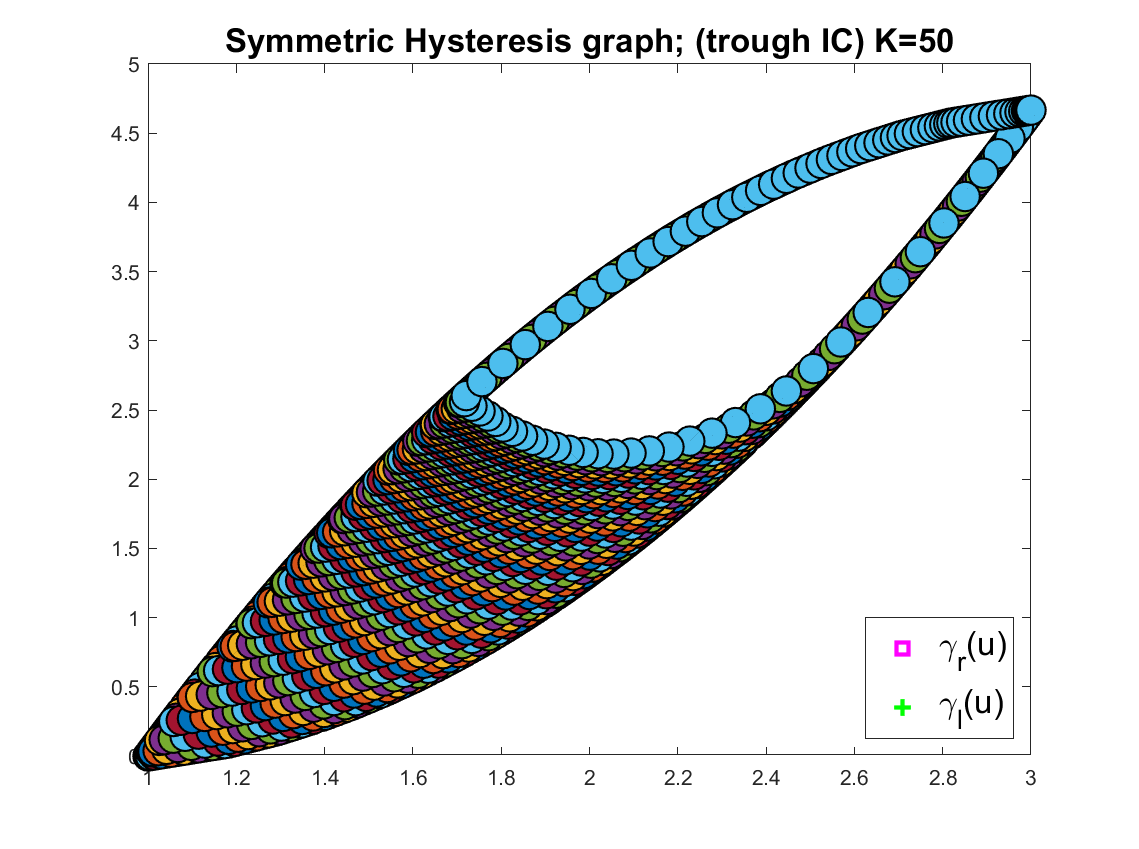}
\includegraphics
[height=30mm]{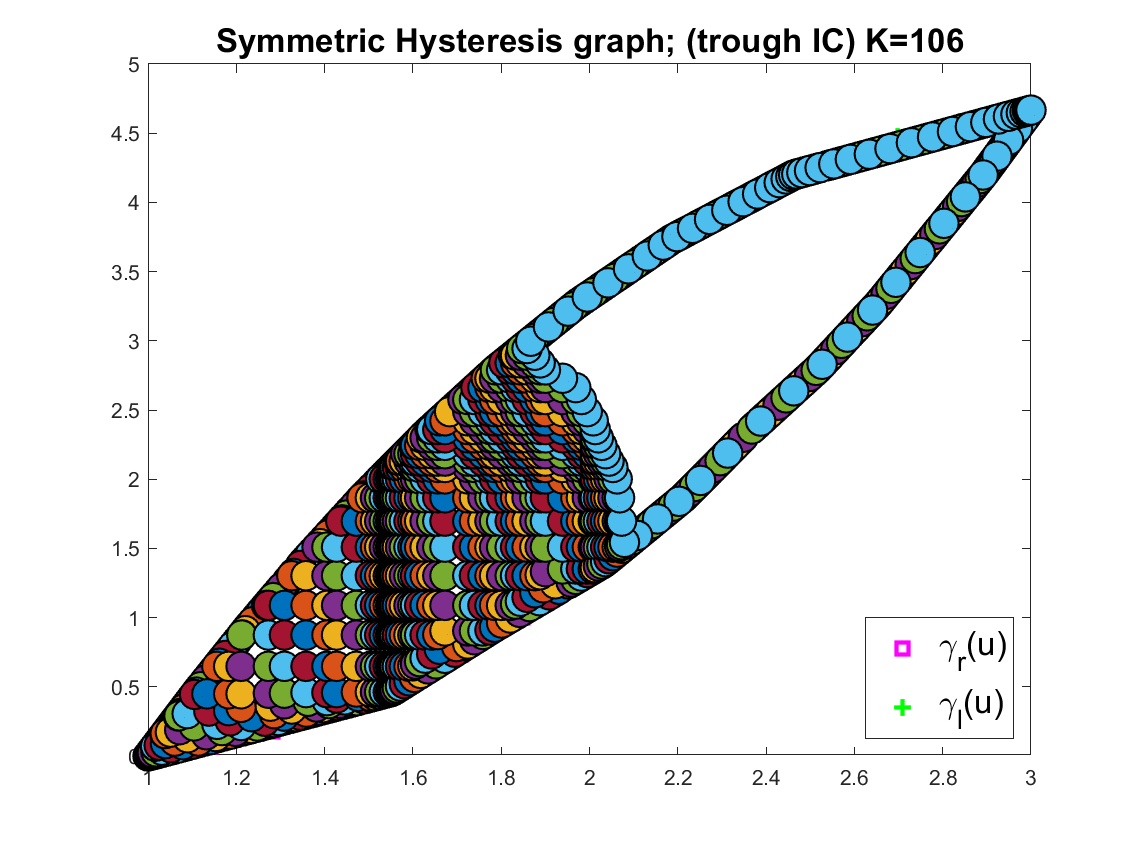}
\includegraphics
[height=30mm]{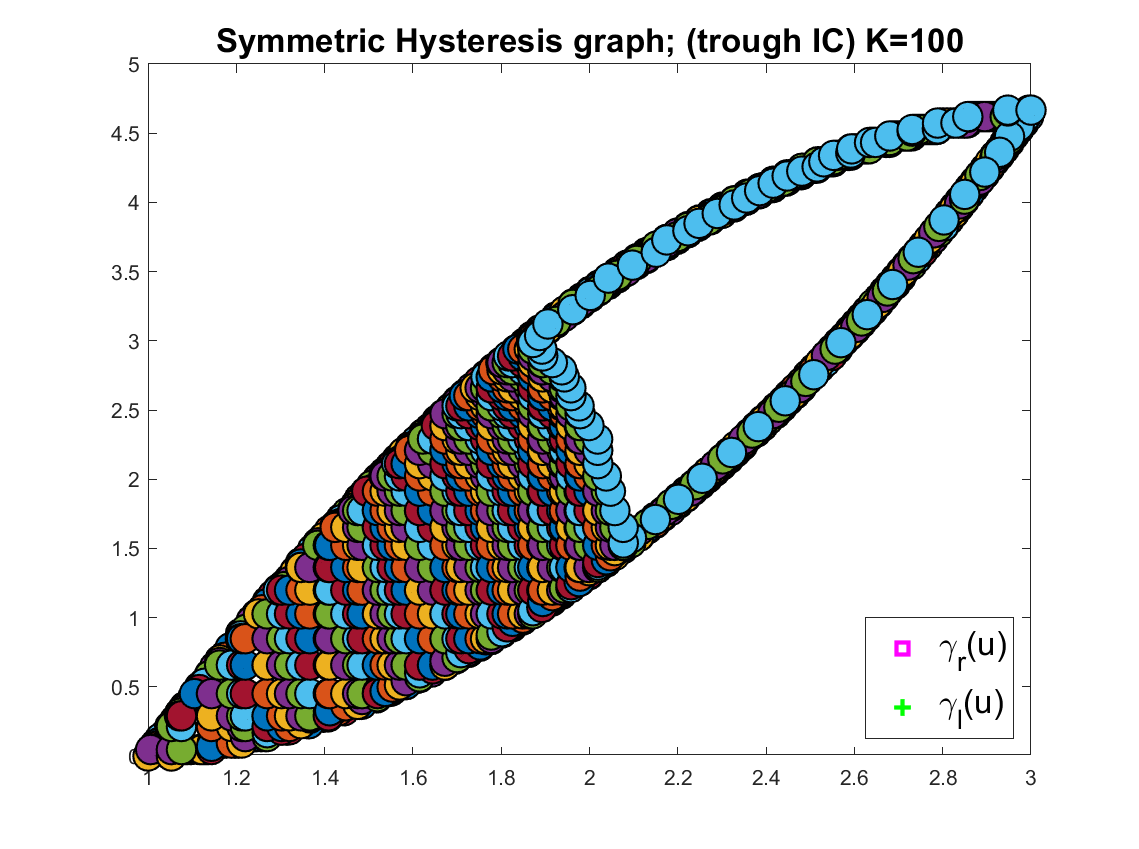}
}
\\
\multicolumn{2}{c}{
(e)\phantom{ala ala ala ala ala ala} (f)\phantom{ala ala ala ala ala} (g)
}
\\
\end{tabular}
\caption{Simulation of adsorption hysteresis models  for the graph from Sec.~\ref{sec:adsorption} at (a) $t=0$ and (b) $t=.5$ obtained with $\pgen$. (c) Trace of $\hH$ from $(U_j^n,W_j^n)_{j,n}$. (d) Results at $t=.5$ from $\pgen,\pklin,\pknon,\ppreisache$. Bottom: the traces $\hH(\pP;)$ for $\pklin,\pknon,\ppreisache$ with $K=50$, $K=106$, $K=100$, respectively.  
\label{fig:ch4}}
\end{figure}

%%%%
\subsection{Complexity of solving transport PDE with hysteresis and extensions to other temporal discretizations}
\label{sec:complexity}

\mpcomment{
It remains for us to state the cost of accounting for hysteresis when solving the transport PDE. }
\mpcomment{
Without hysteresis, at each time step when solving \eqref{eq:PDEFD} we need to find $U$ from the value $a(U_j^n)$ calculated explicitly from the previous the time step. This may take a few iterations of a nonlinear solver at every $j$, or require the use of a lookup table for $a^{-1}$.  Generally we can say this cost is  $f$ flops per each $j$ and $n$, with a rough estimate of $2<f<20$ depending on the approach taken. 
}

\mpcomment{
With hysteresis, the cost of using \knonlinear\ graph depends on the number $K$ of components. Here we must find $U_j^n$, $W_j^n$ as well as all components $V_{j,k}^n,k=1,\ldots, K$. These are found by solving the pointwise ODE \eqref{eq:coupFD} posed at every $j$. The resulting local nonlinear system is solved iteratively, at the following cost.  In each iteration we calculate the resolvent in \eqref{eq:coupFDw} which is just an algebraic formula, with the cost we estimate of $\approx 5K$ per each $j$. The number of iterations is usually about two depending on the solver, thus the additional cost per each $j$ is about $\approx 10K$ per each $j$ and $n$. For the simple graph considered in Fig.~\ref{fig:example:PDE} this amounts to about $f \approx 20$ flops, but for a complex graph such as in Fig.~\ref{fig:calibrate}, the cost might be $f\approx 200$  when $K=42$ or $f\approx 1000$ when $K=287$ per each $j$ and $n$. The modeling precision comes with higher {$K$} but also higher cost. 
}

\mpcomment{
Refinements of the solver are possible. Clearly one can think of introducing clever refinements of time stepping and solver such as adaptivity, since not all components $k=1,\ldots, K$ are  active at every point $j$ and $n$. 
}

\mpcomment{Lastly, we discuss the possibility of using other than fully implicit first order approaches proposed in this paper. Extensions to high order and refinements are clearly possible. However, the use of non-implicit approaches requires regularization of the component graphs which have high Lipschitz constants and require very small time-stepping for accuracy. In turn, higher order schemes are possible but may have limitations due to low temporal regularity of solutions,  even away from shocks. 
}
%% solvers 

%% different temporal discretization

%%%%%%%%%%%%%%%%%%%%%%%%%%%%%%%%%%%%%%%%%%%%%%%%%%%%%%%%%%%%%%%%%%%%%%%%%%%%%%%%%

\section{Towards calibration with secondary curves}
\label{sec:secondary}
As we discussed above, the model \knonlinear\ provides an easy opportunity to enhance the parametrization to account for secondary scanning curves. This is similar to the modeling power of the Preisach model. While a thorough discussion is outside our present scope, we provide a simple example to illustrate the power of $\hH(\pknon;\cdot)$. 

\begin{figure}[ht]
\begin{tabular}{cccc}
\includegraphics
[height=40mm]{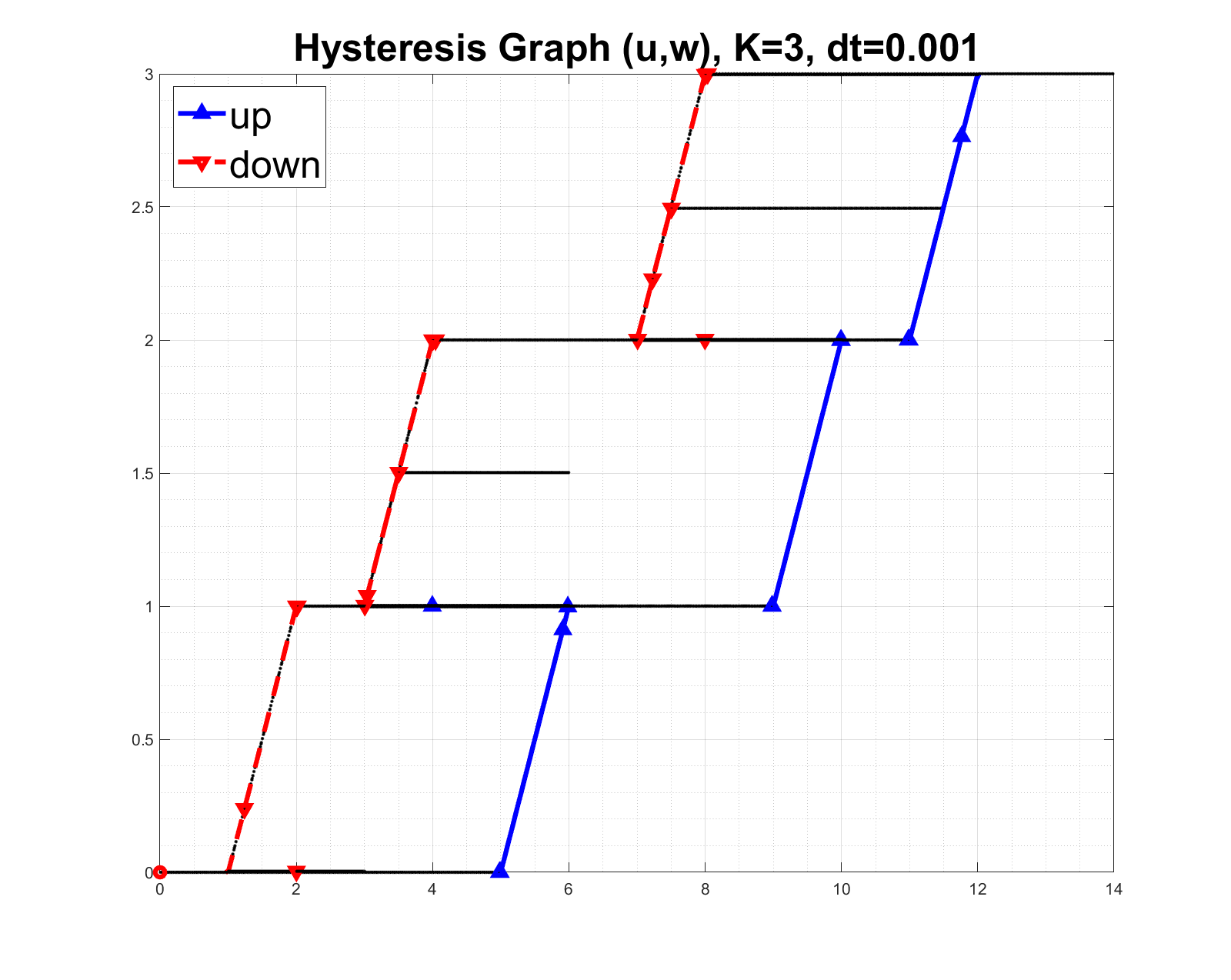}
&
\includegraphics
[height=40mm]{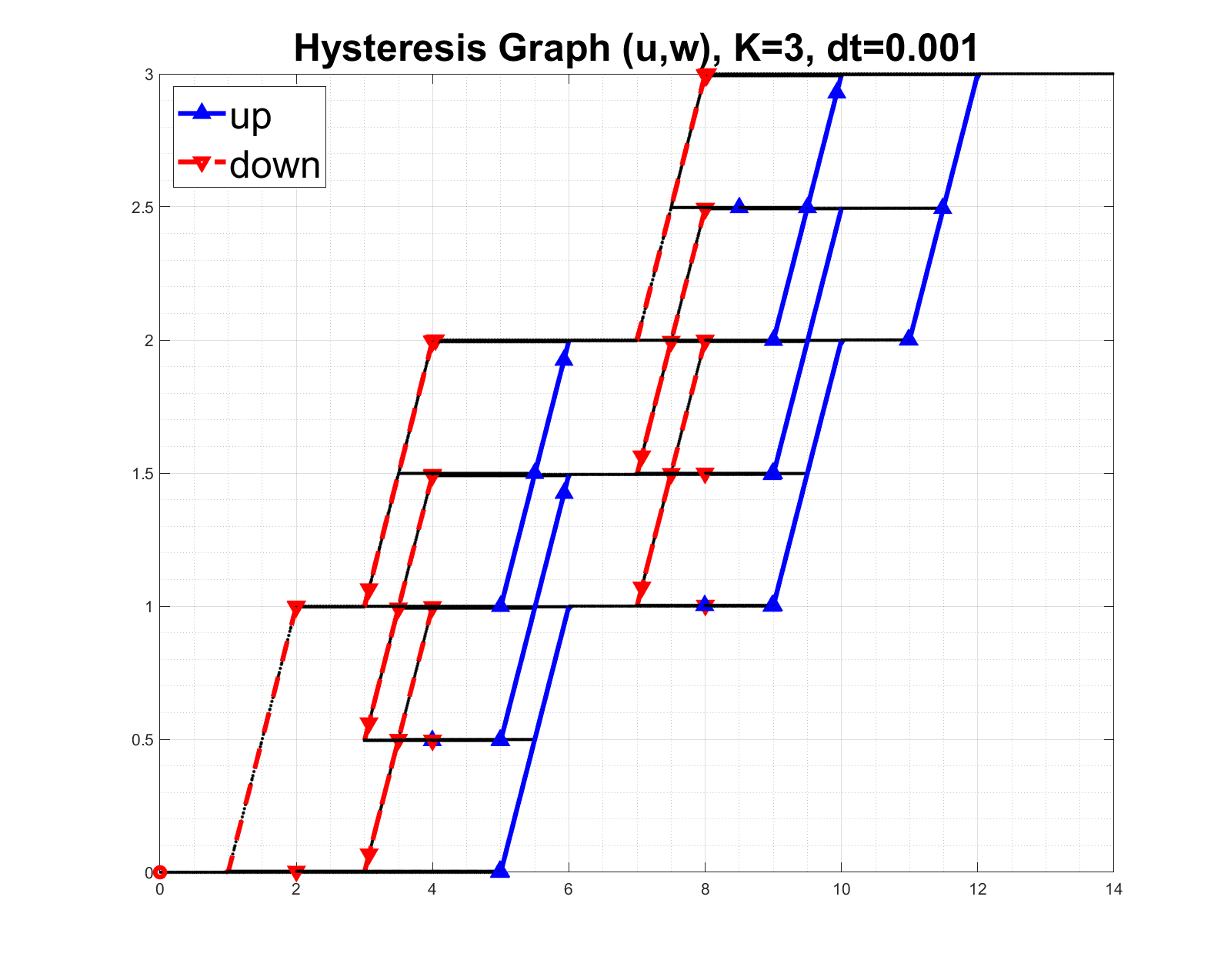}
&
\includegraphics
[height=40mm]{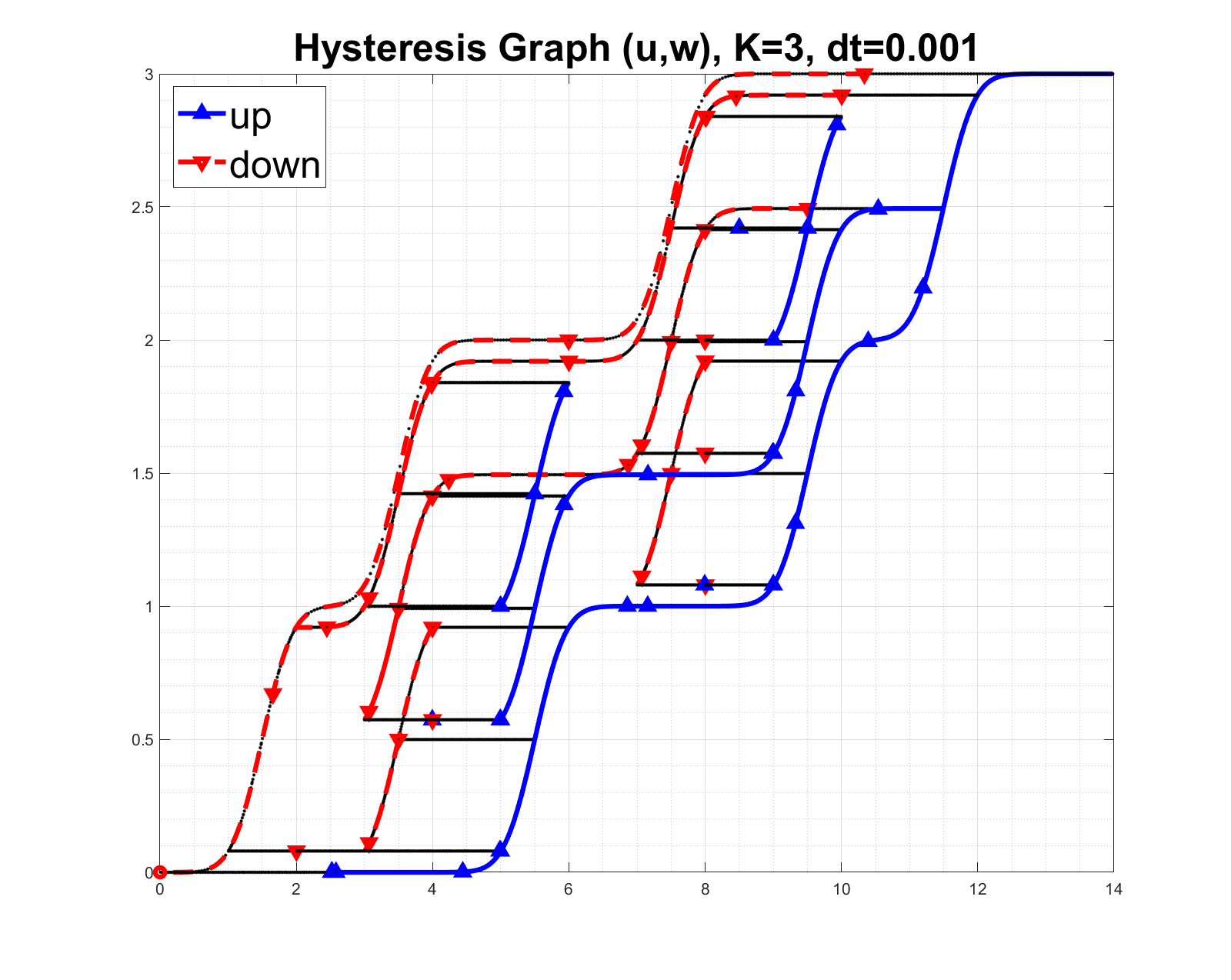}
\\
$\pP^{mon}$&
$\pP^{rich}$&
$\pP^{rich,*}$
\end{tabular}
\caption{Primary and secondary scanning curves for $\pknon$ from \eqref{eq:secondary} from Sec.~\ref{sec:secondary}. 
\label{fig:secondary}}
\end{figure}

Consider three parametrizations of the same primary scanning curves $H$ made of the sides of three parallelograms and shown in Fig.~\ref{fig:secondary}. The set $H$ can be obtained with $\uU^{primary}=[0,14,0]$ so that $H=\hH(\pP;u)$ for $u \in \PL(\uU^{primary})$ 
\bsub
\label{eq:secondary}
\ba
\pP^{mon}
 %twooloopa
 &=& [1,1,5,1;1,3,9,1;1,7,11,1];
 \\
\pP^{rich}
 %twoloop 
 &=& [1,3,5,1;1,7,9,1;1,1,11,1];
 \\
\pP^{rich,*} &=&[1,3,5,\bh_*;1,7,9,\bh_*;1,1,11,\bh_*]; h=1.
 \ea
\esub
with the last using the smooth truncation function \eqref{eq:bsmooth}.  Next we design the input $u\in \PL(\uU^{rich})$ to produce a rich variety of secondary scanning curves, and plot $\hH(\pP;u)$; here
%%%
\bas
\uU^{rich}=
%% utwoloop = 
[0,14,0,6,3,5.5,3,10,7,9.5,7,11.5,3,6,3,5.5,3,10,7,9.5,7,12,7,10,7.5,10,3.5,6,3,6,1,3]
\eas
The results plotted in Fig.~\ref{fig:secondary} show significant difference in secondary curves, without much additional computational effort.

%%%%%%%%%%%

%%%%%%%%%%
\begin{figure}[ht]
\begin{tabular}{cccc}
\includegraphics
[height=30mm]{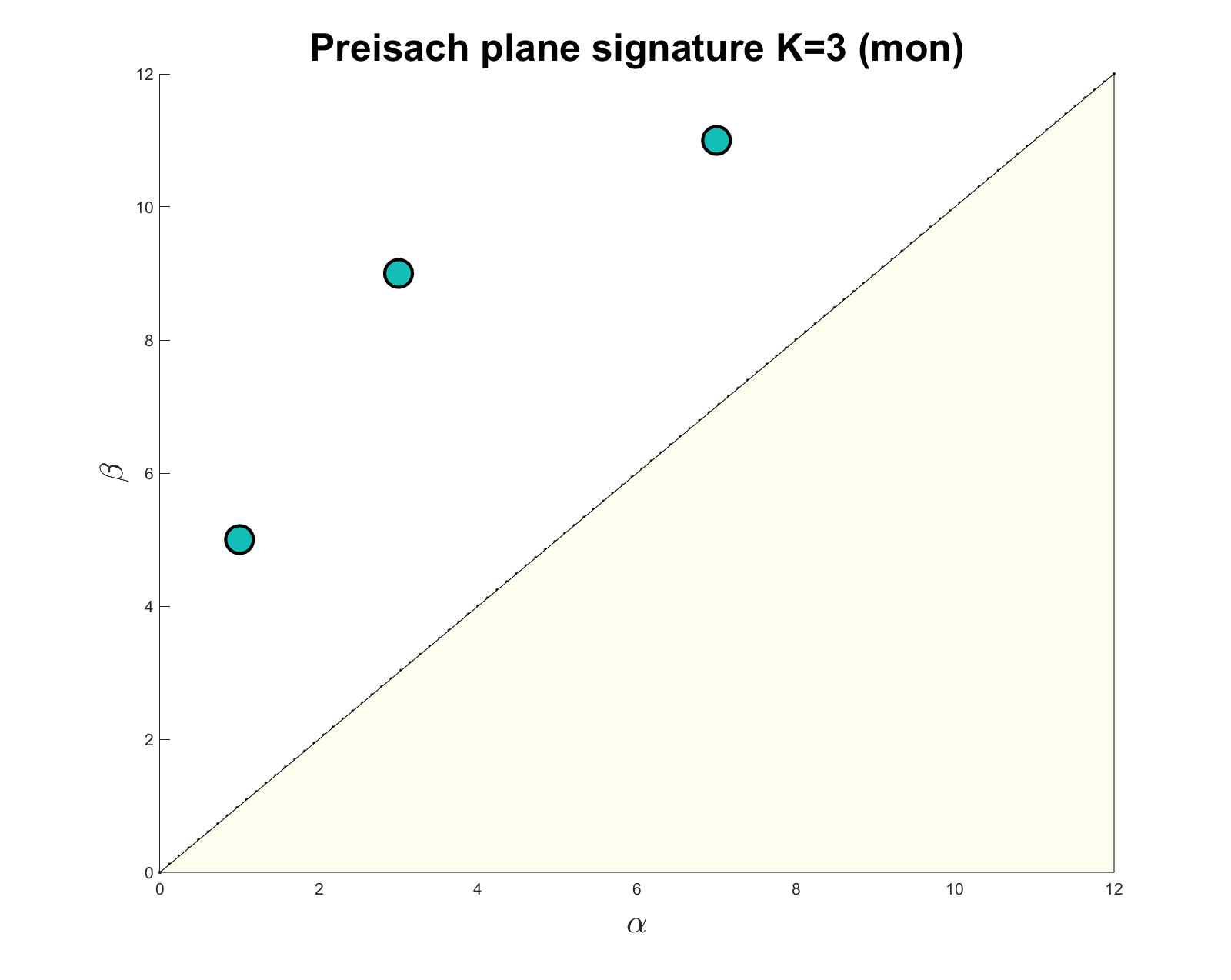}
&
\includegraphics
[height=30mm]{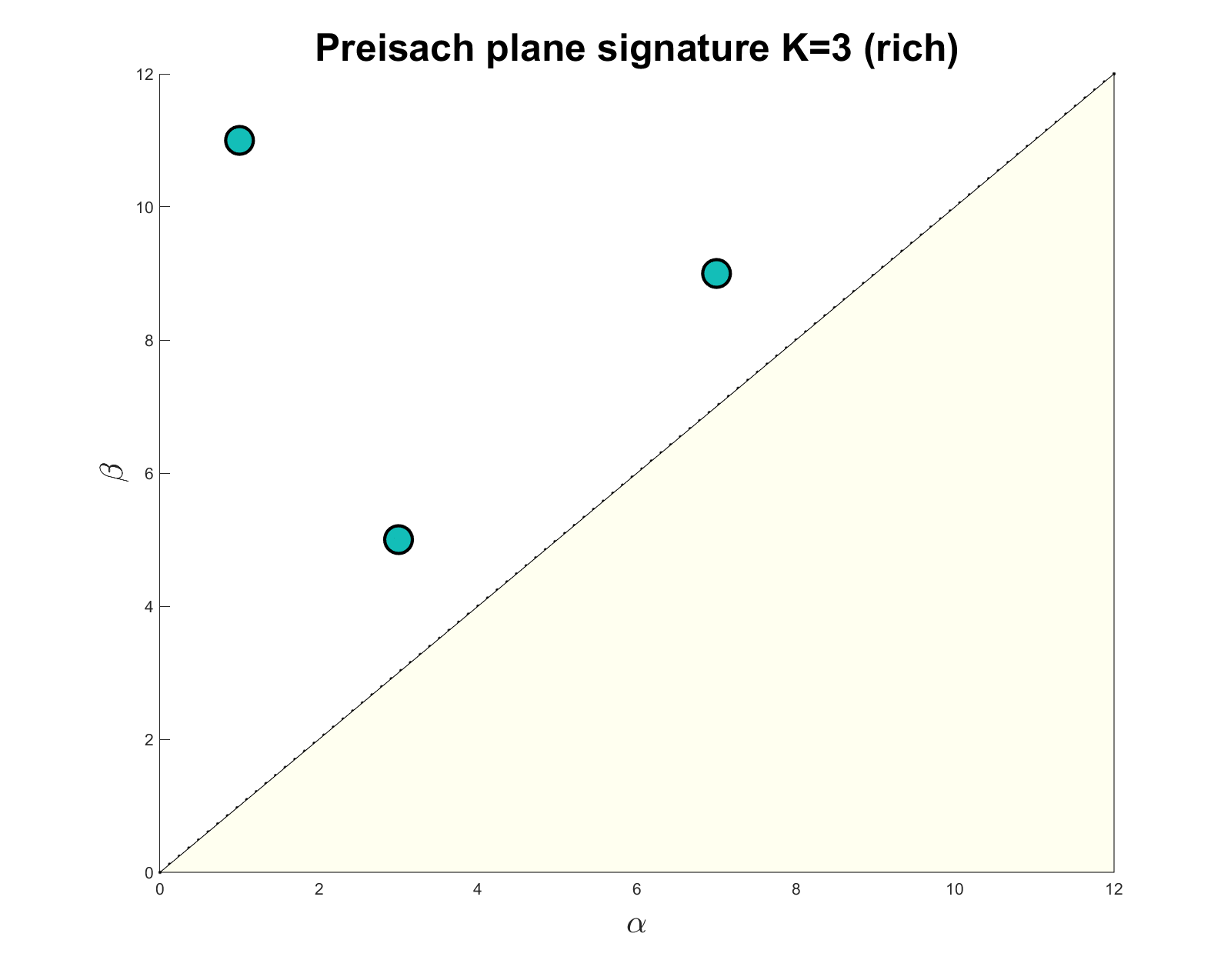}
&
\includegraphics
[height=30mm]{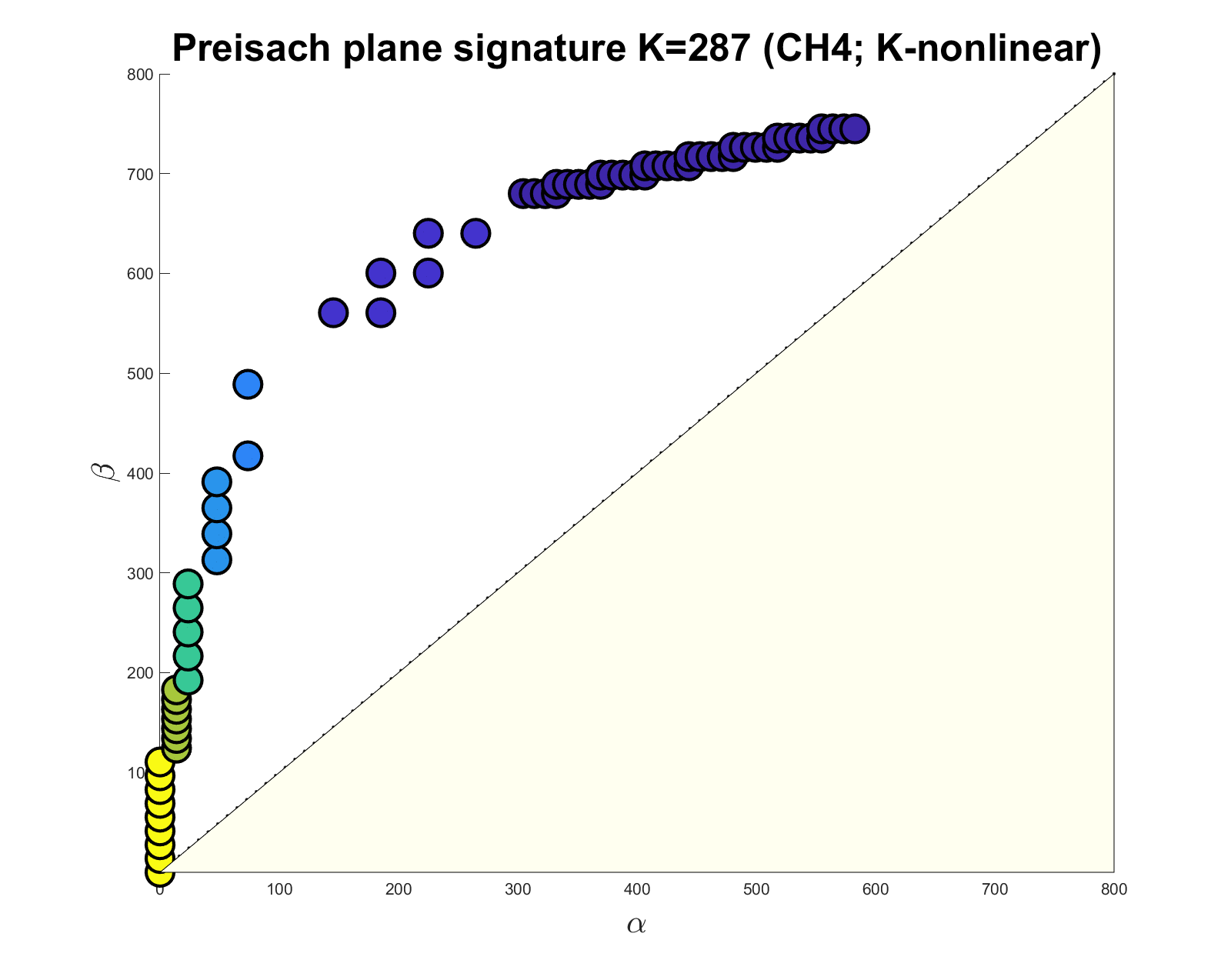}
&
\includegraphics
[height=30mm]{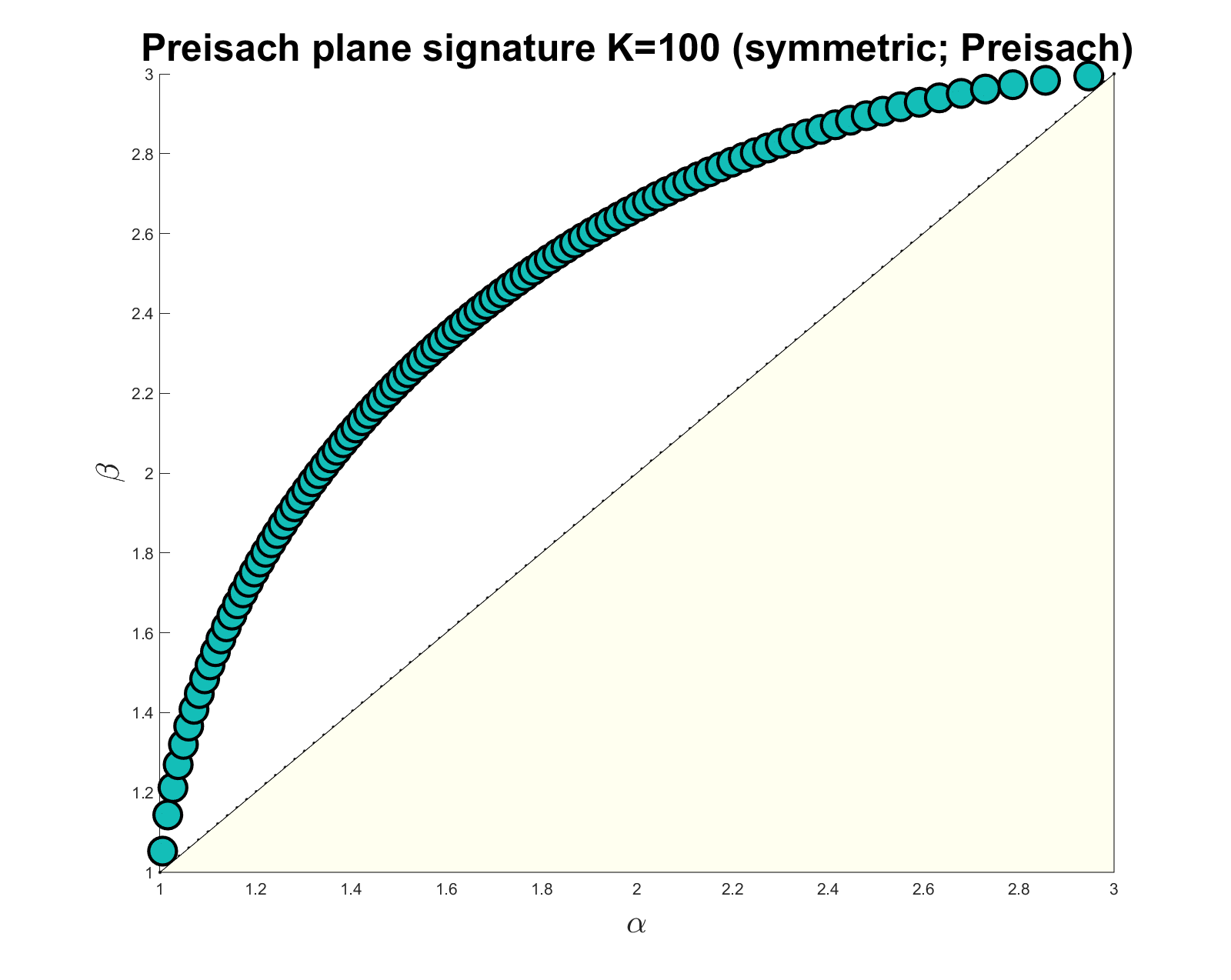}
\\
$\pP^{mon}$&
$\pP^{rich}$&
$\pP^{CH4,K}, K=287$
&
$\pP^{convex,K},K=100$
\end{tabular}
\caption{Preisach plane signatures for $\pP^{mon}$, and $\pP^{rich}$ from Sec.~\ref{sec:secondary}, as well as for $\pknon$ and $\ppreisache$ found  in Sec.~\ref{sec:adsorption} and \ref{sec:convex}. Each illustration presents the collection $(\alpha_k,\beta_k)_k$ in the upper half of the Preisach plane above the line $\beta=\alpha$. The points are colored by the measure $\mu_k$. 
\label{fig:preisach}}
\end{figure}

%%%%%%%%%%%
\section{Summary}
\label{sec:summary}
In this paper we presented a practical view of modeling hysteresis functionals in the case when only limited data is available. In particular, we showed how to calibrate hysteresis graphs of \knonlinear\ type when only the data for primary scanning curves $H$ are available. We compared the use of \nonlinear\ to the generalized play model, and showed that each has advantages and disadvantages, while they all share the theoretical (convergence) properties.  

In particular, (i) \general\ is amenable to the same numerical and well-posedness analysis as \knonlinear. Since it uses one auxiliary equation only, it requires less computational time than \knonlinear, and no calibration. However, by design, it only features horizontal secondary scanning curves. For graphs $\hH$ with a small gap $\fl(u)-\fr(u)$, this may not be significant. In addition, implementation of algorithms involving $\pgen$ requires the use of functions in \eqref{eq:vsol} rather than parameters  as in \eqref{eq:resab}, and is more disruptive to the PDE approximation code and more difficult than that of the models from $\pknon$ or $\ppreisach$ type.  

In turn, (ii) the power and modeling potential of \knonlinear\ is evident from the different examples and from Sec.~\ref{sec:secondary}; this is not surprising since \knonlinear\ includes the discrete version of the Preisach model regularized with $\bheps$. Calibration of \knonlinear\ models can be done with the algorithm provided, and the efforts are not more difficult than approximation of $\fl$ by a piecewise linear function. If $K$ is too large and modeling error not too important, one can allow small flat portions of $H^* \approx H$ and use the simple Preisach model $\ppreisache$ with smaller $K$ instead. 

Our current and future work involves developing further insights into convergence analysis and smoothness, as well as on algorithms for fully implicit schemes for the \kgeneral\ family and when, e.g., operator $A$ is a diffusion operator, as well as when $a(u)$ is less than  strongly monotone. We are also exploring the connection between parametrization with $\pknon$ using secondary scanning curves, and the Preisach plane.

%%%%%%%%%%%%%%%%%%%%%%%%%%%%%%%%%%%%%%%%%%%%%%
\section*{Acknowledgements}
\mpcomment{
We wish to thank the Editor and the anonymous referees for the helpful suggestions which improved this manuscript. 
}
%%%%%%%%%%%%%%%%%%%%%%%%%%%%%%%%%%%%%%%%%%%%%%
\section*{Appendix}

\myskip{
\subsection{Generalization of \kgeneral\  to include graphs \myplan{RES}}

\mpcomment{Use previous argument to extend to minimal sections of maximal monotone $b_i^0(\cdot)$: approximate with Lipschitz monotone Yosida approximations.}
}
%%%%%%%%%%%%%%%%%%%%%%%
\myskip{
\subsection{Minimal selection notes}
%%%%
For the discrete problem, another way to understand the solvability and uniqueness of \eqref{eq:fd}, is
when $\cg$ is the subgradient of a proper convex lower-semicontinuous function $C$, is to see that $V^n=\vs$, the minimizer 
\ba
\label{eq:argmin}
\vs=\ \mathrm{argmin}\ G(v);\;\; G(v)\eqdef \frac{1}{2\tau}\norm{v-V^{n-1}}{}^2+C(v)-(F^n,v).
\ea
%%%
The minimizer in \eqref{eq:argmin} exists because $G(v)$ is convex lsc. 
At the minimizer $v^*$ we have that $\partial G (v^*) =\frac{v^*-V^{n-1}}{\tau}+\cg(v^*)-F^n \ni 0$ which is the same as \eqref{eq:fd}. This point of view is applied, e.g., in analysis of generalization of \eqref{eq:ode} in [Bassetti'2003] as an extension of the method of ``minimizing movements'' considered by DeGiorgi. 

For purposes of illustration, we show that $V^n=\vs$ and is the element of $f-\cg(V^n)$ of minimal norm, when $\cg$ is the constraint graph, in Sec.~\ref{sec:cg}. 

%%%%%%%%%%%%%%  skipping 
\subsubsection{Minimizing movements in discrete approximation of an ODE with a constraint graph}
\label{sec:cg}
We illustrate how the constraint graph $\cab$ works by showing that the finite difference solution $V^n$ given by \eqref{eq:resolvent} is the same as the minimizer $\vs$ given in  \eqref{eq:argmin}, and is the the finite difference analogue of the minimal section characterization in Sec.~\ref{sec:minimal}. Namely, let $V^{n-1} \in \Dom(\cab)=[a,b]$ and let some $F^n \in \R$ be given. We show that $v^* \in \Dom(\cab)$ is the same as the selection $g^*=(F^n-\cab{\vs})^0$ of minimal norm out of $F^n-\cab{\vs}$ so that $\tfrac{\vs-V^{n-1}}{\tau}=g^*$. 

Since $F^n \in \R$, it is either nonnegative, nonpositive, or trivial. Wlog we consider $F^n>0$, and calculate the trial value $\vit=\tau F^n+V^{n-1} > a$. Now if $\vit<b$, then $\cab{\vit} =\{0\}$ and the only element from $F^n-\cab{\vit}=F^n$ is $F^n$, thus $g^*=F^n$ and $\vs=\vit$, and $\frac{\vs-V^{n-1}}{\tau}=g^*=F^n$.

However, if $\vit\geq \beta$, we still seek $\vs$ in $[a,b]$. Now any choice $\vs<b$ has the selection  $\cab{\vs} =\{0\}$, again with $g^*=F^n$, which leads to inconsistency in that $\tfrac{\vs -V^{n-1}}{\tau}\neq g^*$. Therefore we must have $\vs=b$ and the element from $g^*$ must be chosen exactly to make $\vs=\beta$, i.e. it equals   $g^*=\tfrac{\beta -V^{n-1}}{\tau}$. The corresponding selection $c^*$ out of $\cab{\vs}$ is $c^*=F^n-g^*=F^n-\tfrac{b-V^{n-1}}{\tau}$. 

A similar reasoning applies when $F^n<0$. The case $F^n=0$ is the intersection of the two cases. Summarizing and rearranging we get 
\bas
\vs=(I+\tau c)^{-1}\big(\tau f+V\big)
\mathrm{\ iff\ }
\vs = 
\begin{cases}
a, & \tau f+V \leq a, \\
\tau f+V,& a\leq \tau f+V \leq b,\\
b, & b \leq \tau f+V,
\end{cases}
\eas
which is the same as that given by \eqref{eq:resolvent}
upon $f=f^n,V=V^{n-1}$.  

\myplan{large F small f} 

}

%%%%%%%%%%%%%%%%%%%%%%
\myskip{
\subsection{b results and theory}
\label{sec:app-b}

%%%%%%%%%%%%%%%
When $\mu=1$, the output $w(t)$ in \eqref{eq:1generalized} is a solution to the doubly nonlinear ODE. 
\begin{equation}  \label{nonlinear-play}
w(t) \in b(v(t)), \quad
\tfrac{d}{d t} w(t)  + \calbe(v(t)) \ni 0\,,
\quad v(0) = \vz  \in u(0) + [\alpha,\beta]\,\,.
\end{equation}
We can estimate the solution of \eqref{nonlinear-play} with the following chain rule.
\begin{lemma}
Assume $w(\cdot) \in W^{1,1}(a,b;\R)$ and that $\sigma(t) \in \sgn(w(t))$ is a measurable selection. Then 
$$\tfrac{d}{dt} |w(t)| = \sigma(t) w'(t),\quad \text{a.e. } t \in [a,b].$$
\end{lemma}

\begin{proof}
From Appendix A of Kinderlehrer-Stampacchia (1980) we have $\dfrac{d}{dt} |w(t)| = \sgn_0(w(t)) w'(t)$ for a.e. $t \in [a,b]$. Let $S = \{t \in [a,b]: w(t) = 0\}$. If $t \notin S$ then $\sigma(t) =  \sgn_0(w(t))$. For the set of points in $S$, at the accumulation points we have $w(t) = 0$ and $w'(t) = 0$, and there are only a countable number of isolated points.
\end{proof}

\begin{lemma}		\label{b-unique}
Assume that $b(\cdot)$ is a maximal monotone function (necessarily single-valued and continuous). For a strong solution of \eqref{nonlinear-play}, $b(v(t))$ is uniquely determined by $b(v(0))$.
\end{lemma}

\begin{proof}
Set $w_j(t) = b(v_j(t))$ for $j = 1,2$ where
$$w_j'(t) + c(v_j(t) - u(t)) \ni 0.$$
Then $\sgn_0(v_1(t) - v_2(t)) \in \sgn(w_1(t) - w_2(t))$ since $b(\cdot)$ is a function. The function $w(\cdot)$ is a.e. differentiable and satisfies
\begin{eqnarray*}
\tfrac{d}{dt} |w_1(t) - w_2(t)| = (w_1'(t) - w_2'(t)) \sgn_0(v_1(t) - v_2(t)) 
\\
\in -\big(c_{\alpha,\beta}(v_1(t)-u(t))-c_{\alpha,\beta}(v_2(t)-u(t))\big) \sgn_0(v_1(t) - v_2(t)) \subset (-\infty,0],
\end{eqnarray*}
so we have $|w_1(t) - w_2(t)| \le |w_1(0)-w(0)|$ for $t \ge 0$.
\end{proof}

If $b(\cdot)$ is not injective then $v(\cdot)$ is not determined by \eqref{nonlinear-play}. \myskip{
In particular, if $b(\cdot)$ is constant on any interval within the constraint interval, then examples are easy to construct. We shall use a special solution of \eqref{nonlinear-play} that is constructed from the linear play. 
}

\begin{lemma}		\label{bv_strong}
Assume $b(\cdot)$ is a Lipschitz and monotone function. Define $w(t) = b(v(t))$, where $v(\cdot)$ is the strong solution of \eqref{linear-play}. Then $b(v(\cdot))$ is a strong solution of \eqref{nonlinear-play}. 
\end{lemma}
\begin{proof}
Since $b(\cdot)$ is Lipschitz, $w(\cdot)$ is differentiable a.e., and the chain rule gives
\begin{eqnarray*}
w'(t) = b'(v(t))v'(t) \in - b'(v(t))\,c(v(t) - u(t)) \subset - c(v(t) - u(t)),
\end{eqnarray*}
where the last relation follows from $b'(\cdot) \ge 0$.
\end{proof}
}
%%%%%%%%%%

\subsection{Iterative algorithm for parametrization of \knonlinear}
\label{sec:iteration}
We initialize iteration by setting $\hH_i^{(0)}=\hH_i$; these may have curvilinear sides. In each iteration $q=1,2\ldots q^*$, we parametrize each $\hH_i^{(q)}$ 
\ba
\label{eq:viq}
\ver(H_i^{(q)}) = \{(\alpha_{i}^{(q)},w_{i-1}), (\beta_{i}^{(q)},w_{i-1}),  (B_i^{(q)},w_i), (A_i^{(q)},w_i)\},
\ea
and find some $\pP^{(q)}_i$ with algorithm from Sec.~\ref{sec:trapezoid}.
We require continuity of the piecewise linear sides of $\hH_i^{(q)}$, i.e., %%
\ba
\label{eq:coincideq}
A^{(q)}_i=\alpha_{i+1}^{(q)}, \;\;  B^{(q)}_i=\beta_{i+1}^{(q)}, \;\; 0\leq i \leq I-1.
\ea
%%.  

%
Proceeding from $q$-$1 \mapsto q$ involves successive improvements of efficiency, while making sure that \eqref{eq:coincideq} holds. The process $q$-$1 \mapsto q$  is not automatic, leaving a lot of flexibility to adapt the requirements for accuracy and efficiency to their project's needs. For example, in this paper we choose to loop $i=1,2,\ldots I$ from ``bottom'' to the ``top'' of $\hH$, but other strategies are possible. 

\medskip
{\bf Loop $q=1,\ldots q^*$.}

\medskip
{\bf Loop $i=1,\ldots I$.}

\medskip
{\bf Step i.A.} Choose  $\alpha_i^{(q)},\beta_i^{(q)}$. 

If $i=1$, use $\alpha_i^{(q)}=\alpha_i^{(0)}=\alpha_i$, $\beta_i^{(q)}=\beta_i^{(0)}=\beta_i$. 

If $i>1$, use $\alpha_i^{(q)}=A_{i-1}^{(q)}$, $\beta_i^{(q)}=B_{i-1}^{(q)}$. 
(This ensures \eqref{eq:coincideq}). 

{\bf Step i.B.} Determine $A_i^{(q)},B_i^{(q)}$ with algorithm described in Sec.~\ref{sec:trapezoid} applied to the trial trapezoid $H_i^{(q-\tfrac{1}{2})}$ with the vertices 
\bas
\ver(H_i^{(q-\tfrac{1}{2})}) = 
\{ 
(\alpha_i^{(q)},w_{i-1}), 
(\alpha_i^{(q)},w_{i-1}), 
(A_i^{(q-\tfrac{1}{2})},w_{i}), 
(B_i^{(q-\tfrac{1}{2})},w_{i}).
\}
\eas
Here we select $A_i^{(q-\tfrac{1}{2})}$ from the set $A^{(q)}=\{A_i^{(0)}, A_i^{(1)}, \ldots A_i^{(q-1)}\}$
so that the sides of $H_i^{(q)}$ approximate best the curves  $(u,\fl(u))$, $\alpha_i^{(q)}\leq u\leq \max{A^{(q)}}$. Proceed similarly for $B_i^{(q-\tfrac{1}{2})}$. 

Next find the parametrization $\pP_i^{(q)}$ and two new vertices $A_i^{(q)}, B_i^{(q)}$
as in Sec.~\ref{sec:trapezoid}. Adjust these as needed for accuracy and efficiency.

{\bf Step i.C.} Continue to next $i$.

{\bf End loop over $i=1,\ldots I$.}

\medskip
{\bf Check the quality of current approximation $\hH^{(q)} \approx \hH$}: Assess the fit of the sides of $H_i^{(q)}$ to the curves $\fr(\cdot),\fl(\cdot)$, and whether the overall number $K=\sum_{i=1}^I K_i^{(q)}$ is small enough to be practical. 

If not, set $q=q+1$ and continue to the next iteration. If yes, we're done. 

{\bf End loop over $q$}. Set $q^*=q$, and $\hH_i^*=\hH_i^{(q)}$. Set $\pP^*=(\pP^{(q)}_i)_{i=i}^I$. 
%%

%%%%%%%%%%%%%%%%%%%%%
\bibliographystyle{plain}
\section*{References}

\bibliography{referencesPS20,referencesAHF}

\end{document}